\documentclass[reqno,11pt,a4paper]{amsart}
\usepackage[a4paper,left=30mm,right=30mm,top=30mm,bottom=30mm,marginpar=20mm]{geometry}
\usepackage{amsmath}
\usepackage{amssymb}
\usepackage{amsthm}
\usepackage{mathrsfs} 
\usepackage{setspace}
\usepackage{xcolor}
\usepackage{pict2e}
\usepackage{graphicx}

\usepackage{enumitem}

\newlist{property}{enumerate}{1}
\setlist[property]{
    label=(\textasteriskcentered),
    leftmargin=*,
    resume
}

\usepackage{cite}
\usepackage[colorlinks=true, linkcolor=blue, urlcolor=blue, citecolor=green]{hyperref}

\usepackage[english=american]{csquotes}

\newcommand{\wkto}{\mathrel{\smash{{\overset{\raisebox{-0.2ex}{$\scriptstyle\ast$}}{\rightharpoonup}}}}}

\newcommand{\supext}[1]{\scalebox{1.2}{$\wedge$}^{\!#1}\mkern2mu}
\newcommand{\subext}[1]{\scalebox{1.2}{$\wedge$}_{#1}\mkern2mu}

\newcommand{\mvert}{\mathrel{
		\begingroup
		\begin{picture}(7,7)
			\roundcap\roundjoin
			\linethickness{0.6pt}
			\polyline(0,7)(0,0)(7,0)
		\end{picture}
		\endgroup
}}

\newcommand{\Crm}{\mathrm{C}}

\newcommand{\Irm}{\mathrm{I}}

\newcommand{\Lrm}{\mathrm{L}}

\newcommand{\Nrm}{\mathrm{N}}

\newcommand{\Prm}{\mathrm{P}}
\newcommand{\Qrm}{\mathrm{Q}}

\newcommand{\Srm}{\mathrm{S}}

\newcommand{\IPrm}{\mathrm{IP}}
\newcommand{\QPrm}{\mathrm{QP}}
\newcommand{\IRrm}{\mathrm{IR}}

\newcommand{\Dcal}{\mathcal{D}}

\newcommand{\Mcal}{\mathcal{M}}

\newcommand{\Fbf}{\mathbf{F}}

\newcommand{\Mbf}{\mathbf{M}}
\newcommand{\Nbf}{\mathbf{N}}

\newcommand{\Xbf}{\mathbf{X}}

\renewcommand{\Bbb}{\mathbb{B}}

\newcommand{\Qbb}{\mathbb{Q}}

\newcommand{\Sbb}{\mathbb{S}}

\DeclareMathOperator*{\esssup}{ess\,sup}

\DeclareMathOperator{\supp}{supp}
\DeclareMathOperator{\diam}{diam}

\newcommand{\dd}{\;\mathrm{d}}

\newcommand{\N}{\mathbb{N}}
\newcommand{\R}{\mathbb{R}}

\newcommand{\Q}{\mathbb{Q}}
\newcommand{\Z}{\mathbb{Z}}

\newcommand{\FN}{\mathbf{F}}
\newcommand{\FI}{\mathbf{FI}}

\newcommand{\Var}{\mathrm{Var}}

\newcommand{\NLip}{\Nrm^{\mathrm{Lip}}}
\newcommand{\ILip}{\Irm^{\mathrm{Lip}}}
\newcommand{\PLip}{\Prm^{\mathrm{Lip}}}
\newcommand{\IPLip}{\IPrm^{\mathrm{Lip}}}
\newcommand{\Lip}{\mathrm{Lip}}

\newcommand{\Int}{\mathrm{int}\hspace{0.2em}}
\renewcommand{\det}{\mathrm{det}}
\newcommand{\tbf}{\mathbf{t}}
\newcommand{\pbf}{\mathbf{p}}

\newtheorem{definition}{Definition}[section]
\newtheorem{proposition}[definition]{Proposition}
\newtheorem{lemma}[definition]{Lemma}
\newtheorem{theorem}[definition]{Theorem}
\newtheorem{corollary}[definition]{Corollary}

\newtheorem*{theorem*}{Theorem}
\newtheorem*{proposition*}{Proposition}
\numberwithin{equation}{section}

\theoremstyle{definition}
\newtheorem{example}[definition]{Example}
\newtheorem{remark}[definition]{Remark}

\begin{document}
	
\title[Strict Approximation of Currents and Space-Time Trajectories]{On the Strict Approximation of Currents and Space-Time Trajectories}
	
\begin{abstract} This work establishes a wide variety of strict approximation results for space-time trajectories, a class of space-time currents that are used to model the evolution of dislocations in crystalline materials. In particular, we introduce a strict-type notion of convergence which respects quantities associated with geometric transport phenomena and prove approximation results in this topology. These results build on the classical strict approximation results for normal and integral currents developed by Federer and Fleming, which we also review and refine.
    	
\vspace{.5em}
\noindent\textsc{Date:} \today{}
\end{abstract}
	
\author{Paolo Bonicatto}
\address{Universit\`a di Trento, Dipartimento di Matematica, Via Sommarive 14, 38123 Trento, Italy}
\email{paolo.bonicatto@unitn.it}
	
\author{Filip Rindler}
\address{Mathematics Institute, University of Warwick, Coventry, CV4 7AL, United Kingdom}
\email{F.Rindler@warwick.ac.uk}
	
\author{Harry Turnbull}
\address{Mathematics Institute, University of Warwick, Coventry, CV4 7AL, United Kingdom}
\email{harry.turnbull@warwick.ac.uk}

\maketitle

\section{Introduction}

Approximation results have occupied a central place in the development of the theory of currents over the past six decades. Already in their seminal paper \cite{FedererFleming60}, Federer and Fleming proved three results concerning the approximation of normal and integral currents by polyhedral chains; since then, a wide range of results have been obtained by varying the class of currents to be approximated, the class of approximants, and the topology of the approximation (see, for example, \cite{ChambolleFerrariMerlet21,AlmgrenBrowderCaldiniDeLellis24?,White99,GiaquintaModicaSoucek98book1,GiaquintaModicaSoucek98book2}). 

More recently, a theory of space-time currents was introduced in the integral setting in \cite{Rindler23} and extended to normal currents in \cite{BonicattoRindlerTurnbull26}. These developments were motivated by the model of \cite{HudsonRindler22,Rindler26?} for dislocation motion in crystalline solids \cite{HullBacon11book,AndersonHirthLothe17book}. This theory, which describes the transport of geometric structures represented by currents via their `space-time track' traced out as they move, was further developed and used in \cite{Rindler25,BonicattoDelNinRindler24,BonicattoDelNinRindler25,Bonicatto25?,BonicattoDelNin25?}. 

In the present paper we shall develop several approximation results for a particular class of space-time currents which respect their inherent space-time structure. Together, these approximation results form the technical backbone of the homogenisation procedure for dislocation-driven plasticity conducted in \cite{BonicattoRindlerTurnbull26-II?} (which in turn rigorously implements the model of \cite{Rindler26?}). 

In what follows, we fix $d \in \N$ and take $k \in \N_{0}$ with $k \leq d$. Although this paper is largely self-contained, we follow the notational conventions of \cite{BonicattoRindlerTurnbull26} and defer to \cite[Chapter 4]{Federer69book} for a comprehensive treatment of currents. We say that a sequence of $k$-currents $(T_{j})_{j} \subset \Dcal_{k}(\R^{d})$ with $\Mbf(T_{j}) < \infty$ for $j \in \N$ converges \textbf{strictly} to a $k$-current $T \in \Dcal_{k}(\R^{d})$ if $T_{j} \wkto T$ and $\Mbf(T_{j}) \to \Mbf(T)$. We shall restrict our attention to strict-type topologies where we replace the convergence of the mass in the definition above with the convergence of various functionals.

The first strict approximation results for currents appeared in the aforementioned paper \cite{FedererFleming60} of Federer and Fleming; there (see also \cite[4.2.21, 4.2.24]{Federer69book}), they demonstrated that, if $T \in \Nrm_{k}(\R^{d})$ (respectively $T \in \Irm_{k}(\R^{d})$) and $\supp T \subset \Int K$ for some compact $K \subset \R^{d}$, then one can find a sequence $(P_{j})_{j} \subset \Prm_{k}(\R^{d})$ (respectively $(P_{j})_{j} \subset \IPrm_{k}(\R^{d})$) of (integral) polyhedral $k$-chains with $\supp P_{j} \subset K$ for $j \in \N$ such that $P_{j} \to T$ strictly and $\partial P_{j} \to \partial T$ strictly.

To prove an analogous result in the space-time setting, we shall need to use many techniques from the proofs of \cite[4.2.21, 4.2.24]{Federer69book}. So that we can appeal to these techniques in the sequel, we present a complete and direct proof of \cite[4.2.21, 4.2.24]{Federer69book} which differs in several places from \cite{Federer69book}, but uses only classical tools. In fact, we obtain the following refinement of \cite[4.2.21, 4.2.24]{Federer69book} which allows us to additionally control the boundary of the approximants in the following way:

\begin{theorem}\label{tm:AbbrvPolyBdryStrict} Let $K_{1},K_{2} \subset \R^{d}$ be compact with $K_{2} \subset K_{1}$ and take $\varepsilon > 0$. Let $T \in \Nrm_{k}(\R^{d})$. Suppose that $\supp T \subset \Int K_{1}$ and that $\supp \partial T \subset \Int K_{2}$. Then there is $P \in \Prm_{k}(\R^{d})$ with $\supp P \subset K_{1}$, $\supp \partial P \subset K_{2}$, and
\begin{align*}
    \FN_{K_{1}}(P-T) &\leq \varepsilon, & \Mbf(P) &\leq \Mbf(T) + \varepsilon, & \Mbf(\partial P) &\leq \Mbf(\partial T) + \varepsilon. 
\end{align*}
If $T \in \Irm_{k}(\R^{d})$, then we can additionally ensure that $P \in \IPrm_{k}(\R^{d})$.
\end{theorem}

Here, for a compact set $K \subset \R^{d}$ and a normal current $T \in \Nrm_{k}(K)$, we recall that the \textbf{(normal) flat norm} of $T$ (in $K$) is given by
\begin{align}\label{NormalFlatNormDefn}
    \FN_{K}(T) := \inf \bigl\{\Mbf(T-\partial R) + \Mbf(R): R \in \Nrm_{k+1}(K) \bigr\}.
\end{align}

We shall now introduce the class of space-time currents that we shall approximate in the sequel. Recall from \cite[Section 4]{BonicattoRindlerTurnbull26} that a \textbf{space-time trajectory} is a normal current $S \in \Nrm_{1+k}([0,1] \times \R^{d})$ with $\partial S \mvert (0,1) \times \R^{d} = 0$. Notice that we can write
\begin{align*}
    \partial S = \delta_{1} \times \partial^{+}S - \delta_{0} \times \partial^{-}S
\end{align*}
for some $\partial^{\pm}S \in \Nrm_{k}(\R^{d})$ with $\partial(\partial^{\pm}S) = 0$, and so one can interpret $S$ as a ``trajectory'' describing the evolution of $\partial^{-}S$ into $\partial^{+}S$. Intermediate states can be recovered by slicing with respect to the temporal projection $\tbf(t,x):=t$ and projecting these slices onto the spatial component via $\pbf(t,x) := x$; in particular, we define
\begin{align*}
    S(t) := \pbf_*(S|_t) \in \Nrm_{k}(\R^{d}) \hspace{1em}\text{for $\mathscr{L}^{1}$-a.e.\ $t \in [0,1]$}
\end{align*}
where $\pbf_*$ denotes the pushforward under $\pbf$ defined in \cite[4.1.7]{Federer69book}. From the perspective of geometric transport phenomena, the relevant quantities associated to a space-time trajectory $S$ are $\Mbf(S(t))$ for $\mathscr{L}^{1}$-a.e.\ $t \in [0,1]$ and the so-called \textbf{variation} of $S$, which measures the ``area swept out'' by $S$ (with ``absolute multiplicity'') and was defined in \cite{Rindler23,BonicattoRindlerTurnbull26} as follows:
\begin{align*}
    \Var(S) := \int \Vert \pbf \vec{S} \Vert \dd \Vert S \Vert.
\end{align*}
With this in mind, we introduce the notion of \textit{space-time strict convergence} in Section \ref{s:STStrict}. This is the convergence with respect to which we shall establish approximation results for space-time trajectories. 

In general, there is no relationship between the space-time and classical strict convergences (see Section \ref{s:STStrict}) and so one cannot hope to modify the proof of Theorem \ref{tm:AbbrvPolyBdryStrict} to obtain an analogous space-time strict polyhedral approximation result for normal currents. Furthermore, it is often necessary in applications (see \cite{BonicattoRindlerTurnbull26-II?}) to keep track of the so-called \textbf{$\Lrm^{\infty}$-norm}, which is defined for a space-time trajectory $S$ as follows:
\begin{align*}
    \Vert S \Vert_{\Lrm^{\infty}} := \esssup_{t \in [0,1]} \Mbf(S(t)).
\end{align*}
As a uniform quantity, the $\Lrm^{\infty}$-norm does not behave well under the classical approximation procedures that are reviewed in Section \ref{s:ClassicalPolyStrict}. Nevertheless, we obtain the following space-time strict polyhedral approximation result, which is an abbreviated form of the full statement in Theorem \ref{tm:STPoly}:

\begin{theorem}\label{tm:AbbrvSTPoly} Let $K \subset \R^{d}$ be compact and take $\varepsilon > 0$. Let $S \in \Nrm_{1+k}([0,1] \times \R^{d})$. Suppose that $\partial S \mvert (0,1) \times \R^{d} = 0$ and that $\supp S \subset [0,1] \times \Int K$. Then there is $P \in \Prm_{1+k}([0,1] \times \R^{d})$ with $\partial P \mvert (0,1) \times \R^{d} = 0$, $\supp P \subset [0,1] \times K$, 
\begin{align*}
    \Fbf_{[0,1] \times K}(P-S) &\leq \varepsilon, 
\end{align*}
and
\begin{align*}
    \Var(P) &\leq \Var(S)+\varepsilon, & \int \Mbf(P(t)) \dd t &\leq \int \Mbf(S(t)) \dd t + \varepsilon.
\end{align*}
Furthermore,
\begin{align*}
    \Mbf(\partial^{\pm}P) &\leq \Mbf(\partial^{\pm}S) + \varepsilon, & \Vert {P} \Vert_{\Lrm^{\infty}} &\leq \Vert S \Vert_{\Lrm^{\infty}} + \varepsilon.
\end{align*}
\end{theorem}

We then turn our attention to different classes of approximants. We say that a normal current $T \in \Nrm_{k}(\R^{d})$ is \textbf{rational} if there is $n \in \N$ with $n \cdot T \in \Irm_{k}(\R^{d})$, and we use $\Qrm_{k}(\R^{d})$ to denote the space of rational currents. Notice that a polyhedral chain is a rational current if and only if it can be expressed as a linear combination of oriented simplices with rational coefficients, so one can easily adapt Theorem \ref{tm:AbbrvPolyBdryStrict} and Theorem \ref{tm:AbbrvSTPoly} to obtain strict approximation results in this setting by simply adjusting coefficients. However, doing so completely destroys any control on the support of the boundary of the approximants. In Section \ref{s:Q}, we rectify this with a more precise construction that takes inspiration from \cite[5.5]{Federer74} and \cite[Section 5.4.1]{GiaquintaModicaSoucek98book1}. It turns out that this construction preserves any rational pieces of the boundary and, in the context of space-time trajectories, we can use this observation to prove the following result, which is an abbreviated form of the full statement in Theorem \ref{tm:STOneSidedQ}:

\begin{theorem}\label{tm:AbbrvSTOneSidedQ} Let $K \subset \R^{d}$ be compact and take $\varepsilon > 0$. Let $S \in \Nrm_{1+k}([0,1] \times \R^{d})$. Suppose that $\partial S \mvert (0,1) \times \R^{d} = 0$, $\supp S \subset [0,1] \times \Int K$, and that $\partial^{-}S \in \Qrm_{k}(\R^{d})$. Then there is $Q \in \Qrm_{1+k}([0,1] \times \R^{d})$ with $\partial Q \mvert (0,1) \times \R^{d} = 0$, $\partial^{-}Q = \partial^{-}S$, $\supp Q \subset [0,1] \times K$,
\begin{align*}
    \Fbf_{[0,1] \times K}(Q-S) &\leq \varepsilon,
\end{align*}
and
\begin{align*}
    \Var(Q) &\leq \Var(S)+\varepsilon, & \int \Mbf(Q(t)) \dd t &\leq \int \Mbf(S(t)) \dd t + \varepsilon.
\end{align*}
Furthermore,
\begin{align*}
    \Mbf(\partial^{+}Q) &\leq \Mbf(\partial^{+}S) + \varepsilon, & \Vert {Q} \Vert_{\Lrm^{\infty}} &\leq \Vert S \Vert_{\Lrm^{\infty}} + \varepsilon.
\end{align*}
\end{theorem}

The final class of approximants that we consider are the so-called \textit{Lipschitz trajectories}. Roughly speaking, a space-time trajectory is said to be Lipschitz if it possesses `Lipschitz regularity in time' -- this notion is made precise in \cite[Section 4.1]{BonicattoRindlerTurnbull26}. In Section \ref{s:Lipschitz}, we show that space-time trajectories can be approximated by Lipschitz trajectories with the same endpoints:

\begin{theorem}\label{tm:AbbrvSTLipschitz} Let $K \subset \R^{d}$ be compact and take $\varepsilon > 0$. Let $S \in \Nrm_{1+k}([0,1] \times \R^{d})$. Suppose that $\partial S \mvert (0,1) \times \R^{d} = 0$ and that $\supp S \subset [0,1] \times \Int K$. Then there is $V \in \NLip_{1+k}([0,1] \times \R^{d})$ with $\partial V = \partial S$, $\supp V \subset [0,1] \times K$,
\begin{align*}
    \Fbf_{[0,1] \times K}(V-S) &\leq \varepsilon,
\end{align*}
and
\begin{align*}
    \Var(V) &\leq \Var(S)+\varepsilon, &\int \Mbf(V(t)) \dd t &\leq \int \Mbf(S(t)) \dd t + \varepsilon.
\end{align*}
\end{theorem}
The full statement of Theorem \ref{tm:AbbrvSTLipschitz} is given in Theorem \ref{tm:STLipschitz}. Notice that, in this result, we do not control the $\Lrm^{\infty}$-norm of the approximants. 

Each of the approximation results discussed above shares the property that the support of the approximant grows in some way. When using approximation results involving weaker topologies, the classical approach to rectify this expansion of the support is to assume that the currents to be approximated are supported on some \emph{Lipschitz neighbourhood retract} $K \subset \R^{d}$. Then there is an open set $U \subset \R^{d}$ with $K \subset U$ and a Lipschitz retraction $r: U \to K$, and so one can push approximants back into $K$ using $r$. However, in \cite[Section 5]{BonicattoRindlerTurnbull26}, we gave an example of a Lipschitz neighbourhood retract $K \subset \R^{2}$ with the property that any retract $r$ onto $K$ satisfies $\Lip(r) \geq \sqrt{2}$. As a result, applying the procedure outlined above would yield approximants that are not necessarily `strict-close' to the original current in general. In Section \ref{s:Support}, we use a different approach to control the support of approximants on a class of compact sets which include the closures of $\Crm^{0}$ domains (so, in particular, Lipschitz domains).

\subsection*{Acknowledgements} The authors would like to thank Luca Gennaioli for discussions related to this work. This work was supported by the Italian Ministry of University and Research (MUR) through the FIS2 project BRIDGE: ‘‘Beyond Regularity In Differential Geometric Evolutions'', project code FIS-2023-02224 (CUP: E53C25001790001). This project has received funding from the UKRI Frontier Research Guarantee (ERC guarantee) grant EP/Z000297/1 (ERC CONCENTRATE). H.T.\ was supported by the Warwick Mathematics Institute Centre for Doctoral Training and the EPSRC grant EP/W524645/1. The authors have no relevant financial or non-financial interests to disclose. 

\subsection*{AI Declaration} 

AI was used in the preparation of this paper to discuss a few minor ideas and for proofreading, correcting typos, and literature search. This manuscript contains no AI-generated text.

\section{Preliminaries}

Recall from the introduction that we refer to \cite[Chapter 4]{Federer69book} for a comprehensive treatment of currents and that we follow the notational conventions of\cite{BonicattoRindlerTurnbull26}. In particular, we identify $\R^{1+d} \cong \R \times \R^{d}$ and define the orthogonal temporal and spatial projections as $\tbf: \R^{1+d} \to \R$ and $\pbf: \R^{1+d} \to \R^{d}$ where $\tbf(t,x) = t$ and $\pbf(t,x)=x$ for $(t,x) \in \R^{1+d}$ respectively. Furthermore, if $k=0$, then any statements involving $\partial T$ are to be ignored. 

\subsection{Inhomogeneous and Anisotropic Functionals}\label{s:Inhomogeneous} The results discussed in the introduction were written using the mass and variation, which are both \textit{homogeneous} and \textit{isotropic} functionals. In the sequel, we shall prove more general results which establish control of variants of the mass and variation which are both \textit{inhomogeneous} and \textit{anisotropic} (as required, for example, in \cite{BonicattoRindlerTurnbull26-II?, Rindler26?}).

Let $\psi: \R^{d} \times \subext{k}\R^{d} \to [0,\infty)$ be continuous and suppose that $\psi(x,\cdot)$ is positively $1$-homogeneous for $x \in \R^{d}$. The \textbf{$\psi$-mass} of a $k$-current $T \in \Dcal_{k}(\R^{d})$ admitting a Radon--Nikodym decomposition $T = \vec{T} \, \Vert T \Vert$ is given by
\begin{align*}
    \Mbf_{\psi}(T) := \int \psi(x,\vec{T}(x)) \dd \Vert T \Vert(x).
\end{align*}
Following \cite[5.1.1]{Federer69book} and \cite[Section 2]{Federer74}, we say that a continuous map $\psi: \R^{d} \times \subext{k}\R^{d} \to [0,\infty)$ is a \textbf{positive convex integrand} if $\psi(x,\cdot)$ is a norm on $\subext{k}\R^{d}$ for $x \in \R^{d}$.

Let $\psi: \R^{d} \times \subext{k}\R^{d} \to [0,\infty)$ be a positive convex integrand and take a sequence $(T_{j})_{j} \subset \Dcal_{k}(\R^{d})$ with $\Mbf(T_{j}) < \infty$ for $j \in \N$. We say that $T_{j}$ converges \textbf{$\psi$-strictly} to $T \in \Dcal_{k}(\R^{d})$ if 
\begin{align*}
    T_{j} &\wkto T, & \Mbf_{\psi}(T_{j}) \to \Mbf_{\psi}(T).
\end{align*}
Although the notion of $\psi$-mass is classical (see, for example, \cite[Section 2]{Federer74}), we shall briefly discuss the relationship between the strict and the $\psi$-strict convergence. To this end, let $K \subset \R^{d}$ be compact and notice that
there is a constant $C_{\psi,K} > 0$ such that, for $(x,\xi) \in K \times \subext{k}\R^{d}$,
\begin{align}\label{psiVsMass}
    C_{\psi,K}^{-1} \cdot \Vert \xi \Vert \leq \psi(x,\xi) \leq C_{\psi,K} \cdot \Vert \xi \Vert.
\end{align}
Here, $\Vert \xi \Vert$ denotes the mass norm of $\xi$ (see \cite[1.8.1]{Federer69book}). It follows that, for $T \in \Dcal_{k}(\R^{d})$ with $\Mbf(T) < \infty$ and $\supp T \subset K$, we have
\begin{align}\label{psiMassVsMass}
    C_{\psi,K}^{-1} \cdot \Mbf(T) \leq \Mbf_{\psi}(T) \leq C_{\psi,K} \cdot \Mbf(T).
\end{align}
Take a sequence $(T_{j})_{j} \subset \Nrm_{k}(K)$, let $T \in \Nrm_{k}(K)$, and suppose that 
\begin{align*}
    \Fbf_{K}(T_{j}-T) &\to 0, & \Mbf(T_{j}) &\to \Mbf(T).
\end{align*}
One can use inequalities analogous to \eqref{psiMassVsMass} to show that there is no benefit to considering inhomogeneous and anisotropic variants of the flat norm in this setting. However, one cannot use \eqref{psiMassVsMass} to deduce that $\Mbf_{\psi}(T_{j}) \to \Mbf_{\psi}(T)$. Instead, one may hope to use Reshetnyak's continuity theorem (see Appendix \ref{s:Reshetnyak}) to deduce $\psi$-strict convergence from the strict convergence of $(T_{j})_{j}$, but this approach fails because $(\subext{k}\R^{d},\Vert \cdot \Vert)$ is not strictly convex in general.

On the other hand, $(\subext{k}\R^{d}, \vert \cdot \vert)$ is a Hilbert space; this suggests that we can use the approach above when we consider classes of currents where we can instead work with the Euclidean norm $\vert \cdot \vert$ (see \cite[1.7.5]{Federer69book}) on $k$-vectors. Take $T \in \Dcal_{k}(\R^{d})$ admitting a Radon--Nikodym decomposition $T = \vec{T} \, \Vert T \Vert$. We say that $T$ is \textbf{simple} if $\vec{T}(x)$ is simple for $\Vert T \Vert$-a.e.\ $x \in \R^{d}$. Notice that, if $T$ is simple, then $\vert \vec{T}(x) \vert = \Vert \vec{T}(x) \Vert = 1$ for $\Vert T \Vert$-a.e.\ $x \in \R^{d}$ and so
\begin{align*}
    T = \vec{T} \hspace{0.2em} \Vert T \Vert
\end{align*}
is the polar decomposition of $T$ with respect to both $\vert \cdot \vert$ and $\Vert \cdot \Vert$ on $\subext{k}\R^{d}$. In particular, on the class of simple currents, the notion of strict convergence does not change if we work with $\vert \cdot \vert$ in place of $\Vert \cdot \Vert$. We therefore have

\begin{proposition}\label{pp:Simple} Let $\psi: \R^{d} \times \subext{k}\R^{d} \to [0,\infty)$ be continuous and suppose that $\psi(x,\cdot)$ is positively $1$-homogeneous for $x \in \R^{d}$. Let $K \subset \R^{d}$ be compact. Take $(T_{j})_{j} \subset \Dcal_{k}(\R^{d})$ with $\Mbf(T_{j}) < \infty$ and $\supp T_{j} \subset K$ for $j \in \N$ and take $T \in \Dcal_{k}(\R^{d})$. Suppose that $T_{j}$ is simple for $j \in \N$ and that $T$ is simple. If $T_{j} \to T$ strictly, then
\begin{align}\label{cc:Simple}
    \Mbf_{\psi}(T_{j}) \to \Mbf_{\psi}(T).
\end{align}
\end{proposition}

\begin{proof} Since $(\subext{k}\R^{d},\vert \cdot \vert)$ is strictly convex, this result is an immediate consequence of Proposition \ref{pp:ReshetnyakCty}.
\end{proof}

We cannot expect that Proposition \ref{pp:Simple} holds without the assumption that the sequence and its limit are simple. To illustrate this, let $\psi: \R^{d} \times \subext{k}\R^{d} \to [0,\infty)$ be a positive convex integrand; one can use Proposition \ref{pp:ReshetnyakSharp} in conjunction with Remark \ref{rm:Reshetnyak} to deduce that, if \eqref{cc:Simple} holds for all strictly converging sequences of normal currents in $[0,1]^{d}$, then $\psi(x,\cdot)$ is affine on the proper faces of the unit ball of $(\subext{k}\R^{d},\Vert \cdot \Vert)$ for all $x \in [0,1]^{d}$. Unfortunately, this assumption on $\psi$ is too strong for our purposes; in particular, the integrands used in the proof of Proposition \ref{pp:SimpleSTS} do not satisfy this condition.  In Section \ref{s:NormalPolyStrict}, we shall therefore
construct polyhedral approximants of normal currents where the $\psi$-mass is directly controlled. We remark that this procedure requires an extra condition on $\psi$, but it is straightforward to see that this condition is optimal (see Section \ref{s:NormalPolyStrict} for more details). 

Let $K \subset \R^{d}$ be compact, take $(T_{j})_{j} \subset \Dcal_{k}(\R^{d})$ with $\Mbf(T_{j}) < \infty$ and $\supp T_{j} \subset K$ for $j \in \N$, and let $T \in \Dcal_{k}(\R^{d})$. Suppose that 
\begin{align}\label{MeasureWkto}
    T_{j}(\omega) \to T(\omega) \hspace{1em} \text{for all $\omega \in \Crm(K;\supext{k}\R^{d})$}.
\end{align}
It follows from Reshetnyak's lower semicontinuity theorem (which, unlike the continuity theorem, holds for measures taking values in \textit{any} finite-dimensional normed space) that 
\begin{align}\label{psiLSC}
    \Mbf_{\psi}(T) \leq \liminf_{j \to \infty}\Mbf_{\psi}(T_{j}).
\end{align}
In fact, this holds for any continuous $\psi: \R^{d} \times \subext{k}\R^{d} \to [0,\infty)$ such that $\psi(x,\cdot)$ is positively $1$-homogeneous and convex for $x \in \R^{d}$ (see, for example, \cite[Theorem 2.38]{AmbrosioFuscoPallara00book}). It turns out that, since we additionally assume that $\psi(x,\cdot)$ is positive definite for $x \in \R^{d}$, we can relax the condition \eqref{MeasureWkto} and assume only that $T_{j} \wkto T$:

\begin{lemma}\label{lm:psiLSC} Let $\psi: \R^{d} \times \subext{k}\R^{d} \to [0,\infty)$ be a positive convex integrand and let $K \subset \R^{d}$ be compact. Take $(T_{j})_{j} \subset \Dcal_{k}(\R^{d})$ with $\Mbf(T_{j}) < \infty$ and $\supp T_{j} \subset K$ for $j \in \N$. If $T_{j} \wkto T$ for some $T \in \Dcal_{k}(\R^{d})$, then \eqref{psiLSC} holds. 
\end{lemma}

Lemma \ref{lm:psiLSC} is an immediate consequence of the following result:

\begin{lemma}\label{lm:psiDual} Let $\psi: \R^{d} \times \subext{k}\R^{d} \to [0,\infty)$ be a positive convex integrand and let $K \subset \R^{d}$ be compact. Take $T \in \Dcal_{k}(\R^{d})$ with $\Mbf(T) < \infty$ and $\supp T \subset K$. Then
\begin{align*}
    \Mbf_{\psi}(T) = \sup \bigl\{T(\omega): \omega \in \Dcal^{k}(\R^{d}), \hspace{0.2em} \omega(x) \leq \psi(x,\cdot) \hspace{0.2em} \text{for $x \in \R^{d}$} \bigr\}.
\end{align*}
\end{lemma}

Let $\psi: \R^{d} \times \subext{k}\R^{d} \to [0,\infty)$ be a positive convex integrand and take $\omega \in \Dcal^{k}(\R^{d})$. Here, and in what follows, we identify $\supext{k}\R^{d} = (\subext{k}\R^{d})^{*}$ so that, for $x \in \R^{d}$, both $\omega(x)$ and $\psi(x,\cdot)$ are real-valued functions on $\subext{k}\R^{d}$ and the inequality $\omega(x) \leq \psi(x,\cdot)$ can be interpreted pointwise. In the case $\psi(x,\cdot) = \Vert \cdot \Vert$ for $x \in \R^{d}$, Lemma \ref{lm:psiDual} reduces to the familiar relation $\Vert T \Vert(\R^{d}) = \Mbf(T)$. 

Lemma \ref{lm:psiDual} is the content of \cite[2.1]{Federer74}, but for the benefit of the reader we give a proof with more details and in our notation below:

\begin{proof}[Proof of Lemma \ref{lm:psiDual}] First take $\omega \in \Dcal^{k}(\R^{d})$ with $\omega(x) \leq \psi(x,\cdot)$ for $x \in \R^{d}$. Then
\begin{align*}
    T(\omega) = \int \langle \vec{T}(x), \omega(x) \rangle \dd \Vert T \Vert(x) \leq \Mbf_{\psi}(T).
\end{align*}
It follows that
\begin{align*}
    \Mbf_{\psi}(T) \geq \sup \bigl\{T(\omega): \omega \in \Dcal^{k}(\R^{d}), \hspace{0.2em} \omega(x) \leq \psi(x,\cdot ) \hspace{0.5em} \text{for $x \in \R^{d}$} \bigr\}.
\end{align*}
It remains to establish the converse inequality. First notice that we can assume that $\supp T \subset \Int K$ without loss of generality. Fix $\varepsilon > 0$. We appeal to Lusin's theorem to find a compact set $A \subset \supp T$ such that
\begin{align}\label{psiDual1}
    \Vert T \Vert((\supp T) \backslash A) &\leq \varepsilon, & \vec{T} \vert_{A} \hspace{0.5em} &\text{is continuous}, & \Vert \vec{T}(x) \Vert &=1 \hspace{0.5em} \text{for $x \in A$}.
\end{align}
Since $A$ is compact, $\vec{T}$ is uniformly continuous on $A$. We can therefore find $\delta_{1} > 0$ such that, for $x_{1},x_{2} \in A$, 
\begin{align}\label{psiDual2}
    \vert x_{1} - x_{2} \vert \leq \delta_{1} \Rightarrow \Vert \vec{T}(x_{1}) - \vec{T}(x_{2}) \Vert \leq \varepsilon.
\end{align}
Furthermore, since $\psi$ is continuous and $\{\xi \in \subext{k}\R^{d}: \Vert \xi \Vert = 1\}$ is compact, we can find $\delta_{2} > 0$ such that, for $x_{1},x_{2} \in K$,
\begin{align}\label{psiDual3}
    \vert x_{1} - x_{2} \vert \leq \delta_{2} \Rightarrow \vert \psi(x_{1},\xi) - \psi(x_{2},\xi) \vert \leq \varepsilon \cdot \Vert \xi \Vert \hspace{1em} \text{for $\xi \in \subext{k}\R^{d}$}.
\end{align}
Now take $\delta < \min\{\delta_{1},\delta_{2},\mathrm{dist}(A,\R^{d}\backslash K)\}$. Since $A$ is compact, there is $N \in \N$ and $\{a_{j}\}_{j=1}^{N} \subset A$ such that 
\begin{align*}
     A \subset \bigcup_{j=1}^{N} B(a_{j},\delta).
\end{align*}
For $j \in \{1,\dotsc,N\}$, let $\varphi_{j}: \R^{d} \to [0,1]$ be smooth with $\supp\hspace{0.2em} \varphi_{j} \subset B(a_{j},\delta)$. We can additionally ensure that
\begin{align*}
    \sum_{j=1}^{N} \varphi_{j}(x) &= 1 \hspace{1em} \text{for $x \in A$}, & \sum_{j=1}^{N} \varphi_{j}(x) \leq 1 \hspace{1em} \text{for $x \in \R^{d}$}.
\end{align*}
Fix $j \in \{1,\dotsc,N\}$. We use the Hahn-Banach theorem to find $\omega_{j} \in \supext{k}\R^{d}$ with 
\begin{align}\label{psiDual4}
     \omega_{j} &\leq \psi(a_{j},\cdot), & \langle \vec{T}(a_{j}), \omega_{j} \rangle &= \psi(a_{j},\vec{T}(a_{j})).
\end{align}
Take $x \in \overline{B}(a_{j},\delta)$. It follows from \eqref{psiDual3} that, for $\xi \in \subext{k}\R^{d}$,
\begin{align}\label{psiDual5}
    \vert \psi(x,\xi) - \psi(a_{j},\xi) \vert \leq \varepsilon \cdot \Vert \xi \Vert.
\end{align}
If we additionally have $x \in A$, then we can use \eqref{psiDual2} and \eqref{psiDual5} to deduce that
\begin{align}\label{psiDual6}
\begin{aligned}
    \psi(x,\vec{T}(x)) &\leq \psi(a_{j},\vec{T}(x)) + \varepsilon \\
    &\leq \psi(a_{j},\vec{T}(a_{j})) + (C_{\psi,K}+1)\cdot \varepsilon \\
    &= \langle \vec{T}(a_{j}),\omega_{j} \rangle + (C_{\psi,K}+1) \cdot \varepsilon \\
    &\leq \langle \vec{T}(x),\omega_{j}\rangle + (2C_{\psi,K}+1)\cdot \varepsilon
\end{aligned}
\end{align}
where $C_{\psi,K} > 0$ is defined as in \eqref{psiVsMass}.
Now set
\begin{align*}
    \tilde{\omega} := \sum_{j=1}^{N} \varphi_{j}\omega_{j}.
\end{align*}
Then we can use \eqref{psiDual1} and \eqref{psiDual4} to deduce that
\begin{align*}
    \Mbf_{\psi}(T) &\leq \int_{A} \psi(x,\vec{T}(x)) \dd \Vert T \Vert(x) + C_{\psi,K} \cdot \varepsilon 
\end{align*}
and
\begin{align*}
    \int_{A} \langle \vec{T}(x),\tilde{\omega}(x)\rangle \dd \Vert T \Vert(x) \leq T(\tilde{\omega})+ C_{\psi,K} \cdot \varepsilon.
\end{align*}
Combining the inequalities above with \eqref{psiDual6} yields
\begin{align}\label{psiDual7}
    \Mbf_{\psi}(T)\leq T(\tilde{\omega}) + \varepsilon((2C_{\psi,K}+1) \Mbf(T)+2C_{\psi,K}).
\end{align}
Take $j \in \{1,\dotsc,N\}$ and $x \in B(a_{j},\delta)$. It follows from \eqref{psiDual4} and \eqref{psiDual5} that
\begin{align}\label{psiDual8}
    \langle \tilde{\omega}(x),\xi \rangle \leq (1+C_{\psi,K}\cdot\varepsilon)\psi(x,\xi) \hspace{1em} \text{for $\xi \in \subext{k}\R^{d}$}.
\end{align}
Now set
\begin{align*}
    \omega := \frac{1}{1+C_{\psi,K}\cdot \varepsilon}\cdot \tilde{\omega}.
\end{align*}
It follows from \eqref{psiDual8} that $\omega(x) \leq \psi(x,\cdot)$ for $x \in \R^{d}$. We can now use \eqref{psiDual7} to deduce that
\begin{align*}
     \Mbf_{\psi}(T) \leq (1+C_{\psi,K}\cdot \varepsilon) T(\omega) + \varepsilon \cdot ((2C_{\psi,K}+1)\Mbf(T)+ 2C_{\psi,K}).
\end{align*}
The result follows since the choice of $\varepsilon$ was arbitrary.
\end{proof}

Let $f: \R^{d} \to \R^{d}$ be Lipschitz and take $T \in \Nrm_{k}(\R^{d})$. Before proceeding, we recall the results of \cite[Section 5]{AlbertiMarchese16}, from which it follows that
\begin{align*}
    f_{*}T(\omega) = \int \langle D^{T}f(x)\vec{T}(x), \omega(f(x)) \rangle \dd \Vert T \Vert(x) \hspace{1em} \text{for $\omega \in \Dcal^{k}(\R^{d})$},
\end{align*}
where, for $\Vert T \Vert$-a.e.\ $x \in \R^{d}$, the linear map $D^{T}f(x): \mathrm{span}\hspace{0.2em}\vec{T}(x) \to \R^{d}$ is given by
\begin{align*}
    D^{T}f(x)v := \lim_{h \to 0}\frac{f(x+hv)-f(x)}{h} \hspace{1em} \text{for $v \in \mathrm{span}\hspace{0.2em}\vec{T}(x)$}.
\end{align*}

Let $\psi: \R^{d} \times \subext{k}\R^{d} \to [0,\infty)$ be a positive convex integrand. To finish our discussion of the $\psi$-mass, we establish an estimate for the $\psi$-mass under pushforwards. 

\begin{lemma}\label{lm:psiMassPushFwd} Let $\psi: \R^{d} \times \subext{k}\R^{d} \to [0,\infty)$ be a positive convex integrand. Let $f: \R^{d} \to \R^{d}$ be Lipschitz and take $T \in \Nrm_{k}(\R^{d})$. Then
\begin{align*}
    \Mbf_{\psi}(f_{*}T) &\leq \int \psi(f(x),D^{T}f(x)\vec{T}(x)) \dd\Vert T \Vert(x).
\end{align*}
\end{lemma}

\begin{proof} We first disintegrate $\Vert T \Vert$ with respect to the map $f$ (see \cite[Theorem 5.3.1]{AmbrosioGigliSavare05book}), which yields a Borel family $\{\Vert T \Vert_{y}\}_{y \in \R^{d}}$ of probability measures such that
\begin{align}\label{psiMassPushFwd1}
    \int \varphi \dd\Vert T \Vert = \int_{\R^{d}} \int_{f^{-1}(y)} \varphi(x) \dd\Vert T \Vert_{y}(x) \dd f_{\#}\Vert T \Vert(y)
\end{align}
for all bounded Borel $\varphi: \R^{d} \to \R$. It follows that, for $\omega \in \Dcal^{k}(\R^{d})$,
\begin{align*}
    f_{*}T(\omega) = \int_{\R^{d}} \bigg\langle \int_{f^{-1}(y)} D^{T}f(x)\vec{T}(x)\dd \Vert T \Vert_{y}(x), \omega(y) \bigg\rangle \dd f_{\#}\Vert T \Vert(y).
\end{align*}
We therefore have
\begin{align*}
    \Mbf_{\psi}(f_{*}T) = \int_{\R^{d}} \psi\bigg(y,\int_{f^{-1}(y)} D^{T}f(x)\vec{T}(x)\dd \Vert T \Vert_{y}(x)\bigg) \dd f_{\#}\Vert T \Vert(y).
\end{align*}
The result follows after using Jensen's inequality to pass $\psi(y,\cdot)$ into the inner integral for $y \in \R^{d}$ and employing \eqref{psiMassPushFwd1}.
\end{proof}

Now let $\phi: \R^{d} \times \subext{1+k}\R^{d} \to [0,\infty)$ be a positive convex integrand. The \textbf{$\phi$-variation} of a space-time normal current $S \in \Nrm_{1+k}(\R^{1+d})$ on a Borel $I \subset \R$ is given by
\begin{align*}
    \Var_{\phi}(S;I) := \int_{I \times \R^{d}} \phi(x,\pbf \vec{S}(t,x)) \dd\Vert S \Vert(t,x).
\end{align*}
If $\tbf(\supp S) \subset I$, we set
\begin{align*}
    \Var_{\phi}(S) := \Var_{\phi}(S;I).
\end{align*}
Notice that $\Var_{\phi}(S;\cdot)$ is a finite Borel measure on $\R$ and that
\begin{align}\label{phiVarMeasure}
    \int g \dd\Var_{\phi}(S;\cdot) = \int g(t)\phi(x,\pbf\vec{S}(t,x)) \dd \Vert S \Vert(t,x)
\end{align}
for all bounded Borel $g: \R \to \R$. The next lemma will allow us to compute the $\phi$-variation of ``jumps'':

\begin{lemma}\label{lm:DissJumps} Let $\phi: \R^{d} \times \subext{1+k}\R^{d} \to [0,\infty)$ be a positive convex integrand. Let $R \in \Nrm_{1+k}(\R^{d})$ and take $a \in \R$. Then
\begin{align*}
    \Var_{\phi}(\delta_{a} \times R) = \Mbf_{\phi}(R).
\end{align*}
\end{lemma}

On the other hand, we have

\begin{lemma}\label{lm:pPhiMassByDiss} Let $\phi: \R^{d} \times \subext{1+k}\R^{d} \to [0,\infty)$ be a positive convex integrand. Take $S \in \Nrm_{1+k}(\R^{1+d})$ and let $I \subset \R$ be Borel. Then
\begin{align*}
    \Mbf_{\phi}(\pbf_{*}(S \mvert I \times \R^{d})) &\leq \Var_{\phi}(S;I).
\end{align*}
\end{lemma}
\begin{proof} Take $\omega \in \Dcal^{1+k}(\R^{d})$ with $\omega(x) \leq \phi(x,\cdot)$ for $x \in R^{d}$. Then
\begin{align*}
    \pbf_{*}(S \mvert I \times \R^{d})(\omega) &= \int_{I \times \R^{d}} \langle \pbf\vec{S}(t,x),\omega(x) \rangle \dd \Vert S \Vert(t,x) \\
    &\leq \int_{I \times \R^{d}} \phi(x,\pbf\vec{S}(t,x)) \dd \Vert S \Vert(t,x) = \Var_{\phi}(S;I).
\end{align*}
After taking a supremum, the result follows from Lemma \ref{lm:psiDual}.
\end{proof}

Let $\phi: \R^{d} \times \subext{1+k}\R^{d} \to [0,\infty)$ be a positive convex integrand. It turns out that the $\phi$-variation enjoys properties which are analogous to those established for the variation in \cite[Section 3.3]{BonicattoRindlerTurnbull26}. Here, we recall the convention of \cite{BonicattoRindlerTurnbull26} and set $a_{*}S := (a \times \mathrm{id}_{\R^{d}})_{*}S$ when $S \in \Nrm_{1+k}(\R^{1+d})$ and $a: \R \to \R$ is Lipschitz.

\begin{proposition}\label{pp:DissRescale} Let $\phi: \R^{d} \times \subext{1+k}\R^{d} \to [0,\infty)$ be a positive convex integrand. Fix $\sigma,\tau \in \R$ with $\sigma < \tau$, take $S \in \Nrm_{1+k}([\sigma,\tau] \times \R^{d})$, and let $a \in \Lip([\sigma,\tau])$ be injective. Let $I \subset [\sigma,\tau]$ be Borel. Then
\begin{align*}
    \Var_{\phi}(a_{*}S;a(I)) = \Var_{\phi}(S;I)
\end{align*}
\end{proposition}

\begin{proof} The proof proceeds in the same way as \cite[Proposition 3.5]{BonicattoRindlerTurnbull26}.
\end{proof}

\begin{lemma}\label{lm:PushFwdDiss} Let $\phi: \R^{d} \times \subext{1+k}\R^{d} \to [0,\infty)$ be a positive convex integrand. Take $S \in \Nrm_{1+k}(\R^{1+d})$ and let $f: \R^{1+d} \to \R^{1+d}$ be Lipschitz. Let $I \subset \R$ be Borel. Then
\begin{align*}
    \Var_{\phi}(f_{*}S;I) &\leq \int_{f^{-1}(I \times \R^{d})} \phi((\pbf \circ f)(t,x),\pbf(D^{S}f(t,x)\vec{S}(t,x))) \dd\Vert S \Vert(t,x).
\end{align*}
\end{lemma}

\begin{proof} The proof proceeds in the same way as \cite[Lemma 3.7]{BonicattoRindlerTurnbull26}.
\end{proof}

\subsection{Space-Time Strict Convergence}\label{s:STStrict} Let $(S_{j})_{j} \subset \Nrm_{1+k}(\R^{1+d})$ be a sequence of space-time normal currents and take a space-time normal current $S \in \Nrm_{1+k}(\R^{1+d})$. We say that the sequence $(S_{j})_{j}$ converges to $S$ \textbf{space-time strictly} if
\begin{align*}
    S_{j} &\wkto S, & \Var(S_{j}) &\to \Var(S), & \int \vert \Mbf(S_{j}(t)) - \Mbf(S(t)) \vert\dd t &\to 0.
\end{align*}
The following result provides us with a useful criterion for establishing space-time strict convergence:
\begin{proposition}\label{pp:STSSuffCond} Take positive convex integrands
\begin{align*}
    \psi &: \R^{d} \times \subext{k}\R^{d} \to [0,\infty), & \phi &: \R^{d} \times \subext{1+k}\R^{d} \to [0,\infty).
\end{align*}
Let $K \subset \R^{1+d}$ be compact. Take a sequence $(S_{j})_{j} \subset \Nrm_{1+k}(K)$ and suppose that $(\Mbf(\partial S_{j}))_{j}$ is uniformly bounded. Let $S \in \Nrm_{1+k}(K)$ and suppose that $S_{j} \wkto S$,
\begin{align}\label{a:STSSuffCond}
    \limsup_{j \to \infty} \Var_{\phi}(S_{j}) &\leq \Var_{\phi}(S), & \limsup_{j \to \infty} \int \Mbf_{\psi}(S_{j}(t)) \dd t &\leq \int \Mbf_{\psi}(S(t)) \dd t.
\end{align}
Then
\begin{align*}
    \Var_{\phi}(S_{j}) &\to \Var_{\phi}(S), & \int \vert \Mbf_{\psi}(S_{j}(t)) - \Mbf_{\psi}(S(t)) \vert \dd t &\to 0.
\end{align*}
\end{proposition}

\begin{proof} Take a (not relabelled) subsequence of $(S_{j})_{j}$. We can use \cite[Lemma 3.10]{BonicattoRindlerTurnbull26} to deduce that
\begin{equation}\label{STSSuffCond1}
\begin{aligned}
    \limsup_{j \to \infty} \Mbf(S_{j}) &\leq C_{\psi,\pbf(K)}\int \Mbf_{\psi}(S(t)) \dd t + C_{\phi,\pbf(K)} \cdot \Var_{\phi}(S)
\end{aligned}
\end{equation}
where $C_{\phi,\pbf(K)}, C_{\psi,\pbf(K)} > 0$ are defined as in \eqref{psiVsMass}.
Since $(\Mbf(\partial S_{j}))_{j}$ is uniformly bounded by assumption, we can appeal to \cite[Corollary 3.12]{BonicattoRindlerTurnbull26} to pass to a (not relabelled) subsequence of $(S_{j})_{j}$ along which we have $S_{j}(t) \wkto S(t)$ for $\mathscr{L}^{1}$-a.e.\ $t \in \R$. We can then use Lemma \ref{lm:psiLSC} to deduce that 
\begin{align*}
    \Mbf_{\psi}(S(t)) \leq \liminf_{j \to \infty} \Mbf_{\psi}(S_{j}(t)) \hspace{1em} \text{for $\mathscr{L}^{1}$-a.e.\ $t \in \R$}.
\end{align*}
Notice that, for $j \in \N$, we have
\begin{equation}\label{STSSuffCond2}
\begin{aligned}
    \vert \Mbf_{\psi}(S_{j}&(t)) - \Mbf_{\psi}(S(t)) \vert \\
    &= (\Mbf_{\psi}(S_{j}(t)) - \Mbf_{\psi}(S(t)))+2\max\{0,\Mbf_{\psi}(S(t))-\Mbf_{\psi}(S_{j}(t))\}.
\end{aligned}
\end{equation}
We then integrate \eqref{STSSuffCond2} and estimate (using the second assumption of \eqref{a:STSSuffCond} on the first term and the reverse Fatou lemma on the second term) to obtain
\begin{align*}
    \limsup_{j \to \infty}\int \vert &\Mbf_{\psi}(S_{j}(t)) - \Mbf_{\psi}(S(t)) \vert \dd t \\
    &\leq 2 \int \limsup_{j \to \infty}\max\{0,\Mbf_{\psi}(S(t)) - \Mbf_{\psi}(S_{j}(t))\}\dd t \\
    &\leq 2 \int \max\bigg\{0,\Mbf_{\psi}(S(t))-\liminf_{j \to \infty}\Mbf_{\psi}(S_{j}(t))\bigg\} \dd t \\
    &= 0.
\end{align*}
Furthermore, it follows from Reshetnyak's lower semicontinuity theorem (see, for example, \cite[Theorem 2.38]{AmbrosioFuscoPallara00book}) that
\begin{align*}
    \Var_{\phi}(S) &\leq \liminf_{j \to \infty} \Var_{\phi}(S_{j}).
\end{align*}
The result follows since the initial choice of subsequence was arbitrary.
\end{proof}

\begin{example}\label{Ex:1} Let $(S_{j})_{j} \subset \Nrm_{1+k}(\R^{1+d})$. The assumption that $(\Mbf(\partial S_{j}))_{j}$ is uniformly bounded cannot be removed from Proposition \ref{pp:STSSuffCond}. To see this, consider
\begin{align*}
    S_{j} := \sum_{i=0}^{j-1} \bigg[\!\!\bigg[\frac{i}{j},\frac{2i+1}{2j}\bigg]\!\!\bigg] \in \Nrm_{1}(\R^{1+0}) \hspace{1em} \text{for $j \in \N$}.
\end{align*}
One can check that 
\begin{align*}
    \Var(S_{j}) &= 0, & \int \Mbf(S_{j}(t)) \dd t &= \frac{1}{2} \hspace{1em} \text{for $j \in \N$}.
\end{align*}
Notice that $S_{j} \wkto S$ where
\begin{align*}
    S = \frac{1}{2} \cdot [\![0,1]\!] \in \Nrm_{1}(\R^{1+0}).
\end{align*}
One can check that
\begin{align*}
    \Var(S) &= 0, & \int \Mbf(S(t)) \dd t &= \frac{1}{2}.
\end{align*}
On the other hand,
\begin{align*}
    \int \vert \Mbf(S_{j}(t))-\Mbf(S(t)) \vert \dd t = \frac{1}{2} \hspace{1em} \text{for $j \in \N$}.
\end{align*}
In particular, $S_{j} \not\to S$ space-time strictly.
\end{example}

In what follows, we shall need two auxiliary results, the first of which is a generalisation of \cite[Proposition 3.2]{BonicattoRindlerTurnbull26}:

\begin{lemma}\label{lm:Coarea} Let $S \in \Nrm_{1+k}(\R^{1+d})$ and take a bounded Borel $g: \R^{1+d} \to \R$. Then
\begin{align*}
    \int_{\R}\int_{\R^{d}} g(t,x) \dd\Vert S (t)\Vert(x) \dd t = \int_{\R^{1+d}} g(t,x) \Vert \vec{S}(t,x) \mvert \mathrm{d}t \Vert \dd\Vert S \Vert(t,x).
\end{align*}
\end{lemma}

\begin{proof} Let $\varphi: \R^{d} \to [0,\infty)$ and $f: \R \to [0,\infty)$ be smooth and compactly supported. One can check that
\begin{align*}
    \varphi \cdot S(t) = ((\varphi \circ \pbf) \cdot S)(t) \hspace{1em} \text{for $\mathscr{L}^{1}$-a.e.\ $t \in \R$}
\end{align*}
and so 
\begin{align*}
    \Mbf(((\varphi \circ \pbf) \cdot S)(t)) = \int \varphi(x) \dd \Vert S(t)\Vert(x) \hspace{1em} \text{for $\mathscr{L}^{1}$-a.e.\ $t \in \R$}.
\end{align*}
It follows that
\begin{align}\label{Coarea1}
    \int_{\R} f(t) \Mbf(((\varphi \circ \pbf) \cdot S)(t)) \dd t = \int_{\R}\int_{\R^{d}} f(t)\varphi(x) \dd\Vert S(t) \Vert(x) \dd t.
\end{align}
On the other hand, we can use \cite[Proposition 3.2]{BonicattoRindlerTurnbull26} with $(\varphi \circ \pbf) \cdot S$ to deduce that
\begin{align}\label{Coarea2}
    \int_{\R} f(t) \Mbf(((\varphi \circ \pbf) \cdot S)(t)) \dd t = \int_{\R^{1+d}} f(t) \Vert \vec{S}(t,x) \mvert \mathrm{d} t \Vert \varphi(x) \dd \Vert S \Vert(t,x).
\end{align}
Combining \eqref{Coarea1} and \eqref{Coarea2} yields the result in the case that $g(t,x)=f(t)\varphi(x)$ for $(t,x) \in \R^{1+d}$; the general case proceeds by approximation. 
\end{proof}

\begin{lemma}\label{lm:SlicePolar} Let $S \in \Nrm_{1+k}(\R^{1+d})$. Then there is a bounded Borel $\xi: \R^{1+d} \to \subext{k}\R^{d}$ such that
\begin{align}\label{cc:SlicePolar2}
    \xi(t,x) \hspace{0.2em} \Vert \vec{S}(t,x) \mvert \mathrm{d}t \Vert = \pbf(\vec{S}(t,x) \mvert \mathrm{d}t) \hspace{1em} \text{for $\Vert S \Vert$-a.e.\ $(t,x) \in \R^{1+d}$}. 
\end{align}
Moreover,
\begin{align}\label{cc:SlicePolar1}
    \xi(t,\cdot) = \vec{S}(t) \hspace{1em} \text{$\Vert S(t) \Vert$-a.e. for $\mathscr{L}^{1}$-a.e.\ $t \in \R$.}
\end{align}
\end{lemma}

\begin{proof} Let $f: \R \to [0,\infty)$ be smooth and compactly supported and take $\omega \in \Dcal^{k}(\R^{d})$. We can use \cite[4.3.2 (1)]{Federer69book} (see also the proof of \cite[Proposition 3.2]{BonicattoRindlerTurnbull26}) to deduce that
\begin{align*}
    \int f(t) \cdot S(t)(\omega) \dd t = \int f(t) \cdot S\vert_{t}(\pbf^{*}\omega) \dd t = (S\mvert (f \circ \tbf)\hspace{0.2em}\mathrm{d}t)(\pbf^{*}\omega).
\end{align*}
It follows that
\begin{align}\label{SlicePolar1}
\begin{aligned}
    \int_{\R} f(t) \int_{\R^{d}} \langle &\vec{S(t)}(x),\omega(x) \rangle \dd \Vert S(t) \Vert(x) \dd t \\
    &= \int_{\R^{1+d}} f(t) \langle \pbf(\vec{S}(t,x) \mvert \mathrm{d}t),\omega(x) \rangle \dd \Vert S \Vert(t,x).
\end{aligned}
\end{align}
After choosing a representative, we can assume that the map $(t,x) \mapsto \vec{S}(t,x) \mvert \mathrm{d}t$ is bounded and Borel. Define $\xi: \R^{1+d} \to \subext{k}\R^{d}$ by
\begin{align}\label{SlicePolar2}
    \xi(t,x) &:= \frac{\pbf(\vec{S}(t,x) \mvert \mathrm{d}t)}{\Vert \vec{S}(t,x) \mvert \mathrm{d}t \Vert} \hspace{1em} \text{if $\Vert \vec{S}(t,x) \mvert \mathrm{d}t \Vert \neq 0$}, & \xi(t,x) &:= 0 \hspace{1em} \text{otherwise}.
\end{align}
Notice that $\xi$ satisfies \eqref{cc:SlicePolar2} and that $\xi$ is bounded and Borel. We can use Lemma \ref{lm:Coarea} to deduce that 
\begin{align}\label{SlicePolar3}
\begin{aligned}
    \int_{\R^{1+d}} f(t) \langle \pbf(&\vec{S}(t,x) \mvert \mathrm{d}t),\omega(x) \rangle \dd \Vert S \Vert(t,x) \\
    &= \int_{\R} f(t) \int_{\R^{d}} \langle \xi(t,x),\omega(x) \rangle \dd \Vert S(t) \Vert(x) \dd t.
\end{aligned}
\end{align}
Combining \eqref{SlicePolar1} and \eqref{SlicePolar3} and using the fact that the choice of $f$ and $\omega$ was arbitrary yields \eqref{cc:SlicePolar1}. 
\end{proof}

Let $\psi: \R^{d} \times \subext{k}\R^{d} \to [0,\infty)$ be a positive convex integrand and take $S \in \Nrm_{1+k}(\R^{1+d})$. It follows from Lemma \ref{lm:Coarea} and Lemma \ref{lm:SlicePolar} that
\begin{align}\label{psiMassIntegral}
    \int \Mbf_{\psi}(S(t)) \dd t = \int_{\R}\int_{\R^{d}} \psi(x,\pbf(\vec{S}(t,x) \mvert \mathrm{d}t)) \dd\Vert S \Vert(t,x).
\end{align}

Let $K \subset \R^{d}$ be compact. Take $(T_{j})_{j} \subset \Nrm_{k}(K)$, let $T \in \Nrm_{k}(K)$, and suppose that $T_{j} \to T$ strictly. One can show (see \cite[Corollary 13.2]{Rindler26book}) that $\Vert T_{j} \Vert \wkto \Vert T \Vert$. We shall now establish an analogous result in our setting:

\begin{proposition}\label{pp:VarWkStar} Take positive convex integrands
\begin{align*}
    \psi &: \R^{d} \times \subext{k}\R^{d} \to [0,\infty), & \phi &: \R^{d} \times \subext{1+k}\R^{d} \to [0,\infty).
\end{align*}
Let $K \subset \R^{d}$ be compact. Take a sequence $(S_{j})_{j} \subset \Nrm_{1+k}([0,1] \times K)$ and let $S \in \Nrm_{1+k}([0,1] \times K)$. Suppose that $S_{j} \wkto S$ and that \eqref{a:STSSuffCond} holds. Then
\begin{align}\label{cc:VarWkStar1}
    \Var_{\phi}(S_{j};\cdot) \wkto \Var_{\phi}(S;\cdot)
\end{align}
and
\begin{align}\label{cc:VarWkStar2}
    \Mbf_{\psi}(S_{j}(\cdot)) \hspace{0.2em} \mathscr{L}^{1} \wkto \Mbf_{\psi}(S(\cdot)) \hspace{0.2em} \mathscr{L}^{1}.
\end{align}
\end{proposition}

\begin{proof} Take a (not relabelled) subsequence of $(S_{j})_{j}$. Since \eqref{a:STSSuffCond} holds, the sequences $(\Var_{\phi}(S_{j}))_{j}$ and $(\Mbf_{\psi}(S_{j}(\cdot)) \hspace{0.2em} \mathscr{L}^{1})_{j}$ are uniformly bounded. We can therefore find a (not relabelled) subsequence of $(\Var_{\phi}(S_{j}))_{j}$ and $(\Mbf_{\psi}(S_{j}(\cdot)) \hspace{0.2em} \mathscr{L}^{1})_{j}$ and $\nu,\mu \in \Mcal^{+}([0,1])$ such that
\begin{align}\label{VarWkStar1}
    \Var_{\phi}(S_{j};\cdot) &\wkto \nu, & \Mbf_{\psi}(S_{j}(\cdot)) \hspace{0.2em} \mathscr{L}^{1} &\wkto \mu,
\end{align}
in $\Mcal^{+}([0,1])$. Notice that we can estimate as in \eqref{STSSuffCond1} to deduce that the sequence $(\Mbf(S_{j}))_{j}$ is uniformly bounded. Since we also have $S_{j} \wkto S$, it follows that $S_{j}(\omega) \to S(\omega)$ for all $\omega \in \Crm_{0}(\R^{1+d};\supext{1+k}\R^{1+d})$. Let $g: [0,1] \to [0,\infty)$ be continuous. It follows from \eqref{phiVarMeasure} and Reshetnyak's lower semicontinuity theorem that
\begin{align*}
    \int g \dd \Var_{\phi}(S;\cdot) \leq \liminf_{j \to \infty} \int g \dd \Var_{\phi}(S_{j};\cdot) = \int g \dd \nu.
\end{align*}
Similarly, it follows from Lemma \ref{lm:Coarea}, Lemma \ref{lm:SlicePolar}, and Reshetnyak's lower semicontinuity theorem that
\begin{align*}
    \int g(t) \Mbf_{\psi}(S(t)) \dd t \leq \liminf_{j \to \infty} \int g(t) \Mbf_{\psi}(S_{j}(t)) \dd t = \int g \dd\mu.
\end{align*}
In particular, we have $\Var_{\phi}(S;\cdot) \leq \nu$ and $\Mbf_{\psi}(S(\cdot)) \hspace{0.2em} \mathscr{L}^{1} \leq \mu$. On the other hand, testing with $\mathbf{1}_{[0,1]}$ in \eqref{VarWkStar1} yields
\begin{align*}
    \Var_{\phi}(S_{j}) &\to \nu([0,1]), & \int \Mbf_{\psi}(S_{j}(t)) \dd t \to \mu([0,1]).
\end{align*}
It then follows from \eqref{a:STSSuffCond} that
\begin{align*}
    \nu([0,1]) &\leq \Var_{\phi}(S), & \mu([0,1]) &\leq \int \Mbf_{\psi}(S(t)) \dd t.
\end{align*}
We therefore have $\Var_{\phi}(S;\cdot) = \nu$ and $\Mbf_{\psi}(S(\cdot)) \hspace{0.2em} \mathscr{L}^{1} = \mu$. The result follows since the choice of the initial subsequence was arbitrary. 
\end{proof}

Now set
\begin{align}\label{Intervals1}
    A_{j} &:= \bigcup_{i=0}^{j-1} \bigg[\frac{i}{j},\frac{2i+1}{2j}\bigg] \subset [0,1] \hspace{1em} \text{for $j \in \N$}.
\end{align}

The next two examples establish that, in general, there is no relationship between the strict and space-time strict convergences. 

\begin{example}\label{Ex:2} Consider the sequence $(S_{j})_{j} \subset \Nrm_{1}([0,1]^{2})$ where, for $j \in \N$,
\begin{align*}
    S_{j} := [\![0,1]\!] \times (\mathscr{L}^{1} \mvert A_{j}) + (\mathscr{L}^{1} \mvert A_{j}) \times [\![0,1]\!].
\end{align*}
One can check that, for $j \in \N$, we have
\begin{align*}
    \Var(S_{j}) &= \frac{1}{2}, & \Mbf(S_{j}(t)) &= \frac{1}{2} \hspace{1em} \text{for $\mathscr{L}^{1}$-a.e.\ $t \in [0,1]$}.
\end{align*}
Notice that $S_{j} \wkto S$ where
\begin{align*}
    S := \frac{e_{0}+e_{1}}{2} \cdot \mathscr{L}^{2} \mvert [0,1]^{2}.
\end{align*}
and one can check that
\begin{align*}
     \Var(S) &= \frac{1}{2}, & \Mbf(S(t)) &= \frac{1}{2} \hspace{1em} \text{for $\mathscr{L}^{1}$-a.e.\ $t \in [0,1]$}.
\end{align*}
It follows that $S_{j} \to S$ space-time strictly. However, we have
\begin{align*}
    \Mbf(S_{j}) &= \frac{1}{2} + \frac{\sqrt{2}}{4} \hspace{1em} \text{for $j \in \N$}, & \Mbf(S) &= \frac{1}{\sqrt{2}}.
\end{align*}
In particular, $S_{j} \not\to S$ strictly.
\end{example}

\begin{example}\label{Ex:3} Let $v_{1},v_{2} \in \subext{2}\R^{1+3}$ given by
\begin{align*}
    v_{1} &:= e_{0} \wedge e_{1} + e_{2} \wedge e_{3}, & v_{2} := e_{0} \wedge e_{2} + e_{1} \wedge e_{3}.
\end{align*}
Notice that $\star v_{1} = v_{1}$, $\star v_{2} = -v_{2}$ and that $\vert v_{1} \vert = \vert v_{2} \vert = \sqrt{2}$, so we can use \cite[Corollary 3.4]{Morgan85} to deduce that
\begin{align*}
    \Vert v_{1} + v_{2} \Vert = \Vert v_{1} - v_{2} \Vert = \Vert v_{1} \Vert = 2.
\end{align*}
Now set 
\begin{align*}
    R_{1} &:= (v_{1}+v_{2}) \mathscr{L}^{4} \mvert [0,1]^{4}, & R_{2} &= (v_{1}-v_{2}) \mathscr{L}^{4} \mvert [0,1]^{4},
\end{align*}
and notice that 
\begin{align}\label{ThirdExample1}
    \Mbf(R_{1}) = \Mbf(R_{2}) = \Mbf\bigg(\frac{R_{1}+R_{2}}{2}\bigg) = 2.
\end{align}
For $j \in \N$, define 
\begin{align*}
    S_{j} := R_{1} \times (\mathscr{L}^{1} \mvert A_{j}) + R_{2} \times (\mathscr{L}^{1} \mvert [0,1] \backslash A_{j}).
\end{align*}
where $A_{j} \subset [0,1]$ is given by \eqref{Intervals1}. Then $S_{j} \wkto S$, where
\begin{align*}
    S = \bigg(\frac{R_{1}+R_{2}}{2}\bigg) \times (\mathscr{L}^{1} \mvert [0,1]).
\end{align*}
It follows from \eqref{ThirdExample1} that $S_{j} \to S$ strictly. For $j \in \N$ and $\mathscr{L}^{1}$-a.e.\ $t \in [0,1]$, one can check that
\begin{align*}
    S_{j}(t) &= ((e_{1}+e_{2}) \mathscr{L}^{3} \mvert [0,1]^{3}) \times (\mathscr{L}^{1} \mvert A_{j}) \\
    &\qquad + ((e_{1}-e_{2}) \mathscr{L}^{3} \mvert [0,1]^{3}) \times (\mathscr{L}^{1} \mvert [0,1] \backslash A_{j})
\end{align*}
and so $\Mbf(S_{j}(t)) = \sqrt{2}$. Furthermore, for $\mathscr{L}^{1}$-a.e.\ $t \in [0,1]$, we have
\begin{align*}
    S(t) = (e_{1} \mathscr{L}^{3} \mvert [0,1]^{3}) \times (\mathscr{L}^{1} \mvert [0,1]).
\end{align*}
and so $\Mbf(S(t)) = 1$ for $\mathscr{L}^{1}$-a.e.\ $t \in [0,1]$. In particular,
\begin{align*}
    \int \vert \Mbf(S_{j}(t))-\Mbf(S(t)) \vert \dd t = \sqrt{2} - 1 \hspace{1em} \text{for $j \in \N$}
\end{align*}
and so $S_{j} \not\to S$ space-time strictly.
\end{example}

The underlying obstruction in Example \ref{Ex:3} is the same as that encountered in Section \ref{s:Inhomogeneous}: The space $(\subext{1+k}\R^{1+d},\Vert \cdot \Vert)$ is not strictly convex in general.

Furthermore, notice that $(\Mbf(\partial S_{j}))_{j}$ is uniformly bounded in both Example \ref{Ex:2} and Example \ref{Ex:3}, so imposing conditions of this kind alone does not yield a positive result relating the classical and space-time strict convergences. Nevertheless, we have the following partial result:

\begin{proposition}\label{pp:SimpleSTS} Take positive convex integrands
\begin{align*}
    \psi &: \R^{d} \times \subext{k}\R^{d} \to [0,\infty), & \phi &: \R^{d} \times \subext{1+k}\R^{d} \to [0,\infty).
\end{align*}
Let $K \subset \R^{1+d}$ be compact. Take a sequence $(S_{j})_{j} \subset \Nrm_{1+k}(K)$ of simple currents and suppose that $(\Mbf(\partial S_{j}))_{j}$ is uniformly bounded. Let $S \in \Nrm_{1+k}(K)$ be a simple current and suppose that $S_{j} \to S$ strictly. Then
\begin{align*}
    \Var_{\phi}(S_{j}) &\to \Var_{\phi}(S), & \int \vert \Mbf_{\psi}(S_{j}(t)) - \Mbf_{\psi}(S(t)) \vert \dd t &\to 0
\end{align*}
\end{proposition}

\begin{proof} We first appeal to Proposition \ref{pp:Simple} to deduce that
\begin{align*}
    \int \phi(x,\pbf\vec{S_{j}}(t,x)) \dd\Vert S_{j} \Vert(t,x) \to \int \phi(x,\pbf\vec{S}(t,x)) \dd \Vert S \Vert(t,x)
\end{align*}
and
\begin{align*}
    \int \psi(x,\pbf(\vec{S_{j}}(t,x) \mvert \mathrm{d}t)) \dd\Vert S_{j} \Vert(t,x) \to \int \psi(x,\pbf(\vec{S}(t,x) \mvert \mathrm{d}t)) \dd\Vert S \Vert(t,x).
\end{align*}
It follows from the definition of the $\phi$-variation that $\Var_{\phi}(S_{j}) \to \Var_{\phi}(S)$ and we can use \eqref{psiMassIntegral} to deduce that
\begin{align*}
    \int \Mbf_{\psi}(S_{j}(t)) \dd t \to \int \Mbf_{\psi}(S(t)) \dd t.
\end{align*}
The result then follows from Proposition \ref{pp:STSSuffCond}.
\end{proof}

\subsection{The Deformation Theorem} Let $T \in \Nrm_{k}(\R^{d})$ and let $f: \R^{d} \to \R^{d}$ be Lipschitz. One can view the pushforward $f_{*}T$ of $T$ under $f$ as a `deformation' of $T$ using the homotopy formula: If $H: [0,1] \times \R^{d} \to \R^{d}$ is a (Lipschitz) homotopy from $\mathrm{id}_{\R^{d}}$ to $f$, then
\begin{align}\label{Homotopy}
    f_{*}T = T + H_{*}([\![0,1]\!] \times \partial T) + \partial H_{*}([\![0,1]\!] \times T).
\end{align}
By refining this procedure with a well-chosen map adapted to $T$, one can deform $T$ into a polyhedral chain $P \in \Prm_{k}(\R^{d})$ with support on a cubical `skeleton'. This is the content of the \textit{deformation theorem}; see, for example, \cite{Federer69book,KrantzParks08book,Simon14Notes,White99} and also \cite{BonicattoRindlerTurnbull26} for a recent treatment in the case $\partial T = 0$. While this is not a strict approximation result, it is implicitly used in many of the approximation results presented in the sequel, so we give a statement of the theorem below:

\begin{proposition}[Deformation Theorem]\label{pp:Deformation} Let $T \in \Nrm_{k}(\R^{d})$ and fix $\varepsilon > 0$. Then there are $P \in \Srm_{k,\varepsilon}(\R^{d})$, $Q \in \Nrm_{k}(\R^{d})$, and $R \in \Nrm_{k+1}(\R^{d})$ with
\begin{align*}
    T = P + Q + \partial R.
\end{align*}
Furthermore, we have the following conclusions:
\begin{enumerate}
    \item We have the following mass estimates with a constant $\gamma > 0$ depending only on $k$ and $d$:
    \begin{align*}
        \Mbf(P) &\leq \gamma(\Mbf(T)+\varepsilon \cdot \Mbf(\partial T)), & \Mbf(Q) &\leq \gamma\varepsilon \cdot \Mbf(\partial T), \\
        \Mbf(\partial P) &\leq \gamma \cdot \Mbf(\partial T), & \Mbf(R) &\leq \gamma\varepsilon \cdot \Mbf(T).
    \end{align*}
    \item We have
    \begin{align*}
        \supp P, \hspace{0.2em} \supp R &\subset (\supp T)_{2d\varepsilon}, & \supp \partial P, \hspace{0.2em} \supp Q &\subset (\supp \partial T)_{2d\varepsilon}.
    \end{align*}
    \item If $\partial T \in \Irm_{k-1}(\R^{d})$, then $Q \in \Irm_{k}(\R^{d})$.
    \item If $T \in \Irm_{k}(\R^{d})$, then $Q \in \Irm_{k}(\R^{d})$, $R \in \Irm_{k+1}(\R^{d})$, and 
    \begin{align*}
        P \in \IPrm_{k}(\R^{d}) &:= \Irm_{k}(\R^{d}) \cap \Prm_{k}(\R^{d}).
    \end{align*}
    \item If $\partial T \in \Prm_{k-1}(\R^{d})$, then $Q \in \Prm_{k}(\R^{d})$.
    \item If $T \in \Prm_{k}(\R^{d})$, then $Q \in \Prm_{k}(\R^{d})$ and $R \in \Prm_{k+1}(\R^{d})$.
\end{enumerate}
\end{proposition}

Here, for $\varepsilon > 0$, the space $\Srm_{k,\varepsilon}(\R^{d})$ is defined in \cite[Section 2.2]{BonicattoRindlerTurnbull26}. A proof of the deformation theorem in this form can be found in \cite[4.2.9]{Federer69book}.

\section{Strict Polyhedral Approximation of Normal and Integral Currents}\label{s:ClassicalPolyStrict}

In this section, we review some classical results concerning the strict polyhedral approximation of normal and integral currents and use these results to prove Theorem \ref{tm:PolyBdryStrict}.

\subsection{Normal Currents}\label{s:NormalPolyStrict} Let $K \subset 
\R^{d}$ be compact. Here we recall that, for $T \in \Dcal_{k}(\R^{d})$ with $\supp T \subset K$, the \textbf{flat norm} of $T$ (in $K$) is given by
\begin{align}\label{FlatNormDefn}
    \Fbf_{K}(T) := \inf\{\Mbf(T - \partial R)+ \Mbf(R): R \in \Dcal_{k+1}(\R^{d}), \hspace{0.2em} \supp R \subset K\}.
\end{align}
Notice that, if $T \in \Nrm_{k}(K)$, then \eqref{FlatNormDefn} agrees with \eqref{NormalFlatNormDefn}. Our first strict approximation result is the following:

\begin{lemma}\label{lm:SimpleConstPolyStrict} Let $K \subset \R^{d}$ be compact and take $\varepsilon > 0$. Let $\xi \in \subext{k}\R^{d}$ be simple, take an $\mathscr{L}^{d}$-measurable $A \subset \Int K$, and set
\begin{align*}
    T := \xi \cdot \mathscr{L}^{d} \mvert A.
\end{align*}
Then there is $P \in \Prm_{k}(\R^{d})$ with 
\begin{align*}
    \Fbf_{K}(P-T) &\leq \varepsilon, & \Mbf(P) &\leq \Mbf(T) + \varepsilon.
\end{align*}
\end{lemma}
To prove Lemma \ref{lm:SimpleConstPolyStrict}, we shall follow some of the arguments in \cite[4.1.23]{Federer69book}. We first treat a special case which can be interpreted as a generalisation of the fact that
\begin{align*}
    \frac{1}{N} \sum_{j=0}^{N-1} \delta_{j/N} \to \mathscr{L}^{1} \mvert [0,1] \hspace{1em} \text{strictly}.
\end{align*}
\begin{lemma}\label{lm:SimpleConstPolyStrict1} Let $A \subset \R^{d}$ be compact and take $\varepsilon > 0$. Suppose that
\begin{align}\label{Rectangle}
    A = [a_{1},b_{1}] \times \cdots \times [a_{d},b_{d}]
\end{align}
for some $a_{1},\dotsc,a_{d},b_{1},\dotsc,b_{d} \in \R$. Let $T := (e_{1} \wedge \cdots \wedge e_{k}) \cdot \mathscr{L}^{d} \mvert A$. Then there is $P \in \Prm_{k}(\R^{d})$ with $\supp P \subset A$ and
\begin{align*}
    \Fbf_{A}(P-T) &\leq \varepsilon, & \Mbf(P) &= \Mbf(T).
\end{align*}
\end{lemma}

\begin{proof} Let $B := [a_{k+1},b_{k+1}] \times \cdots \times [a_{d},b_{d}]$ and notice that $T = R \times Y$ with
\begin{align*}
    R &= [\![a_{1},b_{1}]\!] \times \cdots \times [\![a_{k},b_{k}]\!] \in \Prm_{k}(\R^{k}), & Y &= \mathscr{L}^{d-k} \mvert B \in \Nrm_{0}(\R^{d-k}).
\end{align*}
Fix $\delta > 0$ (to be chosen later). Since $B$ is bounded, we can find $N \in \N$ and a collection $\{B_{i}\}_{i=1}^{N}$ of bounded Borel sets which partition $B$ and which satisfy $\diam(B_{i}) \leq \delta$ for $1 \leq i \leq N$. For $1 \leq i \leq N$, choose points $x_{i} \in B_{i}$. Let
\begin{align*}
    Q := \sum_{i=1}^{N} \mathscr{L}^{d-k}(B_{i})\delta_{x_{i}} \in \Prm_{0}(\R^{d-k})
\end{align*}
and then set $P := R \times Q$. Then, $P \in \Prm_{k}(\R^{d})$ and $\supp P = \supp R \times \supp Q \subset A$. Since $\Mbf(Q) = \Mbf(Y) = \mathscr{L}^{d-k}(B)$, we have
\begin{align*}
    \Mbf(P) = \Mbf(R)\Mbf(Q) = \Mbf(R)\Mbf(Y) = \Mbf(T).
\end{align*}
and so it suffices to show that $\Fbf_{A}(P-T) \leq \varepsilon$. To this end, fix $1 \leq i \leq N$. Take $\omega \in \Dcal^{0}(\R^{d-k})$ and suppose that 
\begin{align*}
    \sup \bigl\{\Vert \mathrm{d}\omega (x)\Vert: x \in B \bigr\} \leq 1.
\end{align*}
Then
\begin{align*}
    (\mathscr{L}^{d-k}\mvert B_{i} - \mathscr{L}^{d-k}(B_{i})\delta_{x_{i}})(\omega) \leq \int_{B_{i}} \vert \omega(x) - \omega(x_{i}) \vert \dd x \leq \diam(B_{i}) \cdot \mathscr{L}^{d-k}(B_{i})
\end{align*}
and so 
\begin{align*}
    (Y-Q)(\omega) \leq \sum_{i=1}^{N} \diam(B_{i}) \cdot \mathscr{L}^{d-k}(B_{i}) \leq \delta \cdot \mathscr{L}^{d-k}(B).
\end{align*}
It then follows from \cite[4.1.12]{Federer69book} that $\Fbf_{B}(Y-Q) \leq \delta \cdot \mathscr{L}^{d-k}(B)$. Since
\begin{align*}
    \Fbf_{A}(P-T) &\leq \Nbf(R) \cdot \Fbf_{B}(Y-Q),
\end{align*}
where
\begin{align*}
  \Nbf(R) = \Mbf(R) + \Mbf(\partial R),
\end{align*}
the result follows after choosing $\delta$ with $\delta \cdot \Nbf(R) \mathscr{L}^{d-k}(B) \leq \varepsilon$.
\end{proof}

We now relax the assumption \eqref{Rectangle}:

\begin{lemma}\label{lm:SimpleConstPolyStrict2} Let $K \subset \R^{d}$ be compact and take $\varepsilon > 0$. Let $A \subset \Int K$ be $\mathscr{L}^{d}$-measurable and let $T := (e_{1} \wedge \cdots \wedge e_{k}) \cdot \mathscr{L}^{d} \mvert A$. Then there is $P \in \Prm_{k}(\R^{d})$ with $\supp P \subset K$ and 
\begin{align*}
    \Fbf_{K}(P-T) &\leq \varepsilon, & \Mbf(P) &\leq \Mbf(T) + \varepsilon.
\end{align*}
    
\end{lemma}

\begin{proof} Since $A \subset \Int K$ is $\mathscr{L}^{d}$-measurable, there is $N \in \N$ and a $\mathscr{L}^{d}$-disjoint collection $\{A_{i}\}_{i=1}^{N}$ of rectangles (i.e., each $A_{i}$ is a product of $d$ intervals as in \eqref{Rectangle}) such that $A_{i} \subset K$ for $1 \leq i \leq N$ and 
\begin{align*}
    \mathscr{L}^{d}\bigg(A \hspace{0.2em}\triangle \bigcup_{i=1}^{N} A_{i}\bigg) \leq \frac{\varepsilon}{2}.
\end{align*}
For $1 \leq i \leq N$, set $T_{i} := (e_{1} \wedge \cdots \wedge e_{k}) \cdot \mathscr{L}^{d} \mvert A_{i}$. Then
\begin{align*}
    \Mbf\bigg(T -\sum_{i=1}^{N} T_{i}\bigg) \leq \frac{\varepsilon}{2}.
\end{align*}
For $1 \leq i \leq N$, we apply Lemma \ref{lm:SimpleConstPolyStrict1} to $T_{i}$ to obtain $P_{i} \in \Prm_{k}(\R^{d})$ with $\supp P_{i} \subset A_{i}$ and
\begin{align*}
    \Fbf_{A_{i}}(P_{i}-T_{i}) &\leq \frac{\varepsilon}{2N}, & \Mbf(P_{i}) &\leq \Mbf(T_{i}) + \frac{\varepsilon}{2N}.
\end{align*}
The result follows after setting
\begin{align*}
    P &:= \sum_{i=1}^{N} P_{i}. \qedhere
\end{align*}
\end{proof}

\begin{proof}[Proof of Lemma \ref{lm:SimpleConstPolyStrict}] If $\xi = 0$, then $T = 0$ and there is nothing to prove. Otherwise, take $v_{1},\dots,v_{k} \in \R^{d}$ with $\{v_{1},\dots,v_{k}\}$ orthonormal and $\lambda \in \R$ with $\vert \lambda \vert = \vert \xi \vert$ satisfying
\begin{align*}
    \xi = \lambda \cdot  v_{1} \wedge \cdots \wedge v_{k}.
\end{align*}
Extend $\{v_{1},\dotsc,v_{k}\}$ to an orthonormal basis $\{v_{1},\dotsc,v_{d}\}$ of $\R^{d}$, and let $L: \R^{d} \to \R^{d}$ be the unique linear map with $L(v_{i})=e_{i}$ for $1 \leq i \leq d$. Then $L$ is orthogonal and so $\Lip(L) = \Lip(L^{-1}) = \vert\det(L)\vert = 1$. Notice that
\begin{align*}
    L_{*}T &= \lambda \cdot \tilde{T}, & \tilde{T} &= (e_{1} \wedge \cdots \wedge e_{k}) \cdot \mathscr{L}^{d} \mvert L(A).
\end{align*}
Since $L(A) \subset \Int L(K)$, we can appeal to Lemma \ref{lm:SimpleConstPolyStrict2} to find $\tilde{P} \in \Prm_{k}(\R^{d})$ with $\supp \tilde{P} \subset L(K)$ and 
\begin{align*}
    \Fbf_{L(K)}(\tilde{P} - \tilde{T}) &\leq \frac{\varepsilon}{\vert \lambda \vert}, & \Mbf(\tilde{P}) &\leq \Mbf(\tilde{T}) + \frac{\varepsilon}{\vert \lambda \vert}.
\end{align*}
The result follows after setting $P := \lambda \cdot (L^{-1})_{*}\tilde{P}.$
\end{proof}

Let $\psi: \R^{d} \times \subext{k}\R^{d} \to [0,\infty)$ be a positive convex integrand. Recall from Section \ref{s:Inhomogeneous} that we cannot expect to deduce $\psi$-strict approximation results from strict approximation results in general and so we need to work with the $\psi$-mass directly. However, the following result demonstrates that there is an obstruction to performing $\psi$-strict polyhedral approximation procedures in general:

\begin{proposition}\label{pp:FlatSharp} Let $\psi: \R^{d} \times \subext{k}\R^{d} \to [0,\infty)$ be a positive convex integrand. Suppose that, for all compact $K \subset \R^{d}$, for all $\varepsilon > 0$, and for all $T \in \Nrm_{k}(\R^{d})$ with $\supp T \subset \Int K$, there is a simple $R \in \Nrm_{k}(\R^{d})$ with $\supp R \subset K$ and
\begin{align*}
    \Fbf_{K}(R-T) &\leq \varepsilon, & \Mbf_{\psi}(R) &\leq \Mbf_{\psi}(T) + \varepsilon.
\end{align*}
Then, for $(x,\xi) \in \R^{d} \times \subext{k}\R^{d}$, we have
\begin{align}\label{cc:FlatSharp}
    \psi(x,\xi) = \inf\bigg\{\sum_{i=1}^{N} \psi(x,\xi_{i}): N \in \N, \hspace{0.2em} \xi = \sum_{i=1}^{N} \xi_{i}, \hspace{0.2em} \text{$\xi_{i}$ simple}\bigg\}.
\end{align}
\end{proposition}

Let $\psi: \R^{d} \times \subext{k}\R^{d} \to [0,\infty)$ be a positive convex integrand. Motivated by Proposition \ref{pp:FlatSharp}, we define the \textbf{flat envelope} $\psi^{\flat}: \R^{d} \times \subext{k}\R^{d} \to [0,\infty)$ of $\psi$ for $(x,\xi) \in \R^{d} \times \subext{k}\R^{d}$ via
\begin{align}\label{Flat}
    \psi^{\flat}(x,\xi) := \inf\bigg\{\sum_{i=1}^{N} \psi(x,\xi_{i}): N \in \N, \hspace{0.2em} \xi = \sum_{i=1}^{N} \xi_{i}, \hspace{0.2em} \text{$\xi_{i}$ simple}\bigg\}.
\end{align}
Take $(x,\xi) \in \R^{d} \times \subext{k}\R^{d}$. Notice that $\psi(x,\xi) \leq \psi^{\flat}(x,\xi)$ and that
\begin{align}\label{FlatSimple}
    \psi(x,\xi) = \psi^{\flat}(x,\xi) \hspace{1em} \text{if $\xi$ is simple.}
\end{align}
It follows that, if $\psi: \R^{d} \times \subext{k}\R^{d} \to [0,\infty)$ is a positive convex integrand, then $(\psi^{\flat})^{\flat} = \psi^{\flat}$ and so $\psi^{\flat}$ is a flat integrand.

To prove Proposition \ref{pp:FlatSharp}, we shall need the following lemma, the statement of which can be found in \cite[2.7]{Federer74}:
\begin{lemma}\label{lm:FlatIsIntegrand} Let $\psi: \R^{d} \times \subext{k}\R^{d} \to [0,\infty)$ be a positive convex integrand. Then $\psi^{\flat}:  \R^{d} \times \subext{k}\R^{d} \to [0,\infty)$ is a positive convex integrand. 
\end{lemma}

\begin{proof} Fix $x \in \R^{d}$. One can use the fact that $\psi(x,\cdot)$ is a norm on $\subext{k}\R^{d}$ to deduce that $\psi^{\flat}(x,\cdot)$ is a norm on $\subext{k}\R^{d}$. To see that $\psi^{\flat}$ is continuous, fix $(x_{0},\xi_{0}) \in \R^{d} \times \subext{k}\R^{d}$ and let $\varepsilon > 0$. Take $\alpha \in (0,1)$ (to be chosen later). Since $\psi(x_{0},\cdot)$ is a norm on $\subext{k}\R^{d}$, $K := \{\xi \in \subext{k}\R^{d}: \psi(x_{0},\xi) = 1\}$ is compact and so there is $\delta > 0$ such that, for $x \in \R^{d}$,
\begin{align*}
    \vert x - x_{0} \vert \leq \delta \Rightarrow \vert \psi(x,\xi)-\psi(x_{0},\xi)\vert \leq \alpha \hspace{1em} \text{for $\xi \in K$}.
\end{align*}
It follows that 
\begin{align*}
     \vert x - x_{0} \vert \leq \delta \Rightarrow \vert \psi(x,\xi)-\psi(x_{0},\xi)\vert \leq \alpha \cdot \psi(x_{0},\xi) \hspace{1em} \text{for $\xi \in \subext{k}\R^{d}$}.
\end{align*}
One can then use \eqref{Flat} to deduce that
\begin{align}\label{FlatIsIntegrand1}
    \vert x - x_{0} \vert \leq \delta \Rightarrow \vert \psi^{\flat}(x,\xi)-\psi^{\flat}(x_{0},\xi)\vert \leq \alpha \cdot \psi^{\flat}(x_{0},\xi) \hspace{1em} \text{for $\xi \in \subext{k}\R^{d}$}.
\end{align}
Now let $\beta > 0$ (to be chosen later) and take $(x,\xi) \in \R^{d} \times \subext{k}\R^{d}$ with $\vert x - x_{0} \vert \leq \delta$ and $\Vert \xi - \xi_{0} \Vert \leq \beta$. We can then use the fact that $\psi^{\flat}(x,\cdot)$ is a norm on $\subext{k}\R^{d}$ and \eqref{FlatIsIntegrand1} to see
\begin{align*}
    \vert \psi^{\flat}(x,\xi)-\psi^{\flat}(x_{0},\xi_{0}) \vert &\leq \psi^{\flat}(x,\xi-\xi_{0}) + \vert \psi^{\flat}(x,\xi_{0})-\psi^{\flat}(x_{0},\xi_{0})\vert \\
    &\leq (1+\alpha) \cdot \psi^{\flat}(x_{0},\xi-\xi_{0}) + \alpha \cdot \psi^{\flat}(x_{0},\xi_{0}) \\
    &\leq \varepsilon.
\end{align*}
The result follows after choosing $\alpha,\beta > 0$ so that
\begin{align*}
    2\alpha \cdot \psi^{\flat}(x_{0},\xi_{0}) &\leq \varepsilon, & 4\beta \cdot \Lip(\psi^{\flat}(x_{0},\cdot)) &\leq \varepsilon. \qedhere
\end{align*}
\end{proof}

\begin{proof}[Proof of Proposition \ref{pp:FlatSharp}] Let $R \in \Nrm_{k}(\R^{d})$ be simple. Then it follows from \eqref{FlatSimple} that
\begin{align}\label{FlatMassSimple}
    \Mbf_{\psi}(R) = \Mbf_{\psi^{\flat}}(R).
\end{align}
Now take $j \in \N$, let $A \subset \R^{d}$ be a rectangle as in \eqref{Rectangle}, take a compact set $K \subset \R^{d}$ with $A \subset K^{\circ}$, let $\xi \in \subext{k}\R^{d}$, and set
\begin{align*}
    T := \xi \cdot \mathscr{L}^{d} \mvert A.
\end{align*}
We can then find simple $R_{j} \in \Nrm_{k}(\R^{d})$ with $\supp R_{j} \subset K$ and 
\begin{align*}
    \Fbf_{K}(R_{j} - T) &\leq \frac{1}{j}, & \Mbf_{\psi}(R_{j}) &\leq \Mbf_{\psi}(T) + \frac{1}{j}.
\end{align*}
It then follows from \eqref{FlatMassSimple} that
\begin{align}\label{FlatSharp1}
    \limsup_{j \to \infty}\Mbf_{\psi^{\flat}}(R_{j}) \leq \Mbf_{\psi}(T).
\end{align}
Since $\psi^{\flat}$ is a positive convex integrand and positively $1$-homogeneous for $x \in \R^{d}$, we can use Lemma \ref{lm:psiLSC} to deduce that
\begin{align}\label{FlatSharp2}
    \Mbf_{\psi^{\flat}}(T) \leq \liminf_{j \to \infty} \Mbf_{\psi^{\flat}}(R_{j}).
\end{align}
Combining \eqref{FlatSharp1} and \eqref{FlatSharp2} yields
\begin{align*}
    \int_{A} \psi^{\flat}(x,\xi) \dd x \leq \int_{A} \psi(x,\xi) \dd x.
\end{align*}
The result follows since $\psi$ and $\psi^{\flat}$ are continuous (see Lemma \ref{lm:FlatIsIntegrand}) and because the choice of $A \subset \R^{d}$ was arbitrary.
\end{proof}

Taking inspiration from \cite[2.7]{Federer74}, we say that a positive convex integrand $\psi: \R^{d} \times \subext{k}\R^{d} \to [0,\infty)$ is a \textbf{flat integrand} if \eqref{cc:FlatSharp} holds. We end this section with the following strict approximation result:

\begin{proposition}\label{pp:NormalPolyStrict} Let $\psi: \R^{d} \times \subext{k}\R^{d} \to [0,\infty)$ be a flat integrand. Let $K \subset \R^{d}$ be compact and take $\varepsilon > 0$. Let $T \in \Nrm_{k}(\R^{d})$ and suppose that $\supp T \subset \Int K$. Then there is $P \in \Prm_{k}(\R^{d})$ with $\supp P \subset K$ and
\begin{align*}
    \FN_{K}(P-T) &\leq \varepsilon, & \Mbf_{\psi}(P) &\leq \Mbf_{\psi}(T) + \varepsilon.
\end{align*}
\end{proposition}

The statement of Proposition \ref{pp:NormalPolyStrict} can be found in \cite[2.8]{Federer74}; in what follows, we give a proof, beginning with the following special case:

\begin{lemma}\label{lm:NormalPolyStrictIntermediate} Let $\psi: \R^{d} \times \subext{k}\R^{d} \to [0,\infty)$ be a flat integrand. Let $K \subset \R^{d}$ be compact and take $\varepsilon > 0$. Let $\xi \in \subext{k}\R^{d}$, take an $\mathscr{L}^{d}$-measurable $A \subset \Int K$, and set 
\begin{align*}
    T := \xi \cdot \mathscr{L}^{d} \mvert A.
\end{align*}
Then there is $P \in \Prm_{k}(\R^{d})$ with
\begin{align*}
    \Fbf_{K}(P-T) &\leq \varepsilon, & \Mbf_{\psi}(P) &\leq \Mbf_{\psi}(T) + \varepsilon.
\end{align*}
\end{lemma}

\begin{proof} If $T = 0$, then we can take $P = 0$ and there is nothing to show. Otherwise, fix $\alpha > 0$ (to be chosen later). Since $K$ is compact, we can take $\delta > 0$ with the property that, for $x_{1},x_{2} \in K$, 
\begin{align*}
    \vert x_{1} - x_{2} \vert \leq \delta \Rightarrow \vert \psi(x_1,\eta)-\psi(x_{2},\eta)\vert \leq \alpha \Vert \eta \Vert \hspace{1em} \text{for $\eta \in \subext{k}\R^{d}$}.
\end{align*}
Notice that, after possibly decomposing $T$, we can assume without loss of generality that $\diam(A) < \delta$. Take $x_{0} \in A$. Since $\psi$ is a flat integrand, there is $N \in \N$ and a collection $\{\xi_{i}\}_{i=1}^{N} \subset \subext{k}\R^{d}$ of simple $k$-vectors with
\begin{align*}
    \xi &= \sum_{i=1}^{N} \xi_{i}, & \sum_{i=1}^{N} \psi(x_{0},\xi_{i}) &\leq \psi(x_{0},\xi) + \alpha.
\end{align*}
For $x \in A$, we can then estimate
\begin{align*}
    \sum_{i=1}^{N} \psi(x,\xi_{i}) \leq \psi(x,\xi) + \alpha\bigg(1 + \Vert \xi \Vert + \sum_{i=1}^{N} \Vert \xi_{i} \Vert\bigg).
\end{align*}
Now notice that
\begin{align*}
    \sum_{i=1}^{N} \Vert \xi_{i} \Vert \leq C_{\psi,K}(C_{\psi,K}\cdot \Vert \xi \Vert + \alpha).
\end{align*}
where $C_{\psi,K} > 0$ is defined as in \eqref{psiVsMass}. We now choose $\alpha$ so that
\begin{align*}
    \alpha\bigg(1 + \Vert \xi \Vert + C_{\psi,K}(C_{\psi,K}\cdot \Vert \xi \Vert + \alpha)\bigg) \leq \frac{\varepsilon}{2 \cdot \mathscr{L}^{d}(K)}.
\end{align*}
It follows that, for $x \in A$,
\begin{align}\label{NormalPolyStrictIntermediate1}
    \sum_{i=1}^{N} \psi(x,\xi_{i}) \leq \psi(x,\xi) + \frac{\varepsilon}{2\cdot \mathscr{L}^{d}(K)}.
\end{align}
For $1 \leq i \leq N$, set $T_{i} = \xi_{i} \cdot \mathscr{L}^{d} \mvert A$. It follows from \eqref{NormalPolyStrictIntermediate1} that
\begin{align*}
    \sum_{i=1}^{N}\Mbf_{\psi}(T_{i}) \leq \Mbf_{\psi}(T) + \frac{\varepsilon}{2}.
\end{align*}
For $1 \leq i \leq N$, we appeal to Lemma \ref{lm:SimpleConstPolyStrict} together with Proposition \ref{pp:Simple} to find $P_{i} \in \Prm_{k}(\R^{d})$ with $\supp P_{i} \subset K$ such that
\begin{align*}
    \Fbf_{K}(P_{i}-T_{i}) &\leq \frac{\varepsilon}{2N}, & \Mbf_{\psi}(P_{i}) \leq \Mbf_{\psi}(T_{i}) + \frac{\varepsilon}{2N}.
\end{align*}
The result follows after setting
\begin{align*}
    P &:= \sum_{i=1}^{N} P_{i}. \qedhere
\end{align*}
\end{proof}

\begin{proof}[Proof of Proposition \ref{pp:NormalPolyStrict}] Take an intermediate compact set $A \subset \Int K$ with $\supp T \subset \Int A$. Let $\rho: \R^{d} \to [0,\infty)$ be smooth with $\supp \rho \subset B(0,1)$ and $\Vert \rho \Vert_{\Lrm^{1}} = 1$. For $\delta > 0$, define $\rho_{\delta}: \R^{d} \to \R$ by $\rho_{\delta}(x) := \delta^{-d}\rho(x/\delta)$ for $x \in \R^{d}$, and let $T_{\delta} := \rho_{\delta} * T$ as in \cite[4.1.2]{Federer69book}. Fix $\delta > 0$ with $(\supp T)_{\delta} \subset A$. Then (see \cite[4.1.18]{Federer69book}) we have $\supp T_{\delta} \subset A$ and
\begin{align}\label{NormalPolyStrict1}
    \Fbf_{A}(T_{\delta}-T) &\leq \delta \cdot \Nbf(T).
\end{align}
For $\omega \in \Dcal^{k}(\R^{d})$, we have
\begin{align*}
    T_{\delta}(\omega) = \int \bigg\langle \int \vec{T}(x)\rho_{\delta}(x-y) \dd \Vert T \Vert(x), \omega(y) \bigg\rangle \dd y.
\end{align*}
It follows that
\begin{align*}
    \Mbf_{\psi}(T_{\delta}) = \int \psi\bigg(y,\int\vec{T}(x)\rho_{\delta}(x-y) \dd \Vert T \Vert(x)\bigg) \dd y
\end{align*}
and so we can use Jensen's inequality to deduce that
\begin{align*}
    \Mbf_{\psi}(T_{\delta}) \leq \int \int \rho_{\delta}(x-y) \psi(y,\vec{T}(x)) \dd y \dd \Vert T \Vert(x).
\end{align*}
Take $(x,\xi) \in \supp T \times \subext{k}\R^{d}$. Since $\psi$ is continuous, we have
\begin{align}\label{NormalPolyStrict2}
    \int \rho_{\delta}(x-y)\psi(y,\xi) \dd y \to \psi(x,\xi) \hspace{1em} \text{as $\delta \downarrow 0$}.
\end{align}
Furthermore, notice that
\begin{align*}
    \bigg\vert \int \rho_{\delta}(x-y)\psi(y,\xi) \dd y \bigg\vert \leq C_{\psi,K} \Vert \xi \Vert
\end{align*}
where $C_{\psi,K} > 0$ is defined in \eqref{psiVsMass}. We then use \eqref{NormalPolyStrict2} with $\xi = \vec{T}(x)$ at $\Vert T \Vert$-a.e.\ $x \in \supp T$ to deduce that
\begin{align*}
    \int\int \rho_{\delta}(x-y)\psi(y,\vec{T}(x)) \dd y \dd\Vert T \Vert(x) \to \Mbf_{\psi}(T) \hspace{1em} \text{as $\delta \downarrow 0$}
\end{align*}
and so
\begin{align*}
    \limsup_{\delta \downarrow 0} \Mbf_{\psi}(T_{\delta}) \leq \Mbf_{\psi}(T).
\end{align*}
Since \eqref{NormalPolyStrict1} holds, we can use Lemma \ref{lm:psiLSC} to deduce that
\begin{align*}
    \Mbf_{\psi}(T) &\leq \liminf_{\delta \downarrow 0} \Mbf_{\psi}(T_{\delta}).
\end{align*}
Combining,
\begin{align}\label{NormalPolyStrict3}
    \Mbf_{\psi}(T_{\delta}) \to \Mbf_{\psi}(T) \hspace{1em} \text{as $\delta \downarrow 0$}.
\end{align}
Since \eqref{NormalPolyStrict1} and \eqref{NormalPolyStrict3} hold, we can take $\delta > 0$ so that $\supp T_{\delta} \subset A$ and
\begin{align}\label{NormalPolyStrict4}
    \Fbf_{A}(T_{\delta} -T) &\leq \frac{\varepsilon}{3}, & \Mbf_{\psi}(T_{\delta}) &\leq \Mbf_{\psi}(T) + \frac{\varepsilon}{3}.
\end{align}
It is shown in \cite[4.1.2]{Federer69book} that there is a smooth $\eta: \R^{d} \to \subext{k}\R^{d}$ with $\supp \eta \subset \Int K$ such that $T_{\delta} = \eta \cdot \mathscr{L}^{d}$. We can then find a step function $\xi: \R^{d} \to \subext{k}\R^{d}$ with $\supp \xi \subset \Int K$
\begin{align*}
    \Vert \xi - \eta \Vert_{\Lrm^{1}} \leq \frac{\varepsilon}{3} \cdot \min\bigg\{1,\frac{1}{C_{\psi,K}}\bigg\}.
\end{align*}
It follows that
\begin{align}\label{NormalPolyStrict5}
    \Mbf(\xi \cdot \mathscr{L}^{d} - T_{\delta}) &\leq \frac{\varepsilon}{3}, & \Mbf_{\psi}(\xi \cdot \mathscr{L}^{d} - T_{\delta}) \leq \frac{\varepsilon}{3}.
\end{align}
By appealing to linearity and Lemma \ref{lm:NormalPolyStrictIntermediate}, we can find $P \in \Prm_{k}(\R^{d})$ with $\supp P \subset K$ and 
\begin{align}\label{NormalPolyStrict6}
    \Fbf_{K}(P-\xi \cdot \mathscr{L}^{d}) &\leq \frac{\varepsilon}{3}, & \Mbf_{\psi}(P) &\leq \Mbf_{\psi}(\xi \cdot \mathscr{L}^{d}) + \frac{\varepsilon}{3}.
\end{align}
The result follows after combining \eqref{NormalPolyStrict4}, \eqref{NormalPolyStrict5}, and \eqref{NormalPolyStrict6}.
\end{proof}

\subsection{Integral Currents} Let $\psi: \R^{d} \times \subext{k}\R^{d} \to [0,\infty)$ be a positive convex integrand. In Proposition \ref{pp:NormalPolyStrict}, we demonstrated that a normal current $T \in \Nrm_{k}(\R^{d})$ can be $\psi$-strictly approximated by polyhedral chains. However, if we additionally assume that $T$ is integral, then we cannot expect that the proof of Proposition \ref{pp:NormalPolyStrict} will produce an integral approximant (for example, mollification destroys integrality). Nevertheless, by an entirely different procedure, one can $\psi$-strictly approximate integral currents by integral polyhedral chains. We stress that, because integral currents are simple, it suffices to consider strict (as opposed to $\psi$-strict) approximation results in this setting (see Proposition \ref{pp:Simple}).

Before proceeding, we fix our conventions concerning integral currents. We say that a current $T \in \Dcal_{k}(\R^{d})$ with $\Mbf(T) < \infty$ is \textbf{integer-rectifiable} if there is an $(\mathscr{H}^{k},k)$-rectifiable set (see \cite[3.2.14]{Federer69book}) $R \subset \R^{d}$ such that
\begin{enumerate}
    \item $\Vert T \Vert = m \cdot \mathscr{H}^{k} \mvert R$, where $m: \R^{d} \to \Z$ is $\mathscr{H}^{k} \mvert R$-integrable,
    \item $\mathrm{span}\hspace{0.2em} \vec{T}(x) = T_{x}R$ for $\mathscr{H}^{k} \mvert R$-a.e.\ $x \in \R^{d}$. 
\end{enumerate}
We shall use $\IRrm_{k}(\R^{d})$ to denote the space of integer-rectifiable $k$-currents in $\R^{d}$. One can show (see \cite[4.1.28]{Federer69book}) that $\IRrm_{k}(\R^{d})$ agrees with the space of \textit{rectifiable $k$-currents} defined in \cite[4.1.24]{Federer69book}. The space of \textbf{integral currents} is given by
\begin{align*}
    \Irm_{k}(\R^{d}):= \Nrm_{k}(\R^{d}) \cap \IRrm_{k}(\R^{d}).
\end{align*}

Notice that $\Irm_{0}(\R^{d}) = \IRrm_{0}(\R^{d})=\IPrm_{0}(\R^{d})$ and so any polyhedral approximation result is immediate in this case. We shall therefore assume that $k \geq 1$ throughout this subsection. We start with the following \textit{strong} approximation result:

\begin{lemma}\label{lm:RectPlaneStrong} Let $K \subset \R^{d}$ be compact and take $\varepsilon > 0$. Let $T \in \IRrm_{k}(\R^{d})$ and take a $k$-dimensional affine subspace $\Pi \subset \R^{d}$. Suppose that $\supp T \subset \Pi \cap \Int K$. Then there is $P \in \IPrm_{k}(\R^{d})$ with $\supp P \subset \Pi \cap K$ and 
\begin{align*}
    \Mbf(P-T) \leq \varepsilon.
\end{align*}
\end{lemma}

In analogy with the proof of Lemma \ref{lm:SimpleConstPolyStrict}, we shall prove Lemma \ref{lm:RectPlaneStrong} by first considering two cases of increasing generality:

\begin{lemma}\label{lm:RectPlaneStrong1} Let $K \subset \R^{k}$ be compact and take $\varepsilon > 0$. Let $A \subset \R^{k}$ be $\mathscr{L}^{k}$-measurable and suppose that $A \subset \Int K$. Fix $\alpha \in \R$ and set $T := \alpha(e_{1} \wedge \cdots \wedge e_{k}) \cdot \mathscr{L}^{k} \mvert A$. Then there is $P \in \Prm_{k}(\R^{k})$ with $\supp P \subset K$ and 
\begin{align*}
    \Mbf(P-T) \leq \varepsilon.
\end{align*}
Furthermore, if $\alpha \in \Z$, then we can additionally ensure that $P \in \IPrm_{k}(\R^{k})$.
\end{lemma}

\begin{proof} If $\alpha = 0$, then we can take $P := 0$ and there is nothing to show. Otherwise, we can use the fact that $A \subset \Int K$ is $\mathscr{L}^{k}$-measurable to find $N \in \N$ and a $\mathscr{L}^{k}$-disjoint collection $\{A_{i}\}_{i=1}^{N}$ of rectangles (i.e., each $A_{i}$ is a product of $k$ intervals) such that $A_{i} \subset K$ for $1 \leq i \leq N$ and 
\begin{align*}
    \mathscr{L}^{k}\bigg(A \triangle \bigcup_{i=1}^{N} A_{i}\bigg) \leq \frac{\varepsilon}{\vert \alpha \vert}.
\end{align*}
For $1 \leq i \leq N$, set $P_{i} := \alpha(e_{1} \wedge \cdots \wedge e_{k}) \cdot \mathscr{L}^{k} \mvert A_{i}$ and notice that $P_{i} \in \Prm_{k}(\R^{k})$. Let
\begin{align*}
    P := \sum_{i=1}^{N} P_{i}.
\end{align*}
If $\alpha \in \Z$, then notice that $P_{i} \in \IPrm_{k}(\R^{k})$ for $1 \leq i \leq N$. 
\end{proof}

\begin{lemma}\label{lm:RectPlaneStrong2} Let $K \subset \R^{k}$ be compact and take $\varepsilon > 0$. Let $T \in \IRrm_{k}(\R^{k})$ and suppose that $\supp T \subset \Int K$. Then there is $P \in \IPrm_{k}(\R^{k})$ with $\supp P \subset K$ and 
\begin{align*}
    \Mbf(P-T) \leq \varepsilon.
\end{align*}
\end{lemma} 

\begin{proof} Since $T \in \IRrm_{k}(\R^{k})$ and $\supp T \subset \Int K$, there is an $\mathscr{L}^{k}$-integrable $m: \R^{k} \to \Z$ with $\supp m \subset \Int K$ such that $T = m (e_{1} \wedge \cdots \wedge e_{k}) \cdot \mathscr{L}^{k}$. For $\alpha \in \Z$, define
\begin{align*}
    T_{\alpha} := \alpha (e_{1} \wedge \cdots \wedge e_{k}) \cdot \mathscr{L}^{k} \mvert m^{-1}(\{\alpha\}).
\end{align*}
Take $N \in \N$ with
\begin{align*}
    \Mbf\bigg(T- \sum_{\alpha=-N}^{N} T_{\alpha}\bigg) \leq \frac{\varepsilon}{2}.
\end{align*}
Set $P_{0} = 0$. Take $\alpha \in \Z\backslash\{0\}$ with $\vert \alpha \vert \leq N$. Notice that $m^{-1}(\{\alpha\}) \subset \Int K$ is $\mathscr{L}^{k}$-measurable, so we can appeal to Lemma \ref{lm:RectPlaneStrong1} to find $P_{\alpha} \in \IPrm_{k}(\R^{k})$ with $\supp P_{\alpha} \subset K$ and $\Mbf(P_{\alpha} - T_{\alpha}) \leq \varepsilon/4N$. The result follows after setting
\begin{align*}
    P &:= \sum_{\alpha=-N}^{N} P_{\alpha}. \qedhere
\end{align*}
\end{proof}

\begin{proof}[Proof of Lemma \ref{lm:RectPlaneStrong}] Let $\Phi: \Pi \to \R^{k}$ be an isometry. Notice that $\Phi_{*}T \in \IRrm_{k}(\R^{k})$ (see \cite[4.1.30]{Federer69book}) and $\supp \Phi_{*}T \subset \Int \Phi(\Pi \cap K)$. We can therefore appeal to Lemma \ref{lm:RectPlaneStrong2} to find $\tilde{P} \in \IPrm_{k}(\R^{k})$ with $\supp \tilde{P} \subset \Phi(\Pi \cap K)$ and $\Mbf(\tilde{P}-\Phi_{*}T) \leq \varepsilon$. Set $P = (\Phi^{-1})_{*}\tilde{P}$. Then $\supp P \subset K$ and the result follows since $\Lip(\Phi)=1$.
\end{proof}

Let $T \in \IRrm_{k}(\R^{d})$ and let $\Pi \subset \R^{d}$ be a $k$-dimensional affine subspace. We now work towards removing the assumption $\supp T \subset \Pi$ that is present in Lemma \ref{lm:RectPlaneStrong}:

\begin{proposition}\label{pp:RectStrong} Let $K \subset \R^{d}$ be compact and take $\varepsilon > 0$. Let $T \in \IRrm_{k}(\R^{d})$ and suppose that $\supp T \subset \Int K$. Then there is $P \in \IPrm_{k}(\R^{d})$ with $\supp P \subset K$ and a $\Crm^{1}$-diffeomorphism $f: \R^{d} \to \R^{d}$ with $\Lip(f) \leq 1+\varepsilon$ and $\Vert f - \mathrm{id}_{\R^{d}} \Vert_{\infty} \leq \varepsilon$ such that
\begin{align*}
    \Mbf(P-f_{*}T) \leq \varepsilon.
\end{align*}
\end{proposition}

A similar result can be found in \cite[Theorem 4.4]{DePauwVasilyev23}, but since we need to keep track of the support we give a complete proof. The proof of Proposition \ref{pp:RectStrong} below is based on \cite[4.2.19]{Federer69book}, but notice that we only assume that $T \in \IRrm_{k}(\R^{d})$ and we obtain a weaker conclusion that does not control the boundary. We shall need the following result, which is a rephrased version of \cite[3.1.23]{Federer69book}:

\begin{lemma}\label{lm:Flattening} Let $M \subset \R^{d}$ be a $k$-dimensional $\Crm^{1}$-submanifold. Fix $\alpha \in (0,1)$ and take $x_{0} \in M$. Then there is a $k$-dimensional affine subspace $\Pi \subset \R^{d}$ and an $R > 0$ with the following property: For $r \in (0,R)$, there is a $\Crm^{1}$-diffeomorphism $f: \R^{d} \to \R^{d}$ such that
\begin{align}\label{cc:Flattening1}
    \overline{B}(x_{0},\alpha r) \cap M &= \overline{B}(x_{0},\alpha r) \cap f^{-1}(\Pi), & f(x) &= x \hspace{1em} \text{for $x \in \R^{d}\backslash \overline{B}(x_{0},r)$},
\end{align}
and
\begin{align}\label{cc:Flattening2}
    \max\{\Lip(f),\Lip(f^{-1})\} &\leq \frac{1}{\alpha}, & \Vert f- \mathrm{id}_{\R^{d}} \Vert_{\infty} &\leq r.
\end{align}
\end{lemma}

\begin{proof} It is not restrictive to assume that $x_{0} = 0$. Since $M$ is a $k$-dimensional $\Crm^{1}$-submanifold, there is an open neighbourhood $U \subset \R^{d}$ of $0$ and an orthogonal projection $\pi: \R^{d} \to \R^{k}$ such that $\pi(U \cap M) = \pi(U)$, $\pi\vert_{U \cap M}$ is injective, and 
\begin{align*}
    g := (\pi\vert_{U \cap M})^{-1}: \pi(U) \to \R^{d}
\end{align*}
is $\Crm^{1}$ with $Dg(0) = \pi^{*}$. Notice that $g(0)=0$ since $0 \in M$ and $\pi(0)=0$. Let $\varphi: \R \to [0,1]$ be a smooth cut-off function with $\varphi(t) = 1$ for $t \leq \alpha$ and $\varphi(t) = 0$ for $t \geq 1$. Take $R > 0$ (to be chosen later) with $\overline{B}(0,R) \subset U$.

Fix $r \in (0,R)$ and define $f: \R^{d} \to \R^{d}$ where, for $x \in \R^{d}$,
\begin{align}\label{Flattening1}
    f(x) = \left\{\begin{array}{cl}
        x + \varphi(\vert x \vert/r) \cdot(\pi^{*}(\pi(x))-g(\pi(x))) & \text{if $\vert x \vert \leq r$}, \\
        x & \text{if $\vert x \vert > r$}.
    \end{array}\right.
\end{align}
One can check that $f$ is $\Crm^{1}$. Set $\Pi = \pi^{*}(\R^{k})$. 

We now verify \eqref{cc:Flattening1}. First notice that \eqref{Flattening1} ensures that $f(x) = x$ for $x \in \R^{d}\backslash \overline{B}(0,r)$. Take $x \in \overline{B}(0,\alpha r)$. Then $\varphi(\vert x \vert/r)=1$ and so 
\begin{align}\label{Flattening2}
    f(x) - \pi^{*}(\pi(x)) = x - g(\pi(x)).
\end{align}
If $x \in M$, then $x \in U \cap M$ and so $g(\pi(x)) = x$; it follows from \eqref{Flattening2} that $f(x) = \pi^{*}(\pi(x)) \in \Pi$. Conversely, if $f(x) \in \Pi$, then we can use \eqref{Flattening2} to deduce that $x-g(\pi(x)) \in \Pi$; on the other hand, 
\begin{align*}
    \pi(x-g(\pi(x))) = \pi(x)-\pi(x) = 0,
\end{align*}
so we must have $x-g(\pi(x)) = 0$ and thus $x = g(\pi(x)) \in M$. 

We shall now work towards showing that \eqref{cc:Flattening2} holds after choosing $R > 0$. Since $g(0) = 0$ and $Dg(0) = \pi^{*}$, we have
\begin{align*}
    \lim_{v \to 0}\frac{g(v)-\pi^{*}(v)}{\vert v \vert} = 0.
\end{align*}
Furthermore, since $g$ is $\Crm^{1}$, we have 
\begin{align*}
    \lim_{v \to 0} \Vert Dg(v)-\pi^{*}\Vert = 0.
\end{align*}
Take $\delta > 0$ (to be chosen later). We can then adjust $R > 0$ so that
\begin{align}\label{Flattening3}
    \vert g(v) - \pi^{*}(v) \vert &\leq \delta \vert v \vert, & \Vert Dg(v) - \pi^{*} \Vert &\leq \delta
\end{align}
for $v \in \R^{k}$ with $\vert v \vert \leq R$. Fix $x \in \R^{d}$ with $\vert x \vert \leq r$ (recall that we fixed $r \in (0,R)$). For $w \in \R^{d}$, we can use \eqref{Flattening1} to compute
\begin{align*}
    Df(x)(w) - w &= \varphi'\bigg(\frac{\vert x \vert}{r}\bigg)\bigg(\frac{x}{r\vert x \vert} \cdot w\bigg)(\pi^{*}(\pi(x))-g(\pi(x))) \\
    &\qquad + \varphi\bigg(\frac{\vert x \vert}{r}\bigg)(\pi^{*}(\pi(w))-Dg(\pi(x))(\pi(w))).
\end{align*}
Since $\vert x \vert \leq r \leq R$, we can use \eqref{Flattening3} to deduce that, for $w \in \R^{d}$,
\begin{align*}
    \vert Df(x)(w) - w \vert \leq \Lip(\varphi) \cdot \frac{\vert w \vert}{r} \cdot \delta\vert \pi(x) \vert + \delta\vert\pi(w)\vert \leq \delta( \Lip(\varphi) + 1)\vert w \vert.
\end{align*}
We now choose $\delta > 0$ so that
\begin{align*}
    \delta(\Lip(\varphi)+1) \leq 1-\alpha.
\end{align*}
Then $\Vert Df(x)-\mathrm{id}_{\R^{d}} \Vert \leq 1-\alpha$ for $x \in \R^{d}$ with $\vert x \vert \leq r$. Since $f(x)=x$ for $x \in \R^{d}$ with $\vert x \vert > r$, we have 
\begin{align}\label{Flattening4}
    \Vert Df(x)-\mathrm{id}_{\R^{d}} \Vert \leq 1-\alpha \hspace{1em} \text{for $x \in \R^{d}$}.
\end{align}
It follows that $\Lip(f-\mathrm{id}_{\R^{d}}) \leq 1 - \alpha$. In particular, for $x_{1},x_{2} \in \R^{d}$, we have
\begin{align*}
    \vert f(x_{2})-f(x_{1})-(x_{2}-x_{1})\vert \leq (1-\alpha)\vert x_{2}-x_{1}\vert \leq (\alpha^{-1}-1)\vert x_{2}-x_{1}\vert
\end{align*}
and so
\begin{align}\label{Flattening5}
    \alpha\vert x_{2} - x_{1} \vert \leq \vert f(x_{2})-f(x_{1}) \vert \leq \alpha^{-1}\vert x_{2}-x_{1} \vert.
\end{align}
It follows that $\max\{\Lip(f),\Lip(f^{-1})\} \leq \alpha^{-1}$. For the remaining conclusion of \eqref{cc:Flattening2}, take $x \in \R^{d}$. If $\vert x \vert \geq r$, then $f(x)=x$ and there is nothing to show. Otherwise, we can use \eqref{Flattening1}, \eqref{Flattening3}, and the fact that $\delta \leq 1$ to estimate
\begin{align*}
    \vert f(x) - x \vert \leq \vert \pi^{*}(\pi(x))-g(\pi(x))\vert \leq \delta \vert \pi(x) \vert \leq r.
\end{align*}

Finally, we have seen that $f$ is $\Crm^{1}$, and it follows from \eqref{Flattening5} that $f$ is a bijection. Since $Df(x)$ is invertible for all $x \in \R^{d}$ (see \eqref{Flattening4}), $f$ is a $\Crm^{1}$-diffeomorphism. 
\end{proof}

Take $\mu \in \Mcal^{+}(\R^{d})$. Before proceeding, we recall that, for $x \in \R^{d}$, the \textbf{$k$-density} of $\mu$ at $x$ is given by
\begin{align}
    \Theta^{k}(\mu,x) := \lim_{r \downarrow 0} \frac{\mu(\overline{B}(x,r))}{\omega_{k}r^{k}}
\end{align}
whenever this limit exists.

\begin{proof}[Proof of Proposition \ref{pp:RectStrong}] Since $T \in \IRrm_{k}(\R^{d})$, we can find a countable collection $\{M_{i}\}_{i=1}^{\infty}$ of $k$-dimensional $\Crm^{1}$-submanifolds in $\R^{d}$ and an $N \subset \R^{d}$ with $\Vert T \Vert(N) = 0$ such that, for $x \in \R^{d}\backslash N$, there is $i(x) \in \N$ with $x \in M_{i(x)}$ and
\begin{align}\label{RectStrong1}
    \Theta^{k}(\Vert T \Vert \mvert (\R^{d}\backslash M_{i(x)}),x) = 0,
\end{align}
We can furthermore assume without loss of generality that $\R^{d}\backslash N \subset \supp T$ and that, for $x \in \R^{d}\backslash N$, we have
\begin{align}\label{RectStrong2}
    \Theta^{k}(\Vert T \Vert,x) \in (0,\infty).
\end{align}
Fix $\alpha \in (0,1)$ (to be chosen later).  Take $x \in \R^{d}\backslash N$. It follows from \eqref{RectStrong1} and \eqref{RectStrong2} that
\begin{align*}
    \lim_{r \downarrow 0} \frac{\Vert T \Vert(\overline{B}(x,\alpha r))}{\Vert T \Vert(\overline{B}(x,r))} &= \alpha^{k}, & \lim_{r \downarrow 0} \frac{\Vert T \Vert \mvert (\R^{d} \backslash M_{i(x)})(\overline{B}(x,\alpha r))}{\Vert T \Vert(\overline{B}(x,\alpha r))} &= 0,
\end{align*}
and so
\begin{align*}
    \lim_{r \to 0} \frac{\Vert T \Vert(\overline{B}(x,r)\backslash( M_{i(x)} \cap \overline{B}(x,\alpha r)))}{\Vert T \Vert(\overline{B}(x,r))} = \lim_{r \downarrow 0} \frac{\Vert T \Vert(\overline{B}(x,r)\backslash \overline{B}(x,\alpha r))}{\Vert T \Vert(\overline{B}(x,r))} =1 - \alpha^{k}.
\end{align*}
We can therefore find $R_{1}(x) > 0$ such that, for $r \in (0,R_{1}(x))$, 
\begin{align*}
    \Vert T \Vert(\overline{B}(x,r)\backslash( M_{i(x)} \cap \overline{B}(x,\alpha r))) < 2(1-\alpha^{k})\Vert T \Vert(\overline{B}(x,r)).
\end{align*}
Let $\Pi(x) \subset \R^{d}$ and $R_{2}(x) > 0$ be the $k$-dimensional affine subspace and the constant from Lemma \ref{lm:Flattening} (applied with $M_{i(x)}$, $\alpha$, and $x$) respectively. Since $\supp T \subset \Int K$, there is $\delta > 0$ with $(\supp T)_{\delta} \subset \Int K$. Set
\begin{align*}
    R(x) := \min\{R_{1}(x),R_{2}(x),\delta/2,1-\alpha^{k}\}
\end{align*}
and notice that the following cover of $\R^{d}\backslash N$ is a fine cover:
\begin{align*}
    \bigl\{ \overline{B}(x,r): x \in \R^{d}\backslash N, \hspace{0.2em} 0 < r < R(x) \bigr\}.
\end{align*}
We can therefore appeal to the Vitali-Besicovitch covering theorem (see, for example, \cite[Theorem 2.19]{AmbrosioFuscoPallara00book}) to find a collection $\{\overline{B}(x_{j},r_{j})\}_{j=1}^{\infty}$ of disjoint closed balls satisfying $r_{j} < R(x_{j})$ for $j \in \N$ and
\begin{align}\label{RectStrong3}
    \Vert T \Vert\bigg(\R^{d}\backslash \bigg(\bigcup_{j=1}^{\infty} \overline{B}(x_{j},r_{j})\bigg)\bigg) = 0.
\end{align}

Fix $j \in \N$. Since $r_{j} < R_{2}(x_{j})$ and $R_{2}(x_{j})$ is given by Lemma \ref{lm:Flattening}, there is a $\Crm^{1}$-diffeomorphism $f_{j}: \R^{d} \to \R^{d}$ such that 
\begin{align*}
    \overline{B}(x_{j},\alpha r_{j}) \cap M_{i(x_{j})} &= \overline{B}(x_{j},\alpha r_{j}) \cap f_{j}^{-1}(\Pi(x_{j})), & f_{j}(x) &= x \hspace{1em} \text{for $x \in \R^{d}\backslash \overline{B}(x_{j},r_{j})$},
\end{align*}
and
\begin{align*}
    \max\{\Lip(f_{j}),\Lip(f_{j}^{-1})\} &\leq \frac{1}{\alpha}, & \Vert f_{j} - \mathrm{id}_{\R^{d}} \Vert_{\infty} &\leq r_{j}.
\end{align*}
Since $ \overline{B}(x_{j},\alpha r_{j}) \cap M_{i(x_{j})} = \overline{B}(x_{j},\alpha r_{j}) \cap f_{j}^{-1}(\Pi(x_{j}))$, we have
\begin{align*}
    \supp (f_{j})_{*}(T \mvert (\overline{B}(x_{j},\alpha r_{j}) \cap M_{i(x_{j})})) &\subset f_{j}(\overline{B}(x_{j},\alpha r_{j}) \cap f_{j}^{-1}(\Pi(x_{j}))) \\
    &\subset f_{j}(\overline{B}(x_{j},\alpha r_{j})) \cap \Pi(x_{j}).
\end{align*}
Since $\Vert f_{j} - \mathrm{id}_{\R^{d}} \Vert_{\infty} \leq r_{j}$, we have $f_{j}(\overline{B}(x_{j},\alpha r_{j})) \subset \overline{B}(x_{j},2r_{j})$. Since $r_{j} \leq \delta/2$ and $x_{j} \in \supp T$, we deduce that $\overline{B}(x_{j},2r_{j}) \subset \Int K$ and so
\begin{align*}
    \supp (f_{j})_{*}(T \mvert (\overline{B}(x_{j},\alpha r_{j}) \cap M_{i(x_{j})})) \subset \Pi(x_{j}) \cap \Int K.
\end{align*}
Since $r_{j} \leq R_{1}(x_{j})$, we can take $\varepsilon_{j}$ satisfying 
\begin{align}\label{RectStrong4}
     \alpha^{-k}\cdot\Vert T \Vert(\overline{B}(x_{j}, r_{j})\backslash(M_{i(x_{j})} \cap \overline{B}(x_{j},\alpha r_{j}))) + \varepsilon_{j} \leq 2(\alpha^{-k}-1)\Vert T \Vert(\overline{B}(x_{j},r_{j})).
\end{align}
We now apply Lemma \ref{lm:RectPlaneStrong} to $(f_{j})_{*}(T \mvert (M_{i(x_{j})} \cap \overline{B}(x_{j},\alpha r_{j})))$ to obtain $P_{j} \in \IPrm_{k}(\R^{d})$ with $\supp P_{j} \subset \Pi(x_{j}) \cap K$ and
\begin{align}\label{RectStrong5}
    \Mbf(P_{j} - (f_{j})_{*}(T \mvert (M_{i(x_{j})} \cap \overline{B}(x_{j},\alpha r_{j})))) \leq \varepsilon_{j}.
\end{align}
We can then use \eqref{RectStrong4}, \eqref{RectStrong5}, and the fact that $\Lip(f_{j}) \leq \alpha^{-1}$ to estimate
\begin{equation}\label{RectStrong6}
\begin{aligned}
    \Mbf(P_{j} - (f_{j})_{*}&(T \mvert \overline{B}(x_{j},r_{j}))) \\
    &\leq \varepsilon_{j} + \alpha^{-k} \cdot \Vert T \Vert(\overline{B}(x_{j}, r_{j})\backslash(M_{i(x_{j})} \cap \overline{B}(x_{j},\alpha r_{j}))) \\
    &\leq 2(\alpha^{-k}-1)\Vert T \Vert(\overline{B}(x_{j},r_{j})).
\end{aligned}
\end{equation}
Since the balls $\{\overline{B}(x_{j},r_{j})\}_{j=1}^{\infty}$ are disjoint and $\Vert T \Vert(\R^{d}) < \infty$, we can find $N \in \N$ with
\begin{align}\label{RectStrong7}
    \sum_{j=N+1}^{\infty} \Vert T \Vert(\overline{B}(x_{j},r_{j})) \leq (\alpha^{-k}-1)\Mbf(T).
\end{align}
Set
\begin{align*}
    P := \sum_{j=1}^{N} P_{j}
\end{align*}
and notice that $\supp P \subset K$. Let $f = f_{1} \circ \cdots \circ f_{N}$. Then $f = f_{j}$ on $\overline{B}(x_{j},r_{j})$ for $1 \leq j \leq N$ and 
\begin{align}\label{RectStrong8}
    f(x) = x \hspace{1em} \text{for $x \in \R^{d}\backslash\bigg(\bigcup_{j=1}^{N} \overline{B}(x_{j},r_{j})\bigg)$}.
\end{align}
Since \eqref{RectStrong3} holds, we have
\begin{align*}
    f_{*}T = \sum_{j=1}^{N} (f_{j})_{*}(T \mvert \overline{B}(x_{j},r_{j})) + \sum_{j=N+1}^{\infty} T \mvert \overline{B}(x_{j},r_{j}).
\end{align*}
We can then use \eqref{RectStrong6} and \eqref{RectStrong7} to obtain 
\begin{equation}\label{RectStrong9}
\begin{aligned}
    \Mbf(P-f_{*}T) &\leq 2(\alpha^{-k}-1)\Vert T \Vert\bigg(\bigcup_{j=1}^{N} \overline{B}(x_{j},r_{j})\bigg) + (\alpha^{-k}-1)\Mbf(T) \\
    & \leq 3(\alpha^{-k}-1)\Mbf(T).
\end{aligned}
\end{equation}
Since the balls $\{\overline{B}(x_{j},r_{j})\}_{j=1}^{\infty}$ are disjoint, we have
\begin{align*}
    \Vert f - \mathrm{id}_{\R^{d}} \Vert_{\infty} \leq \max_{1 \leq j \leq N} r_{j} \leq 1 - \alpha^{k}.
\end{align*}

We shall now show that $\Lip(f) \leq \alpha^{-1}$. Take $x_{1},x_{2} \in \R^{d}$. Since \eqref{RectStrong8} holds, we can assume that there are $j^{(1)},j^{(2)} \in \{1,\dotsc,N\}$ with
\begin{align*}
    x_{i} \in \overline{B}(x_{j^{(i)}},r_{j^{(i)}}) \hspace{1em} \text{for $i \in \{1,2\}$}.
\end{align*}
Since $\Lip(f_{j}) \leq \alpha^{-1}$ for $1 \leq j \leq N$, we can assume that $j^{(1)} \neq j^{(2)}$. Since the balls $\{\overline{B}(x_{j},r_{j})\}_{j=1}^{N}$ are disjoint, we can find $y \in \R^{d}$ with
\begin{align*}
    \vert x_{2} - x_{1} \vert &= \vert x_{2} - y \vert + \vert y - x_{1} \vert, & y &\not\in \bigcup_{j=1}^{N} \overline{B}(x_{j},r_{j}).
\end{align*}
We can then use \eqref{RectStrong8} and the fact that $\Lip(f_{j}) \leq \alpha^{-1}$ for $1 \leq j \leq N$ to deduce that
\begin{align*}
    \vert f(x_{2}) - f(x_{1}) \vert &\leq \vert f(x_{2})-y \vert + \vert y - f(x_{1}) \vert \\
    &\leq \alpha^{-1}\vert x_{2} - y \vert + \alpha^{-1}\vert y - x_{1} \vert \\
    &= \alpha^{-1}\vert x_{2} - x_{1} \vert.
\end{align*}
The result follows after choosing $\alpha$ so that
\begin{align*}
    3(\alpha^{-k}-1)\Mbf(T) &\leq \varepsilon, & 1-\alpha^{k} &\leq \varepsilon, & \alpha^{-1} &\leq 1+ \varepsilon. \qedhere
\end{align*}
\end{proof}

Let $K \subset \R^{d}$ be compact and take $T \in \Irm_{k}(K)$. In what follows, the \textbf{(integral) flat norm} of $T$ (on $K$) is given by
\begin{align*}
    \FI_{K}(T) := \inf \bigl\{\Mbf(T-\partial R) + \Mbf(R): R \in \Irm_{k+1}(K) \bigr\}.
\end{align*}
With Proposition \ref{pp:RectStrong} in hand, we are now ready to state and prove a strict approximation theorem for integral currents:

\begin{proposition}\label{pp:IntegralPolyStrict} Let $K \subset \R^{d}$ be compact and take $\varepsilon > 0$. Let $T \in \Irm_{k}(\R^{d})$ and suppose that $\supp T \subset \Int K$. Then there is $P \in \IPrm_{k}(\R^{d})$ with $\supp P \subset K$ and 
\begin{align*}
    \FI_{K}(P-T) &\leq \varepsilon, & \Mbf(P) &\leq \Mbf(T) + \varepsilon.
\end{align*}
\end{proposition}

\begin{proof} Fix $\alpha > 0$ (to be chosen later) with $\alpha \leq \mathrm{dist}(\supp T; \R^{d}\backslash \Int K)$. We start by applying Proposition \ref{pp:RectStrong} to find $P \in \IPrm_{k}(\R^{d})$ with $\supp P \subset K$ and a  $\Crm^{1}$-diffeomorphism $f: \R^{d} \to \R^{d}$ with $\Lip(f) \leq 1 + \alpha$ and $\Vert f - \mathrm{id}_{\R^{d}} \Vert_{\infty} \leq \alpha$ such that
\begin{align}\label{IntegralPolyStrict1}
    \Mbf(P-f_{*}T) &\leq \alpha.
\end{align}
It follows that
\begin{align}\label{IntegralPolyStrict2}
    \FI_{K}(f_{*}T - T) &\leq \alpha(1+\alpha)^{k}\cdot \Nbf(T), & \Mbf(f_{*}T) &\leq (1+\alpha)^{k}\cdot\Mbf(T). 
\end{align}
Combining \eqref{IntegralPolyStrict1} and \eqref{IntegralPolyStrict2} yields
\begin{align*}
    \FI_{K}(P-T) &\leq \alpha + \alpha(1+\alpha)^{k}\cdot \Nbf(T), & \Mbf(P) &\leq \alpha + (1+\alpha)^{k}\Mbf(T).
\end{align*}
The result follows after choosing $\alpha$ so that
\begin{align*}
    \alpha + \alpha(1+\alpha)^{k}\cdot \Nbf(T) &\leq \varepsilon, & \alpha + (1+\alpha)^{k}\Mbf(T) &\leq \Mbf(T) + \varepsilon. \qedhere
\end{align*}
\end{proof}

We end this section with the following corollary, which is an immediate consequence of Proposition \ref{pp:IntegralPolyStrict} and Proposition \ref{pp:Simple}:

\begin{corollary}\label{cr:IntegralPolyStrict} Let $\psi: \R^{d} \times \subext{k}\R^{d} \to [0,\infty)$ be a positive convex integrand. Let $K \subset \R^{d}$ be compact and take $\varepsilon > 0$. Let $T \in \Irm_{k}(\R^{d})$ and suppose that $\supp T \subset \Int K$. Then there is $P \in \IPrm_{k}(\R^{d})$ with $\supp P \subset K$ and 
\begin{align*}
    \FI_{K}(P-T) &\leq \varepsilon, & \Mbf_{\psi}(P) &\leq \Mbf_{\psi}(T) + \varepsilon.
\end{align*}
\end{corollary}

In the sequel, we shall restrict our attention to flat integrands, but this is not restrictive in the integral setting. For example, if $\psi: \R^{d} \times \subext{k}\R^{d} \to [0,\infty)$ is a positive convex integrand, then we can use \eqref{FlatSimple} to deduce that we can work with the flat envelope $\psi^{\flat}$ without loss of generality.

\subsection{Currents with Polyhedral Boundaries} Let $\psi: \R^{d} \times \subext{k}\R^{d} \to [0,\infty)$ be a flat integrand. We shall now demonstrate that a normal current $T \in \Nrm_{k}(\R^{d})$ with $\partial T \in \Prm_{k-1}(\R^{d})$ can be $\psi$-strictly approximated by polyhedral chains with the \textit{same} boundary. To do so, we adopt the approach of \cite[Section 1.3.4, Theorem 6]{GiaquintaModicaSoucek98book2}; we also show that a similar approach yields the analogous result for integral currents. The proof actually demonstrates that, in the case that $\partial T = 0$, the approximant is \textit{homologous} to the original current and so, for this reason, we rephrase the result here as in \cite[Proposition 2.3]{BonicattoRindlerTurnbull26}, where it is stated without proof. 

\begin{proposition}\label{pp:PolyBdryPolyStrict} Let $\psi: \R^{d} \times \subext{k}\R^{d} \to [0,\infty)$ be a flat integrand. Let $K \subset \R^{d}$ be compact and take $\varepsilon > 0$. Let $T \in \Nrm_{k}(\R^{d})$. Suppose that $\supp T \subset \Int K$ and that $\partial T \in \Prm_{k-1}(\R^{d})$. Then there are $P \in \Prm_{k}(\R^{d})$ and $R \in \Nrm_{k+1}(\R^{d})$ with $\supp P$, $\supp R \subset K$ and 
\begin{align*}
    T &= P + \partial R, & \Mbf_{\psi}(P) &\leq \Mbf_{\psi}(T) + \varepsilon, & \Mbf(R) &\leq \varepsilon.
\end{align*}
If $T \in \Irm_{k}(\R^{d})$, then we can additionally ensure that $P \in \IPrm_{k}(\R^{d})$ and that $R \in \Irm_{k+1}(\R^{d})$. 
\end{proposition}

\begin{proof} Fix $\alpha > 0$ (to be chosen later). Take an intermediate compact set $A \subset \Int K$ with $\supp T \subset \Int A$, and apply Proposition \ref{pp:NormalPolyStrict} to $T$ on $A$ to obtain $\tilde{P} \in \Prm_{k}(\R^{d})$ with $\supp \tilde{P} \subset A$ and
\begin{align}\label{PolyBdryPolyStrict1}
    \Fbf_{A}(\tilde{P}-T) &\leq \alpha, & \Mbf_{\psi}(\tilde{P}) &\leq \Mbf_{\psi}(T) + \alpha.
\end{align}
It follows that there are $X \in \Nrm_{k+1}(A)$ and $Y \in \Nrm_{k}(A)$ with 
\begin{align*}
    \tilde{P} - T &= \partial X + Y, & \Mbf(X)+\Mbf(Y) &\leq \alpha.
\end{align*}
Notice that $\partial Y \in \Prm_{k-1}(\R^{d})$. Take $\delta \in (0,1)$ (to be chosen later) with $A_{2d\delta} \subset K$ and apply the deformation theorem to $Y$, noting that we have conclusion (5), to obtain $P^{Y} \in \Prm_{k}(\R^{d})$ and $R^{Y} \in \Nrm_{k+1}(\R^{d})$ with $\supp P^{Y}$, $\supp R^{Y} \subset K$,
\begin{align*}
    Y &= P^{Y} + \partial R^{Y},
\end{align*}
and
\begin{align}\label{PolyBdryPolyStrict2}
    \Mbf(P^{Y}) &\leq \gamma(\Mbf(Y) + 2\delta \cdot \Mbf(\partial Y)), & \Mbf(R^{Y}) &\leq \gamma\delta \cdot \Mbf(Y)
\end{align}
Then
\begin{align*}
    \tilde{P}-T = P^{Y} + \partial(R^{Y} + X).
\end{align*}
Set $P := \tilde{P} - P^{Y}$ and $R := -(R^{Y}+X)$. Then $T = P + \partial R$ and \eqref{PolyBdryPolyStrict1}, \eqref{PolyBdryPolyStrict2}, \eqref{psiMassVsMass} yield
\begin{align*}
    \Mbf_{\psi}(P) &\leq \Mbf_{\psi}(T) + \alpha + C_{\psi,K} \cdot \gamma(\Mbf(Y) + 2\delta \cdot \Mbf(\partial Y)),
\end{align*}
where $C_{\psi,K}$ is defined as in \eqref{psiVsMass}, and
\begin{align*}
    \Mbf(R) &\leq \gamma\delta \cdot \Mbf(Y) + \Mbf(X)
\end{align*}
Since $\Mbf(X) + \Mbf(Y) \leq \alpha$, we have
\begin{align*}
    \Mbf_{\psi}(P) &\leq \Mbf_{\psi}(T) + \alpha(C_{\psi,K}\cdot\gamma +1) + 2C_{\psi,K}\cdot\gamma\delta \cdot \Mbf(\partial Y)
\end{align*}
and
\begin{align*}
    \Mbf(R) &\leq \alpha(\gamma\delta + 1).
\end{align*}
Now choose $\alpha $ and $\delta$ so that
\begin{align*}
    \alpha &\leq \min\bigg\{\frac{\varepsilon}{(\gamma+1)},\frac{\varepsilon}{2(C_{\psi,K}\cdot\gamma + 1)}\bigg\}, & 4C_{\psi,K}\cdot \gamma\delta \cdot\Mbf(\partial Y) &\leq \varepsilon.
\end{align*}
If $T \in \Irm_{k}(\R^{d})$, then we work as above, applying Corollary \ref{cr:IntegralPolyStrict} instead of Proposition \ref{pp:NormalPolyStrict} and recalling conclusion (4) of the deformation theorem.
\end{proof}

\subsection{Boundary Strict Approximation} We are now in a position to establish the boundary strict polyhedral approximation result for normal and integral currents discussed in the introduction, where a simplified version of the result was presented as Theorem \ref{tm:AbbrvPolyBdryStrict}. Recall that boundary strict approximation results for normal and integral currents are classical (see, for example, \cite[4.2.21, 4.2.24]{Federer69book}), but that the novel feature of our result is that we establish control on the support of the boundary of the approximants. 

Let $K \subset \R^{d}$ be compact and take $T \in \Nrm_{k}(\R^{d})$. Recall from \cite{BonicattoRindlerTurnbull26} that the \textbf{(normal) homogeneous flat norm} of $T$ (on $K$) is given by
\begin{align*}
    \FN^{\circ}_{K}(T) := \inf \bigl\{\Mbf(R) : R \in \Nrm_{k+1}(K), \hspace{0.2em} \partial R = T \bigr\}.
\end{align*}
If $T \in \Irm_{k}(\R^{d})$, then we also define the \textbf{(integral) homogeneous flat norm} of $T$ (on $K$) as follows:
\begin{align*}
    \FI^{\circ}_{K}(T) := \inf \bigl\{\Mbf(R) : R \in \Irm_{k+1}(K), \hspace{0.2em} \partial R = T \bigr\}.
\end{align*}

\begin{theorem}\label{tm:PolyBdryStrict} Take flat integrands
\begin{align*}
    \psi^{1} &: \R^{d} \times \subext{k}\R^{d} \to [0,\infty), & \psi^{2} &: \R^{d} \times \subext{k-1}\R^{d} \to [0,\infty).
\end{align*}
Let $K_{1},K_{2} \subset \R^{d}$ be compact with $K_{2} \subset K_{1}$ and take $\varepsilon > 0$. Let $T \in \Nrm_{k}(\R^{d})$. Suppose that $\supp T \subset \Int K_{1}$ and that $\supp \partial T \subset \Int K_{2}$. Then there is $P \in \Prm_{k}(\R^{d})$ such that $\supp P \subset K_{1}$, $\supp \partial P \subset K_{2}$, and the following conclusions hold:
\begin{enumerate}
    \item We have 
    \begin{align*}
        \FN_{K_{1}}(P-T) &\leq \varepsilon, & \Mbf_{\psi^{1}}(P) &\leq \Mbf_{\psi^{1}}(T) + \varepsilon.
    \end{align*}
    \item We have
    \begin{align*}
        \FN^{\circ}_{K_{2}}(\partial P - \partial T) &\leq \varepsilon, & \Mbf_{\psi^{2}}(\partial P) &\leq \Mbf_{\psi^{2}}(\partial T) + \varepsilon.
    \end{align*}
    \item If $\partial T \in \Irm_{k-1}(\R^{d})$, then we can additionally ensure that $\partial P \in \IPrm_{k-1}(\R^{d})$ and
    \begin{align*}
        \FI^{\circ}_{K_{2}}(\partial P - \partial T) &\leq \varepsilon
    \end{align*}
    \item If $T \in \Irm_{k}(\R^{d})$, then we can additionally ensure that $P \in \IPrm_{k}(\R^{d})$ and
    \begin{align*}
        \FI_{K_{1}}(P-T) &\leq \varepsilon
    \end{align*}
\end{enumerate}
\end{theorem}

To obtain Theorem \ref{tm:AbbrvPolyBdryStrict} from Theorem \ref{tm:PolyBdryStrict}, it suffices to take $\psi^{1}$ and $\psi^{2}$ to be the mass norms on $\subext{k}\R^{d}$ and $\subext{k-1}\R^{d}$ respectively.

\begin{proof} Fix $\alpha > 0$ (to be chosen later). Take an intermediate compact set $A_{1} \subset \Int K_{1}$ with $\supp T \subset \Int A_{1}$ and apply Proposition \ref{pp:NormalPolyStrict} to $T$ on $A_{1}$ to obtain $P_{1} \in \Prm_{k}(\R^{d})$ with $\supp P_{1} \subset A_{1}$ and 
\begin{align*}
    \FN_{A_{1}}(P_{1}-T) &\leq \alpha, & \Mbf_{\psi^{1}}(P_{1}) &\leq \Mbf_{\psi^{1}}(T) + \alpha.
\end{align*}
It follows that there are $X \in \Nrm_{k+1}(A_{1})$ and $Y \in \Nrm_{k}(A_{1})$ with 
\begin{align*}
    P_{1} - T &= \partial X + Y, & \Mbf(X)+\Mbf(Y) &\leq \alpha.
\end{align*}
Now take an intermediate compact set $A_{2} \subset \Int K_{2}$ with $\supp \partial T \subset \Int A_{2}$ and apply Proposition \ref{pp:PolyBdryPolyStrict} to $\partial T$ on $A_{2}$ to obtain $P_{2} \in \Prm_{k-1}(\R^{d})$ and $R \in \Nrm_{k}(\R^{d})$ with $\supp P_{2}$, $\supp R \subset A_{2}$ and 
\begin{align}\label{PolyBdryStrict1}
    \partial T &= P_{2} + \partial R, & \Mbf_{\psi^{2}}(P_{2}) &\leq \Mbf_{\psi^{2}}(\partial T) + \alpha, & \Mbf(R) &\leq \alpha.
\end{align}
Notice that
\begin{align*}
    \partial P_{1} - P_{2} = \partial(Y + R)
\end{align*}
Since $\partial(Y+R) \in \Prm_{k-1}(\R^{d})$ and $\supp(Y + R) \subset A_{1} \cup A_{2} \subset \Int K_{1}$, we can apply Proposition \ref{pp:PolyBdryPolyStrict} once more to obtain $\tilde{P} \in \Prm_{k}(\R^{d})$ with $\supp \tilde{P} \subset K_{1}$ and 
\begin{align*}
    \partial \tilde{P} &= \partial(Y+R), & \Mbf(\tilde{P}) &\leq \Mbf(Y+R) + \alpha.
\end{align*}
It follows that $\Mbf(\tilde{P}) \leq 3\alpha$. Set $P := P_{1} - \tilde{P}$ and observe that
\begin{align}\label{PolyBdryStrict2}
    \Fbf_{K_{1}}(P-T) &\leq 4\alpha, & \Mbf_{\psi^{1}}(P) &\leq \Mbf_{\psi^{1}}(T) + (1+3C_{\psi^{1},K_{1}})\alpha
\end{align}
where $C_{\psi^{1},K_{1}} > 0$ is defined as in \eqref{psiVsMass}. Now choose $\alpha$ so that 
\begin{align*}
    \alpha \leq \varepsilon\cdot\min\bigg\{\frac{1}{4},\frac{1}{1+3C_{\psi^{1},K_{1}}}\bigg\}.
\end{align*}
Notice that conclusion (1) follows from \eqref{PolyBdryStrict2}; we have conclusion (2) because $\partial P = P_{2}$ and \eqref{PolyBdryStrict1} holds. For conclusion (3), it suffices to follow the argument above, noting that we have the additional conclusion of Proposition \ref{pp:PolyBdryPolyStrict} in this case. Finally, for conclusion (4), we work as above but apply Corollary \ref{cr:IntegralPolyStrict} instead of Proposition \ref{pp:NormalPolyStrict}.
\end{proof}

\section{Space-Time Strict Polyhedral Approximation}\label{s:STStrictPoly}

In this section, we establish the following general form of Theorem \ref{tm:AbbrvSTPoly}:

\begin{theorem}\label{tm:STPoly} Take flat integrands
\begin{align*}
    \psi &: \R^{d} \times \subext{k}\R^{d} \to [0,\infty), & \phi &: \R^{d} \times \subext{1+k}\R^{d} \to [0,\infty).
\end{align*}
Let $K \subset \R^{d}$ be compact and take $\varepsilon > 0$. Let $S \in \Nrm_{1+k}([0,1] \times \R^{d})$. Suppose that $\partial S \mvert (0,1) \times \R^{d} = 0$ and that $\supp S \subset [0,1] \times \Int K$. Then there is $P \in \Prm_{1+k}([0,1] \times \R^{d})$ such that $\partial P \mvert (0,1) \times \R^{d} = 0$, $\supp P \subset [0,1] \times K$ and the following conclusions hold:
\begin{enumerate}
    \item We have
    \begin{align*}
        \FN_{[0,1] \times K}(P-S) &\leq \varepsilon
    \end{align*}
    and
    \begin{align*}
        \Var_{\phi}(P) &\leq \Var_{\phi}(S) + \varepsilon, & \int \Mbf_{\psi}(P(t)) \dd t \leq \int \Mbf_{\psi}(S(t)) \dd t + \varepsilon.
    \end{align*}
    \item We have
    \begin{align*}
        \FN_{K}^{\circ}(\partial^{\pm}P - \partial^{\pm}S) &\leq \varepsilon, & \Mbf_{\psi}(\partial^{\pm}P) &\leq \Mbf_{\psi}(\partial^{\pm}S) + \varepsilon.
    \end{align*}
    \item We have
    \begin{align*}
        \esssup_{t \in [0,1]} \Mbf_{\psi}(P(t)) &\leq \esssup_{t \in [0,1]}\Mbf_{\psi}(S(t)) + \varepsilon.
    \end{align*}
    \item If $\partial^{\pm}S \in \Irm_{k}(\R^{d})$, then we can additionally ensure that $\partial^{\pm}P \in \IPrm_{k}(\R^{d})$ and that
    \begin{align*}
        \FI_{K}^{\circ}(\partial^{\pm}P-\partial^{\pm}S) &\leq \varepsilon.
    \end{align*}
    \item If $S \in \Irm_{1+k}([0,1] \times \R^{d})$, then we can additionally ensure that $P \in \IPrm_{1+k}([0,1] \times \R^{d})$ and that
    \begin{align*}
        \FI_{[0,1] \times K}(P-S) &\leq \varepsilon.
    \end{align*}
    \item If $\partial^{\pm}S \in \Prm_{k}(\R^{d})$, then we can additionally ensure that $\partial^{\pm}P = \partial^{\pm}S$.
\end{enumerate}
\end{theorem}

\subsection{Piecewise-Constant Space-Time Trajectories} To prove Theorem \ref{tm:STPoly}, we start by considering the special case in which the space-time trajectory to be approximated is piecewise-constant in time:

\begin{lemma}\label{lm:PCPolySTStrict} Take flat integrands
\begin{align*}
    \psi &: \R^{d} \times \subext{k}\R^{d} \to [0,\infty), & \phi &: \R^{d} \times \subext{1+k}\R^{d} \to [0,\infty).
\end{align*}
Let $K \subset \R^{d}$ be compact and take $\varepsilon > 0$. Let $S \in \Nrm_{1+k}([0,1] \times \R^{d})$. Suppose that $\partial S \mvert (0,1) \times \R^{d} = 0$ and that $\supp S \subset [0,1] \times \Int K$. Let $N \in \N$ and take $\{t_{j}\}_{j=0}^{N} \subset [0,1]$ with $0 = t_{0} < \cdots < t_{N} = 1$. Take collections 
\begin{align*}
    \{T_{j}\}_{j=1}^{N} &\subset \Nrm_{k}(\R^{d}), & \{R_{j}\}_{j=0}^{N} &\subset \Nrm_{k+1}(\R^{d})
\end{align*}
and suppose that 
\begin{align*}
    S = \sum_{j=0}^{N} \delta_{t_{j}} \times R_{j} + \sum_{j=1}^{N} [\![t_{j-1},t_{j}]\!] \times T_{j}.
\end{align*}
Then there is $P \in \Prm_{1+k}([0,1] \times \R^{d})$ with $\partial P \mvert (0,1) \times \R^{d} = 0$, $\supp P \subset [0,1] \times K$ such that conclusions (1)-(6) of Theorem \ref{tm:STPoly} hold.
\end{lemma}

\begin{proof} For $j \in \N$, we can use the fact that $\partial S = \delta_{1} \times \partial^{+}S - \delta_{0} \times \partial^{-}S$ to deduce that $\partial T_{j} = 0$ for $1 \leq j \leq N$ and
\begin{align*}
    \partial R_{j} := \left\{\begin{array}{cl}
        T_{1} - \partial^{-}S & \text{if $j=0$}, \\
        T_{j+1} - T_{j} & \text{if $1 \leq j \leq N-1$}, \\
        \partial^{+}S - T_{N} & \text{if $j=N$}.
    \end{array}\right.
\end{align*}
Furthermore, since $\supp S \subset [0,1] \times \Int K$, we have $\supp R_{j} \subset \Int K$ for $0 \leq j \leq N$ and $\supp T_{j} \subset \Int K$ for $1 \leq j \leq N$. Fix $\alpha > 0$ (to be chosen later) and take an intermediate compact set $A \subset \Int K$ with $\supp \partial^{\pm}S \subset \Int A$ and $\supp T_{j} \subset \Int A$ for $1 \leq j \leq N$. We apply Proposition \ref{pp:PolyBdryPolyStrict} to $\partial^{\pm}S$ on $A$ to obtain $P^{\pm} \in \Prm_{k}(\R^{d})$ and $Y^{\pm} \in \Nrm_{k+1}(\R^{d})$ with 
$\supp P^{\pm}$, $\supp Y^{\pm} \subset A$ and
\begin{align}\label{PCPolySTStrict1}
    \partial^{\pm}S &= P^{\pm}+\partial Y^{\pm}, & \Mbf_{\psi}(P^{\pm}) &\leq \Mbf_{\psi}(\partial^{\pm}S) + \alpha, & \Mbf(Y^{\pm}) &\leq \alpha.
\end{align}
Similarly, for $1 \leq j \leq N$, we apply Proposition \ref{pp:PolyBdryPolyStrict} to $T_{j}$ on $A$ to obtain $P_{j}^{T} \in \Prm_{k}(\R^{d})$ and $Y_{j} \in \Nrm_{k+1}(\R^{d})$ with $\supp P_{j}^{T}$, $\supp Y_{j} \subset A$ and
\begin{align}\label{PCPolySTStrict2}
    T_{j} &= P^{T}_{j} + \partial Y_{j}, & \Mbf_{\psi}(P^{T}_{j}) &\leq \Mbf_{\psi}(T_{j}) + \alpha, & \Mbf(Y_{j}) &\leq \alpha.
\end{align}
For $0 \leq j \leq N$, we define $Z_{j} \in \Nrm_{k+1}(\R^{d})$ by
\begin{align*}
    Z_{j} := \left\{\begin{array}{cl}
        Y^{-} - Y_{1} & \text{if $j=0$}, \\
        Y_{j} - Y_{j+1} & \text{if $1 \leq j \leq N-1$}, \\
        Y_{N} - Y^{+} & \text{if $j=N$}.
    \end{array}\right.
\end{align*}
Notice that, for $0 \leq j \leq N$, 
\begin{align}\label{PCPolySTStrict3}
    \partial(R_{j}+Z_{j})= \left\{\begin{array}{cl}
        P_{1}^{T} - P^{-} & \text{if $j=0$}, \\
        P_{j+1}^{T} - P_{j}^{T} & \text{if $1 \leq j \leq N-1$}, \\
        P^{+} - P^{T}_{N} & \text{if $j=N$}.
    \end{array}\right.
\end{align}
In particular, we have $\partial(R_{j}+Z_{j}) \in \Prm_{k}(\R^{d})$ for $0 \leq j \leq N$. Since $\supp R_{j}+Z_{j} \subset \Int K$ for $0 \leq j \leq N$, we can appeal to Proposition \ref{pp:PolyBdryPolyStrict} to obtain $P^{R}_{j} \in \Prm_{k+1}(\R^{d})$ and $X_{j} \in \Nrm_{k+2}(\R^{d})$ with $\supp P^{R}_{j}$, $\supp X_{j} \subset K$ with 
\begin{align}\label{PCPolySTStrict4}
    R_{j}+Z_{j} &= P^{R}_{j} + \partial X_{j}, & \Mbf_{\phi}(P^{R}_{j}) &\leq \Mbf_{\phi}(R_{j}+Z_{j}) + \alpha, & \Mbf(X_{j}) &\leq \alpha.
\end{align}
Now define
\begin{align*}
    P := \sum_{j=0}^{N} \delta_{t_{j}} \times P^{R}_{j} + \sum_{j=1}^{N} [\![t_{j-1},t_{j}]\!] \times P^{T}_{j}.
\end{align*}
Notice that $P \in \Prm_{1+k}(\R^{1+d})$ and $\supp P \subset [0,1] \times K$. One can use \eqref{PCPolySTStrict1}, \eqref{PCPolySTStrict2}, \eqref{PCPolySTStrict3}, and \eqref{PCPolySTStrict4} to deduce that 
\begin{align*}
    P - S = \partial X + Y
\end{align*}
where $X \in \Nrm_{k+2}(\R^{1+d})$ and $Y \in \Nrm_{k+1}(\R^{1+d})$ are given by
\begin{align*}
    X &= \sum_{j=1}^{N} [\![t_{j-1},t_{j}]\!] \times Y_{j} - \sum_{j=0}^{N} \delta_{t_{j}} \times X_{j}, & Y = \delta_{0} \times Y^{-} - \delta_{1} \times Y^{+}.
\end{align*}
It follows that 
\begin{align}\label{PCPolySTStrict4.5}
    \partial P \mvert (0,1) \times \R^{d} = 0 \hspace{1em} \text{with $\partial^{\pm}P = P^{\pm}$}.
\end{align}

Notice that $\supp X \subset [0,1] \times K$ and that $\supp Y \subset \{0,1\} \times K$. We can use \eqref{PCPolySTStrict1}, \eqref{PCPolySTStrict2}, and \eqref{PCPolySTStrict4} to deduce that
\begin{align}\label{PCPolySTStrict5}
    \Mbf(X) &\leq (N+2)\alpha, & \Mbf(Y) &\leq 2\alpha.
\end{align}
We then use \eqref{PCPolySTStrict4} and Lemma \ref{lm:DissJumps} to estimate
\begin{equation}\label{PCPolySTStrict6}
\begin{aligned}
    \Var_{\phi}(P) &\leq \sum_{j=0}^{N} \Mbf_{\phi}(R_{j}) + \Mbf_{\phi}(Z_{j}) + \alpha \\
    &\leq \sum_{j=0}^{N} \Mbf_{\phi}(R_{j}) +(2C_{\phi,K}+1)\alpha = \Var_{\phi}(S) + (N+1)(2C_{\phi,K}+1)\alpha
\end{aligned}
\end{equation}
where $C_{\phi,K} > 0$ is defined as in \eqref{psiVsMass}. On the other hand, we can use \eqref{PCPolySTStrict2} to deduce that
\begin{equation}\label{PCPolySTStrict7}
\begin{aligned}
    \int \Mbf_{\psi}(P(t)) \dd t &= \sum_{j=1}^{N} \vert t_{j} - t_{j-1} \vert \cdot \Mbf_{\psi}(P^{T}_{j}) \\
    &\leq \sum_{j=1}^{N} \vert t_{j} - t_{j-1} \vert \cdot (\Mbf_{\psi}(T_{j}) + \alpha) = \int \Mbf_{\psi}(S(t)) \dd t + \alpha.
\end{aligned}
\end{equation}
and
\begin{equation}\label{PCPolySTStrict8}
\begin{aligned}
    \esssup_{t \in [0,1]} \Mbf_{\psi}(P(t))& = \max_{1 \leq j \leq N} \Mbf_{\psi}(P^{T}_{j}) \\
    &\leq \max_{1 \leq j \leq N} \Mbf_{\psi}(T_{j}) + \alpha = \esssup_{t \in [0,1]} \Mbf_{\psi}(S(t)) + \alpha.
\end{aligned}
\end{equation}

We now choose $\alpha$ so that
\begin{align*}
    \alpha \leq \varepsilon \cdot \min\bigg\{\frac{1}{4},\frac{1}{2(N+2)},\frac{1}{(N+1)(2C_{\phi,K}+1)}\bigg\}.
\end{align*}
Conclusion (1) now follows from \eqref{PCPolySTStrict5}, \eqref{PCPolySTStrict6}, and \eqref{PCPolySTStrict7}; conclusion (2) follows from \eqref{PCPolySTStrict4.5} and \eqref{PCPolySTStrict1}, while conclusion (3) follows from \eqref{PCPolySTStrict8}. For conclusion (4), it suffices to follow the proof above, noting that the additional conclusion of Proposition \ref{pp:PolyBdryPolyStrict} ensures that $P^{\pm} \in \IPrm_{k}(\R^{d})$ and $Y^{\pm} \in \Irm_{k+1}(\R^{d})$. Similarly, for conclusion (5), we follow the proof above, noting that the additional conclusion of Proposition \ref{pp:PolyBdryPolyStrict} holds throughout. Finally, to verify conclusion (6), notice that we can take $Y^{\pm} = 0$ in this case.
\end{proof}

\subsection{Proof of Theorem \ref{tm:STPoly}}

We shall need the following lemma:

\begin{lemma}\label{lm:Sweep} Take $\sigma,\tau \in \R$ and let $H: [0,1] \times \R^{1+d}$ be the affine homotopy from $\mathrm{id}_{\R^{1+d}}$ to $\tau \mathbf{1}_{\R} \times \mathrm{id}_{\R^{d}}$. Let $R \in \Nrm_{k}(\R^{1+d})$ and suppose that $\supp R \subset \{\sigma\} \times \R^{d}$. Then
\begin{align*}
    H_{*}([\![0,1]\!] \times R) = [\![\sigma,\tau]\!] \times \pbf_{*}R.
\end{align*} 
\end{lemma}

\begin{proof} First notice that, for $s \in [0,1]$ and $(t,x) \in \R^{1+d}$, we have
\begin{align*}
    H(s,(t,x)) = (1-s)(t,x) + s(\tau,x) = ((1-s)t+s\tau,x)
\end{align*}
Let $H_{\sigma}: [0,1] \times \R^{1+d} \to \R^{1+d}$ where
\begin{align*}
    H_{\sigma}(s,(t,x)) = H(s,(\sigma,x)) \hspace{1em} \text{for $s \in [0,1]$, $(t,x) \in \R^{1+d}$.}
\end{align*}
Notice that $H_{\sigma} = \overline{H}_{\sigma,\tau} \circ \pi$, where
$\overline{H}_{\sigma,\tau}:[0,1] \times \R^{d} \to \R \times \R^{d}$ is given by
\begin{align*}
    \overline{H}_{\sigma,\tau}(s,x) = ((1-s)\sigma + s\tau,x) \hspace{1em} \text{for $s \in [0,1]$, $x \in \R^{d}$}
\end{align*}
and $\pi: [0,1] \times \R^{1+d} \to [0,1] \times \R^{d}$ is given by
\begin{align*}
    \pi(s,(t,x)) = (s,x) \hspace{1em} \text{for $s \in [0,1]$, $(t,x) \in \R^{1+d}$.}
\end{align*}

Since $\supp R \subset \{\sigma\} \times \R^{d}$, we can use \cite[4.1.15]{Federer69book} to deduce that
\begin{align*}
    H_{*}([\![0,1]\!] \times R) = (H_{\sigma})_{*}([\![0,1]\!] \times R) = (\overline{H}_{\sigma,\tau} \circ \pi)_{*}([\![0,1]\!] \times R).
\end{align*}
One can check that
\begin{align*}
    \pi_{*}([\![0,1]\!] \times R) = [\![0,1]\!] \times \pbf_{*}R.
\end{align*}
The result then follows after observing that
\begin{align*}
    (\overline{H}_{\sigma,\tau})_{*}([\![0,1]\!] \times \pbf_{*}R) &= [\![\sigma,\tau]\!] \times \pbf_{*}R. \qedhere
\end{align*}
\end{proof}

\begin{proof}[Proof of Theorem \ref{tm:STPoly}] Fix $\alpha > 0$ (to be chosen later). Let $N \in \N$ and take $\{t_{j}\}_{j=0}^{N} \subset [0,1]$ with $0 = t_{0} < \cdots < t_{N} = 1$ such that $\tbf_{\#}\Vert S \Vert(\{t_{j}\}) = 0$ for $1 \leq j \leq N-1$ and $\vert t_{j} - t_{j-1} \vert \leq \alpha$ for $1 \leq j \leq N$. Take $t \in \R$ and notice that
\begin{align}\label{STPoly1}
    (t\mathbf{1}_{\R} \times \mathrm{id}_{\R^{d}})_{*}R = \delta_{t} \times \pbf_{*}R \hspace{1em} \text{for $R \in \Nrm_{k}(\R^{1+d})$}.
\end{align}
Let $H_{t}:[0,1] \times \R^{1+d} \to \R^{1+d}$ be the affine homotopy from $\mathrm{id}_{\R^{1+d}}$ to $t\mathbf{1}_{\R} \times \mathrm{id}_{\R^{d}}$.

For $1 \leq j \leq N$, take $a_{j} \in [t_{j-1},t_{j}]$ with $\tbf_{\#}\Vert S \Vert(\{a_{j}\}) = 0$,
\begin{align}\label{STPoly2}
    \Mbf_{\psi}(S(a_{j})) &\leq \frac{1}{t_{j}-t_{j-1}}\int_{t_{j-1}}^{t_{j}} \Mbf_{\psi}(S(t)) \dd t.
\end{align}
and
\begin{align}\label{STPoly3}
    \Mbf_{\psi}(S(a_{j})) \leq \esssup_{t \in [0,1]} \Mbf_{\psi}(S(t)).
\end{align}

For $1 \leq j \leq N$, set
\begin{align}\label{STPoly4}
    S^{-}_{j} &= S \mvert [t_{j-1},a_{j}] \times \R^{d}, & S^{+}_{j} &= S \mvert [a_{j},t_{j}] \times \R^{d},
\end{align}
Since $\tbf_{\#}\Vert S \Vert(\{a_{j}\}) = 0$ for $1 \leq j \leq N$ and $\tbf_{\#}\Vert S \Vert(\{t_{j}\}) = 0$ for $1 \leq j \leq N-1$, we can use \cite[Proposition 3.4]{BonicattoRindlerTurnbull26} to deduce that
\begin{equation}\label{STPoly5}
\begin{aligned}
    \partial S^{-}_{j} &= S\vert_{a_{j}} - S\vert_{t_{j-1}} \hspace{1em} \text{for $2 \leq j \leq N$}, \\
    \partial S^{+}_{j} &= S\vert_{t_{j}} - S\vert_{a_{j}} \hspace{1em} \text{for $1 \leq j \leq N-1$}.
\end{aligned}
\end{equation}
Similarly, it follows from \cite[Proposition 3.4]{BonicattoRindlerTurnbull26} that
\begin{equation}\label{STPoly6}
\begin{aligned}
    \partial S_{1}^{-} &= S\vert_{a_{1}} - \delta_{0} \times \partial^{-}S \\
    \partial S_{N}^{+} &= \delta_{1} \times \partial^{+}S - S\vert_{a_{N}}.
\end{aligned}
\end{equation}
Now fix $j \in \N$ with $1 \leq j \leq N$. We can use the homotopy formula \eqref{Homotopy} and \eqref{STPoly1} to deduce that
\begin{align}\label{STPoly7}
    \delta_{t_{j-1}} \times \pbf_{*}S^{-}_{j} - S^{-}_{j} &= (H_{t_{j-1}})_{*}([\![0,1]\!] \times \partial S^{-}_{j})+ \partial (H_{t_{j-1}})_{*}([\![0,1]\!] \times S^{-}_{j})
\end{align}
and
\begin{align}\label{STPoly8}
    \delta_{t_{j}} \times \pbf_{*}S^{+}_{j} - S^{+}_{j} &= (H_{t_{j}})_{*}([\![0,1]\!] \times \partial S^{+}_{j}) + \partial (H_{t_{j}})_{*}([\![0,1]\!] \times S^{+}_{j}).
\end{align}
Since \eqref{STPoly5} and \eqref{STPoly6} hold, we can appeal to Lemma \ref{lm:Sweep} to deduce that
\begin{equation}\label{STPoly9}
\begin{aligned}
    (H_{t_{j-1}})_{*}([\![0,1]\!] \times \partial S^{-}_{j}) &= -[\![t_{j-1},a_{j}]\!] \times S(a_{j}) \\ (H_{t_{j}})_{*}([\![0,1]\!] \times \partial S^{+}_{j}) &= -[\![a_{j},t_{j}]\!]
    \times S(a_{j}).
\end{aligned}
\end{equation}
Combining \eqref{STPoly7}, \eqref{STPoly8}, and \eqref{STPoly9} yields
\begin{equation}
\begin{aligned}\label{STPoly10}
    \delta_{t_{j-1}} &\times \pbf_{*}S^{-}_{j} + \delta_{t_{j}} \times \pbf_{*}S_{j}^{+} + [\![t_{j-1},t_{j}]\!] \times S(a_{j}) \\
    &= S \mvert [t_{j-1},t_{j}] \times \R^{d} + \partial((H_{t_{j-1}})_{*}([\![0,1]\!] \times S_{j}^{-}) + (H_{t_{j}})_{*}([\![0,1]\!] \times S_{j}^{+}))
\end{aligned}
\end{equation}
Summing \eqref{STPoly10} over $j \in \N$ with $1 \leq j \leq N$ and recalling \eqref{STPoly4} yields
\begin{align}\label{STPoly11}
    \sum_{j=0}^{N} \delta_{t_{j}} \times R_{j} + \sum_{j=1}^{N} [\![t_{j-1},t_{j}]\!] \times S(a_{j}) = S + \partial X
\end{align}
where, for $0 \leq j \leq N$, $R_{j} \in \Nrm_{k+1}(\R^{d})$ is given by
\begin{align}\label{STPoly12}
    R_{j} = \left\{\begin{array}{cl}
        \pbf_{*}(S \mvert [0,a_{1}] \times \R^{d}) & \text{if $j=0$}, \\
        \pbf_{*}(S \mvert [a_{j},a_{j+1}] \times \R^{d}) & \text{if $1 \leq j \leq N-1$}, \\
        \pbf_{*}(S \mvert [a_{N},1] \times \R^{d}) & \text{if $j=N$}.
    \end{array}\right.
\end{align}
and $X \in \Nrm_{k+2}(\R^{1+d})$ is given by
\begin{align*}
    X = \sum_{j=1}^{N}(H_{t_{j-1}})_{*}([\![0,1]\!] \times S_{j}^{-}) + (H_{t_{j}})_{*}([\![0,1]\!] \times S_{j}^{+}).
\end{align*}
Notice that $\supp R_{j} \subset \Int K$ for $0 \leq j \leq N$ and that $\supp X \subset [0,1] \times \Int K$. For $1 \leq j \leq N$, we can use \cite[4.1.9]{Federer69book} to estimate
\begin{align*}
    \Mbf((H_{t_{j-1}})_{*}([\![0,1]\!] \times S_{j}^{-})) &\leq \vert a_{j}-t_{j-1}\vert \cdot \Mbf(S_{j}^{-}) \leq \alpha \cdot \Mbf(S^{-}_{j})\\
    \Mbf((H_{t_{j}})_{*}([\![0,1]\!] \times S_{j}^{+})) &\leq \vert t_{j} - a_{j} \vert \cdot \Mbf(S^{+}_{j}) \leq \alpha \cdot \Mbf(S^{+}_{j}).
\end{align*}
It follows that 
\begin{align}\label{STPoly13}
    \Mbf(X) \leq \alpha \cdot \Mbf(S).
\end{align}
Now set 
\begin{align*}
    \tilde{S} = \sum_{j=0}^{N} \delta_{t_{j}} \times R_{j} + \sum_{j=1}^{N} [\![t_{j-1},t_{j}]\!] \times S(a_{j})
\end{align*}
Notice that $\supp \tilde{S} \subset [0,1] \times \Int K$ and we can use \eqref{STPoly11} to deduce that $\tilde{S}-S = \partial X$. In particular, we have
\begin{align}\label{STPoly14}
    \partial \tilde{S} = \partial S.
\end{align}
We can use Lemma \ref{lm:DissJumps}, \eqref{STPoly12}, and Lemma \ref{lm:pPhiMassByDiss} to estimate
\begin{align*}
    \Var_{\phi}(\tilde{S}) \leq \Var_{\phi}(S;[0,a_{1}]) + \sum_{j=1}^{N-1} \Var_{\phi}(S;[a_{j},a_{j+1}]) + \Var_{\phi}(S;[a_{N},1])
\end{align*}
Since we also have $\tbf_{\#}\Vert S \Vert(\{a_{j}\}) = 0$ for $1 \leq j \leq N$, it follows that
\begin{align}\label{STPoly15}
    \Var_{\phi}(\tilde{S}) &\leq \Var_{\phi}(S).
\end{align}
We now use \eqref{STPoly2} to estimate
\begin{align}\label{STPoly16}
    \int \Mbf_{\psi}(\tilde{S}(t)) \dd t = \sum_{j=1}^{N} \vert t_{j} -t_{j-1} \vert \cdot \Mbf_{\psi}(S(a_{j})) \leq \int\Mbf_{\psi}(S(t)) \dd t.
\end{align}
Similarly, we use \eqref{STPoly3} to estimate
\begin{align}\label{STPoly17}
    \esssup_{t \in [0,1]} \Mbf_{\psi}(\tilde{S}(t)) = \max_{1 \leq j \leq N} \Mbf_{\psi}(S(a_{j})) \leq \esssup_{t \in [0,1]} \Mbf_{\psi}(S(t)).
\end{align}
We now choose $\alpha$ so that
\begin{align*}
    \alpha \cdot \max\{1, \Mbf(S)\} \leq \frac{\varepsilon}{2}.
\end{align*}
Finally, we apply Lemma \ref{lm:PCPolySTStrict} to $\tilde{S}$ with $\varepsilon/2$ in place of $\varepsilon$. Conclusions (1) and (3) now follow from \eqref{STPoly13}, \eqref{STPoly15}, \eqref{STPoly16}, \eqref{STPoly17}, and Lemma \ref{lm:PCPolySTStrict}, while conclusions (2), (4), and (6) follow from \eqref{STPoly14} and Lemma \ref{lm:PCPolySTStrict}. To verify conclusion (5), notice that we can follow the argument above whilst additionally ensuring that $S(a_{j}) \in \Irm_{k}(\R^{d})$ for $1 \leq j \leq N$ (see \cite[4.3.6]{Federer69book}). 
\end{proof}

\section{Rational Approximation}\label{s:Q}

Recall from the introduction that a normal current $T \in \Nrm_{k}(\R^{d})$ is \textbf{rational} if there is $n \in \N$ such that $n \cdot T \in \Irm_{k}(\R^{d})$. We shall use $\Qrm_{k}(\R^{d})$ to denote the space of rational currents. In this section, we shall prove various strict-type rational approximation results for normal currents and space-time trajectories. 

\subsection{Polyhedral Chains} To obtain strict-type approximation results concerning normal currents and space-time trajectories, notice that it suffices, by the results of Section \ref{s:ClassicalPolyStrict} and Section \ref{s:STStrictPoly} respectively, to prove approximation results for polyhedral chains. To this end, take $P \in \Prm_{k}(\R^{d})$. Then there is $N \in \N$, a collection $\{P_{i}\}_{i=1}^{N}$ of oriented $k$-simplices in $\R^{d}$, and $\alpha \in \R^{N}$ with
\begin{align*}
    P = \sum_{j=1}^{N} \alpha_{j}P_{j}
\end{align*}
Let $\Xbf: \Prm_{k}(\R^{d}) \to [0,\infty)$ be a seminorm on $\Prm_{k}(\R^{d})$. Notice that, for any $\beta \in \Q^{N}$, we have
\begin{align*}
    \sum_{j=1}^{N}\beta_{j}P_{j} \in \QPrm_{k}(\R^{d}) := \Qrm_{k}(\R^{d}) \cap \Prm_{k}(\R^{d})
\end{align*}
and
\begin{align}\label{NaiveRationalApprox}
    \Xbf\bigg(\sum_{j=1}^{N} \beta_{j}P_{j} - P\bigg) \leq \max_{1 \leq j \leq N} \Xbf(P_{j}) \cdot \vert \beta - \alpha \vert_{1}.
\end{align}
Since $\Q^{N}$ is dense in $\R^{N}$, we see from \eqref{NaiveRationalApprox} that we can approximate $P$ by elements of $\QPrm_{k}(\R^{d})$ arbitrarily well in $\Xbf$. However, this procedure does not control the support of the boundary of the approximants in any way. We shall now see how to refine the procedure outlined above to ensure that the support of the boundary of the approximants is controlled; the following result is similar to \cite[Section 5.4.1, Proposition 2]{GiaquintaModicaSoucek98book1} and \cite[5.5]{Federer74}:

\begin{proposition}\label{pp:PolyQX} Let $\Xbf: \Prm_{k}(\R^{d}) \to [0,\infty)$ be a seminorm and take $\varepsilon > 0$. Let $P \in \Prm_{k}(\R^{d})$ and suppose that $\partial P \mvert A \in \Qrm_{k-1}(\R^{d})$ for some Borel $A \subset \R^{d}$. Then there is $Q \in \QPrm_{k}(\R^{d})$ with $\supp Q \subset \supp P$ and
\begin{align*}
    \Xbf(Q-P) &\leq \varepsilon, & \supp \partial Q &\subset \supp \partial P, & \partial Q \mvert A &= \partial P \mvert A.
\end{align*}
\end{proposition}

\begin{proof} Since $P \in \Prm_{k}(\R^{d})$, there is $N \in \N$, an $\mathscr{H}^{k}$-disjoint collection $\{P_{j}\}_{j=1}^{N}$ of oriented $k$-simplices (which are identified with their associated currents as in \cite[4.1.11]{Federer69book}) in $\R^{d}$, and $\alpha \in \R^{N}$ with
\begin{align*}
    P = \sum_{j=1}^{N} \alpha_{j}P_{j}.
\end{align*}
Take $M \in \N$ and let $\{R_{i}\}_{i=1}^{M}$ be an $\mathscr{H}^{k-1}$-disjoint collection of oriented $(k-1)$-simplices obtained by taking a common finite simplicial subdivision of all $(k-1)$ faces of the simplices in the collection $\{P_{j}\}_{j=1}^{N}$. Fix $j \in \{1,\dotsc,N\}$. Notice that there is $\{\theta_{ji}\}_{i=1}^{M} \subset \{0,\pm 1\}$ such that
\begin{align*}
    \partial P_{j} = \sum_{i=1}^{M} \theta_{ji}R_{i}.
\end{align*}
Then, for $\beta \in \R^{N}$, we have
\begin{align*}
    \partial\bigg(\sum_{j=1}^{N}\beta_{j}P_{j}\bigg) = \sum_{i=1}^{M}\bigg(\sum_{j=1}^{N} \beta_{j}\theta_{ji}\bigg) R_{i}.
\end{align*}
It follows that, for $\beta \in \R^{N}$,
\begin{align}\label{PolyQX1}
    \supp \partial \bigg(\sum_{j=1}^{N}\beta_{j}P_{j}\bigg) = \bigcup_{i \in I^{+}(\beta)} \supp R_{i}
\end{align}
where, for $\beta \in \R^{N}$,
\begin{align*}
    I^{+}(\beta) = \bigg\{i \in \{1,\dotsc,M\}: \bigg\vert \sum_{j=1}^{N} \beta_{j}\theta_{ji} \bigg\vert > 0\bigg\}.
\end{align*}
Define
\begin{align*}
    I := \bigg\{i \in \{1,\dotsc,M\}: \sum_{j=1}^{N}\alpha_{j}\theta_{ji} \in \Q\bigg\}
\end{align*}
and let $V \subset \R^{N}$ be the following affine subspace of $\R^{N}$:
\begin{align*}
    V := \bigg\{\beta \in \R^{N}: \sum_{j=1}^{N} \beta_{j}\theta_{ji} = \sum_{j=1}^{N}\alpha_{j}\theta_{ji} \hspace{1em} \text{for $i \in I$}\bigg\}.
\end{align*}
Fix $\beta \in V$ (to be chosen later) and define
\begin{align*}
    Q := \sum_{j=1}^{N} \beta_{j}P_{j}.
\end{align*}
Notice that $\supp Q \subset \supp P$. Furthermore, since $\beta \in V$ and $0 \in \Q$, we have $I^{+}(\beta) \subset I^{+}(\alpha)$; it then follows from \eqref{PolyQX1} that
\begin{align*}
    \supp \partial Q \subset \supp \partial P.
\end{align*}
Since $\partial P \mvert A \in \Qrm_{k-1}(\R^{d})$, we have
\begin{align*}
    \sum_{i=1}^{M}\bigg(\sum_{j=1}^{N} \alpha_{j}\theta_{ji}\bigg)R_{i} \mvert A \in \Qrm_{k-1}(\R^{d}).
\end{align*}
It follows that, if $i \in \{1,\dotsc,M\}$ and $i \not\in I$, then $R_{i} \mvert A = 0$. We therefore have
\begin{align*}
    \partial Q \mvert A = \sum_{i \in I}\sum_{j=1}^{N} \beta_{j}\theta_{ji} R_{i} \mvert A = \sum_{i \in I}\sum_{j=1}^{N} \alpha_{j}\theta_{ji} R_{i} \mvert A = \partial P \mvert A.
\end{align*}

It therefore remains to choose $\beta$ in such a way that $Q \in \QPrm_{k}(\R^{d})$ and $\Xbf(Q-P) \leq \varepsilon$. To this end, let $\Theta: \R^{I} \to \R^{N}$ be the linear map with
\begin{align*}
    (\Theta f)_{j}=\sum_{i \in I} \theta_{ji}f_{i} \hspace{1em} \text{for $1 \leq j \leq N$}.
\end{align*}
Here, we recall that $\R^{I}$ is the space of functions $f: I \to \R$; we write $f_{i} := f(i)$ for $f \in \R^{I}$ and $i \in I$. Notice that
\begin{align*}
    V = \{\beta \in \R^{N}: \Theta^{*}\beta = \Theta^{*}\alpha\}.
\end{align*}
Since $\Theta^{*}\alpha(I) \subset \Q$ and $\{\theta_{ji}\}_{i=1}^{M} \subset \Q$ for $1 \leq j \leq N$, we can find $\tilde{\beta} \in \Q^{N}$ with $\Theta^{*}\tilde{\beta} = \Theta^{*}\alpha$. It follows that $V = \tilde{\beta} + \ker \Theta^{*}$. Since $\{\theta_{ji}\}_{i=1}^{M} \subset \Q$ for $1 \leq j \leq N$, we can find a rational basis $\Gamma \subset \Q^{N}$ for $\ker \Theta^{*}$ and, without loss of generality, we can assume that $\Gamma \subset \Z^{N}$. Let $\Phi: \R^{\vert \Gamma \vert} \to \ker \Theta^{*}$ be a linear map with
\begin{align*}
    \{\Phi(e_{i})\}_{i=1}^{\vert \Gamma \vert} = \Gamma,
\end{align*}
where $\{e_{i}\}_{i=1}^{\vert \Gamma \vert}$ is the standard basis of $\R^{\vert \Gamma \vert}$. Take $C > 0$ with
\begin{align}\label{PolyQX2}
    \vert \Phi(v) \vert_{1} \leq C \vert v \vert_{1} \hspace{1em} \text{for $v \in \R^{\vert \Gamma \vert}$}.
\end{align}
Since $\alpha \in V$, there is $v \in \R^{\vert \Gamma \vert}$ with $\alpha = \tilde{\beta} + \Phi(v)$. Take $w \in \Qbb^{\vert \Gamma \vert}$ with 
\begin{align}\label{PolyQX3}
    C \cdot \vert w - v \vert_{1} \cdot \max_{1 \leq j \leq N} \Xbf(P_{j}) \leq \varepsilon.
\end{align}
Now choose $\beta := \tilde{\beta} + \Phi(w)$. Since $w \in \Qbb^{\vert \Gamma \vert}$ and $\tilde{\beta} \in \Qbb^{N}$, we have $\beta \in \Qbb^{N}$ and so $Q \in \QPrm_{k}(\R^{d})$. Finally, since $\Phi(w-v) = \beta - \alpha$, we can use \eqref{PolyQX2} and \eqref{PolyQX3} to deduce that 
\begin{align*}
    \Xbf(P-Q) = \Xbf\bigg(\sum_{j=1}^{N}(\beta_{j}-\alpha_{j})P_{j}\bigg)
    \leq \vert \alpha - \beta\vert_{1} \max_{1 \leq j \leq N} \Xbf(P_{j}) \leq \varepsilon.
\end{align*}
The result follows.
\end{proof}

\subsection{Normal Currents} By combining Proposition \ref{pp:PolyQX} with results from Section \ref{s:ClassicalPolyStrict}, one can obtain many strict-type rational approximation results for normal currents. In what follows, we give two examples of such results:

\begin{proposition}\label{pp:NormalQBdryStrict} Take flat integrands
\begin{align*}
    \psi^{1} &: \R^{d} \times \subext{k}\R^{d} \to [0,\infty), & \psi^{2} &: \R^{d} \times \subext{k-1}\R^{d} \to [0,\infty).
\end{align*}
Let $K_{1},K_{2} \subset \R^{d}$ be compact with $K_{2} \subset K_{1}$ and take $\varepsilon > 0$. Let $T \in \Nrm_{k}(\R^{d})$. Suppose that $\supp T \subset \Int K_{1}$ and that $\supp \partial T \subset \Int K_{2}$. Then there is $Q \in \QPrm_{k}(\R^{d})$ such that $\supp Q \subset K_{1}$, $\supp \partial Q \subset K_{2}$, and
\begin{align*}
    \FN_{K_{1}}(Q-T) &\leq \varepsilon, & \Mbf_{\psi^{1}}(Q) &\leq \Mbf_{\psi^{1}}(T) + \varepsilon, \\
    \Fbf_{K_{2}}(\partial Q - \partial T) &\leq \varepsilon, & \Mbf_{\psi^{2}}(\partial Q) &\leq \Mbf_{\psi^{2}}(\partial T) + \varepsilon.
\end{align*}
\end{proposition}

\begin{proof} Fix $\alpha > 0$ (to be chosen later). We first apply Theorem \ref{tm:PolyBdryStrict} to $T$ to find $P \in \Prm_{k}(\R^{d})$ with $\supp P \subset K_{1}$, $\supp \partial P \subset K_{2}$, and 
\begin{align*}
    \FN_{K_{1}}(P-T) &\leq \alpha, & \Mbf_{\psi^{1}}(P) &\leq \Mbf_{\psi^{1}}(T) + \alpha, \\
    \Fbf_{K_{2}}(\partial P - \partial T) &\leq \alpha, & \Mbf_{\psi^{2}}(\partial P) &\leq \Mbf_{\psi^{2}}(\partial T) + \alpha.
\end{align*}
We apply Proposition \ref{pp:PolyQX} to $P$ with $\Xbf = \Nbf$ and $A = \varnothing$ to obtain $Q \in \QPrm_{k}(\R^{d})$ with
\begin{align*}
    \supp Q &\subset \supp P \subset K_{1}, & \supp \partial Q &\subset \supp\partial P \subset K_{2}, & \Nbf(Q-P) &\leq \alpha.
\end{align*}
Then
\begin{align*}
    \FN_{K_{1}}(Q-T) &\leq 2\alpha, \\
    \FN_{K_{2}}(\partial Q - \partial T) &\leq 2\alpha,
\end{align*}
and 
\begin{align*}
    \Mbf_{\psi^{1}}(Q) &\leq \Mbf_{\psi^{1}}(T) +  (C_{\psi^{1},K_{1}}+1)\alpha, \\
    \Mbf_{\psi^{2}}(\partial Q) &\leq \Mbf_{\psi^{2}}(\partial T) + (C_{\psi^{2},K_{2}}+1)\alpha,
\end{align*}
where $C_{\psi^{1},K_{1}},C_{\psi^{2},K_{2}} > 0$ are defined as in \eqref{psiVsMass}. The result follows after choosing $\alpha$ so that
\begin{align*}
    \alpha &\leq \varepsilon \cdot \min\bigg\{\frac{1}{2},\frac{1}{ C_{\psi^{1},K_{1}}+1},\frac{1}{C_{\psi^{2},K_{2}}+1}\bigg\}. \qedhere
\end{align*}
\end{proof}

\begin{proposition}\label{pp:QBdryQStrict} Let $\psi: \R^{d} \times \subext{k}\R^{d} \to [0,\infty)$ be a flat integrand. Let $K \subset \R^{d}$ be compact and take $\varepsilon > 0$. Let $T \in \Nrm_{k}(\R^{d})$. Suppose that $\supp T \subset \Int K$ and that $\partial T \in \Qrm_{k-1}(\R^{d})$. Then there is $Q \in \Qrm_{k}(\R^{d})$ with $\supp Q \subset K$ and
\begin{align*}
    \partial Q &= \partial T, & \Fbf_{K}(Q-T) &\leq \varepsilon, & \Mbf_{\psi}(Q) &\leq \Mbf_{\psi}(T) + \varepsilon.
\end{align*}
\end{proposition}

\begin{proof} Fix $\alpha > 0$ (to be chosen later). It is not restrictive to assume that $\partial T \in \Irm_{k-1}(\R^{d})$. We first apply Theorem \ref{tm:PolyBdryStrict}, noting that we have conclusion (3) in this case, to find $P \in \Prm_{k}(\R^{d})$ with $\partial P \in \IPrm_{k-1}(\R^{d})$, $\supp P \subset K$, and
\begin{align*}
    \FN_{K}(P-T) &\leq \alpha, & \Mbf_{\psi}(P) &\leq \Mbf_{\psi}(T) + \alpha, &
    \FI^{\circ}_{K}(\partial P - \partial T) &\leq \alpha.
\end{align*}
We then apply Proposition \ref{pp:PolyQX} to $P$ with $\Xbf = \Mbf$ and $A = \R^{d}$ to obtain $\tilde{Q} \in \QPrm_{k}(\R^{d})$ with $\supp \tilde{Q} \subset K$ and
\begin{align*}
    \partial \tilde{Q} &= \partial P, & \Mbf(\tilde{Q}-P) &\leq \alpha.
\end{align*}
Since $\FI^{\circ}_{K}(\partial P - \partial T) \leq \alpha$, there is $R \in \Irm_{k}(\R^{d})$ with $\supp R \subset K$ and 
\begin{align*}
    \partial R &= \partial P - \partial T, & \Mbf(R) &\leq \alpha.
\end{align*}
Now set $Q = \tilde{Q} - R$. Notice that $\partial Q = \partial T$; since $\tilde{Q} \in \Qrm_{k}(\R^{d})$ and $R \in \Irm_{k}(\R^{d})$, we have $Q \in \Qrm_{k}(\R^{d})$. Combining the estimates above yields
\begin{align*}
    \Fbf_{K}(Q-T) &\leq \Mbf(R) + \Mbf(\tilde{Q}-P) + \Fbf_{K}(P-T) \leq 3\alpha
\end{align*}
and
\begin{align*}
    \Mbf_{\psi}(Q) &\leq \Mbf_{\psi}(P) + C_{\psi,K}\cdot \Mbf(Q-P) \leq \Mbf_{\psi}(T) + (2C_{\psi,K}+1)\alpha,
\end{align*}
where $C_{\psi,K} > 0$ is defined as in \eqref{psiVsMass}. The result follows after choosing $\alpha$ so that
\begin{align*}
    \alpha &\leq \varepsilon \cdot \min \bigg\{\frac{1}{3},\frac{1}{2C_{\psi,K}+1}\bigg\}. \qedhere
\end{align*}
\end{proof}

Take $T \in \Irm_{k}(\R^{d})$ with $\partial T = 0$. We define the \textbf{global (normal)} and the \textbf{global (integral) homogeneous flat norm} of $T$ respectively as follows:
\begin{align*}
    \FN^{\circ}(T) &:= \inf \bigl\{\Mbf(R): R \in \Nrm_{k+1}(\R^{d}), \hspace{0.2em} \partial R = T \bigr\}, \\
    \FI^{\circ}(T) &:= \inf \bigl\{\Mbf(R): R \in \Irm_{k+1}(\R^{d}), \hspace{0.2em} \partial R = T \bigr\}.
\end{align*}
It turns out that Proposition \ref{pp:QBdryQStrict} can be used to recover the classical result of Federer (see \cite[5.8]{Federer74} or \cite[Section 1.3.4, Corollary 2]{GiaquintaModicaSoucek98book2}) concerning the relationship between these quantities:

\begin{corollary}\label{cr:FIvsFN} Suppose that $k \geq 1$. Take $T \in \Irm_{k}(\R^{d})$ with $\partial T = 0$. Then
\begin{align*}
    \FN^{\circ}(T) = \lim_{n \to \infty} \frac{\FI^{\circ}(n \cdot T)}{n}.
\end{align*}
\end{corollary}

\begin{proof} Since
\begin{align*}
    n \cdot \FN^{\circ}(T) = \FN^{\circ}(n \cdot T) \leq \FI^{\circ}(n \cdot T) \hspace{1em} \text{for $n \in \N$},
\end{align*}
we have
\begin{align*}
    \FN^{\circ}(T) &\leq \inf_{n \in \N} \frac{\FI^{\circ}(n \cdot T)}{n}.
\end{align*}
For the reverse inequality, take $R \in \Nrm_{k+1}(\R^{d})$ with $\partial R = T$ and 
\begin{align*}
    \Mbf(R) = \FN^{\circ}(T).
\end{align*}
Let $\varepsilon > 0$. Since $\partial R \in \Irm_{k}(\R^{d})$, we can appeal to Proposition \ref{pp:QBdryQStrict} to find $Q \in \Qrm_{k+1}(\R^{d})$ with
\begin{align*}
    \partial Q &= \partial R = T, & \Mbf(Q) &\leq \Mbf(R) + \varepsilon
\end{align*}
Take $n \in \N$ with $n \cdot Q \in \Irm_{k+1}(\R^{d})$. It follows that
\begin{align*}
    \FI^{\circ}(n \cdot T) \leq \Mbf(n \cdot Q) \leq n (\Mbf(R) + \varepsilon)
\end{align*}
and so
\begin{align*}
    \inf_{n \in \N} \frac{\FI^{\circ}(n \cdot T)}{n} \leq \FN^{\circ}(T) + \varepsilon.
\end{align*}
Since the sequence $(\FI^{\circ}(n \cdot T))_{n}$ is subadditive, the result follows from Fekete's lemma.
\end{proof}

\subsection{Space-Time Trajectories} In what follows, we present two examples of strict-type rational approximation results for space-time trajectories that are in some sense analogous to Proposition \ref{pp:NormalQBdryStrict} and Proposition \ref{pp:QBdryQStrict} respectively:

\begin{proposition}\label{pp:STQ} Take flat integrands
\begin{align*}
    \psi &: \R^{d} \times \subext{k}\R^{d} \to [0,\infty), & \phi &: \R^{d} \times \subext{1+k}\R^{d} \to [0,\infty). 
\end{align*}
Let $K \subset \R^{d}$ be compact and take $\varepsilon > 0$. Let $S \in \Nrm_{1+k}([0,1] \times \R^{d})$. Suppose that $\partial S \mvert (0,1) \times \R^{d} = 0$ and that $\supp S \subset [0,1] \times \Int K$. Then there is $Q \in \QPrm_{1+k}([0,1] \times \R^{d})$ such that $\partial Q \mvert (0,1) \times \R^{d} = 0$, $\supp Q \subset [0,1] \times K$, and the following conclusions hold:
\begin{enumerate}
    \item We have
    \begin{align*}
        \FN_{[0,1] \times K}(Q-S) &\leq \varepsilon
    \end{align*}
    and
    \begin{align*}
        \Var_{\phi}(Q) &\leq \Var_{\phi}(S) + \varepsilon, & \int \Mbf_{\psi}(Q(t)) \dd t \leq \int \Mbf_{\psi}(S(t)) \dd t + \varepsilon.
    \end{align*}
    \item We have
    \begin{align*}
        \FN_{K}(\partial^{\pm}Q - \partial^{\pm}S) &\leq \varepsilon, & \Mbf_{\psi}(\partial^{\pm}Q) &\leq \Mbf_{\psi}(\partial^{\pm}S) + \varepsilon.
    \end{align*}
    \item We have
    \begin{align*}
        \esssup_{t \in [0,1]} \Mbf_{\psi}(Q(t)) &\leq \esssup_{t \in [0,1]}\Mbf_{\psi}(S(t)) + \varepsilon.
    \end{align*}
\end{enumerate}
\end{proposition}

\begin{proof} We start by defining $\Xbf: \Prm_{1+k}(\R^{1+d}) \to [0,\infty)$ via
\begin{align}\label{STQ1}
    \Xbf(P) = \Var(P) + \Vert P \Vert_{\Lrm^{\infty}} + \Mbf(\partial P) \hspace{1em} \text{for $P \in \Prm_{1+k}(\R^{1+d})$}.
\end{align}
We shall show in Corollary \ref{cr:PolyLInfty} that $\Vert P \Vert_{\Lrm^{\infty}} < \infty$ for $P \in \Prm_{1+k}(\R^{1+d})$ and so $\Xbf$ is well-defined. One can show that $\Xbf$ is a seminorm on $\Prm_{1+k}(\R^{1+d})$.

Fix $\alpha > 0$ (to be chosen later). We apply Theorem \ref{tm:STPoly} to $S$ to find $P \in \Prm_{1+k}([0,1] \times \R^{d})$ such that $\supp P \subset [0,1] \times K$, $\partial P \mvert (0,1) \times \R^{d} = 0$ and conclusions (1)-(3) hold with $Q$ replaced by $P$ and $\varepsilon$ replaced with $\alpha$. We then apply Proposition \ref{pp:PolyQX} to $P$ with $A = \varnothing$ to obtain $Q \in \QPrm_{1+k}(\R^{1+d})$ with
\begin{align*}
    \supp Q & \subset \supp P \subset [0,1] \times K, & \supp \partial Q &\subset \supp \partial P \subset \{0,1\} \times K, 
\end{align*}
and
\begin{align*}
    \Xbf(Q-P) &\leq \alpha.
\end{align*}
It follows from \cite[Lemma 3.10]{BonicattoRindlerTurnbull26} that $\Mbf(Q-P) \leq \alpha$ and so
\begin{align}\label{STQ2}
    \Fbf_{[0,1] \times K}(Q-S) &\leq 2\alpha.
\end{align}
Notice that
\begin{align}\label{STQ3}
    \Var_{\phi}(Q) \leq \Var_{\phi}(S) + (C_{\phi,K}+1)\alpha
\end{align}
where $C_{\phi,K} > 0$ is defined as in \eqref{psiVsMass}. Similarly, we have
\begin{align}\label{STQ4}
    \int \Mbf_{\psi}(Q(t)) \dd t \leq \int \Mbf_{\psi}(S(t)) \dd t + (C_{\psi,K}+1)\alpha
\end{align}
and
\begin{align}\label{STQ5}
    \esssup_{t \in [0,1]} \Mbf_{\psi}(Q(t)) \leq \esssup_{t \in [0,1]} \Mbf_{\psi}(S(t)) + (C_{\psi,K}+1)\alpha
\end{align}
where $C_{\psi,K} > 0$ is defined as in \eqref{psiVsMass}. Finally, we have
\begin{align}\label{STQ6}
    \FN_{K}(\partial^{\pm}Q-\partial^{\pm}S) &\leq 2\alpha
\end{align}
and
\begin{align}\label{STQ7}
    \Mbf_{\psi}(\partial^{\pm}Q) &\leq \Mbf_{\psi}(\partial^{\pm}S) + (C_{\psi,K}+1)\alpha,
\end{align}
so the result follows after choosing $\alpha$ so that
\begin{align*}
    \alpha &\leq \varepsilon \cdot \min\bigg\{\frac{1}{2},\frac{1}{C_{\psi,K}+1},\frac{1}{C_{\phi,K}+1}\bigg\}. \qedhere
\end{align*}
\end{proof}

An abbreviated version of the following result was presented in the introduction as Theorem \ref{tm:AbbrvSTOneSidedQ}:

\begin{theorem}\label{tm:STOneSidedQ} Take flat integrands
\begin{align*}
    \psi &: \R^{d} \times \subext{k}\R^{d} \to [0,\infty), & \phi &: \R^{d} \times \subext{1+k}\R^{d} \to [0,\infty). 
\end{align*}
Let $K \subset \R^{d}$ be compact and take $\varepsilon > 0$. Let $S \in \Nrm_{1+k}([0,1] \times \R^{d})$. Suppose that $\partial S \mvert (0,1) \times \R^{d} = 0$, that $\supp S \subset [0,1] \times \Int K$, and that $\partial^{-}S \in \Qrm_{k}(\R^{d})$. Then there is $Q \in \Qrm_{1+k}([0,1] \times \R^{d})$ such that $\partial Q \mvert (0,1) \times \R^{d} = 0$, $\supp Q \subset [0,1] \times K$, conclusions (1)-(3) of Proposition \ref{pp:STQ} hold, and $\partial^{-}Q = \partial^{-}S$.
\end{theorem}

\begin{proof} Fix $\alpha > 0$ (to be chosen later). It is not restrictive to assume that $\partial^{-}S \in \Irm_{k}(\R^{d})$. We first apply Theorem \ref{tm:STPoly}, noting that we have conclusion (4) in this case, to find $P \in \Prm_{1+k}([0,1] \times \R^{d})$ such that $\supp P \subset [0,1] \times K$, $\partial P \mvert (0,1) \times \R^{d} = 0$, 
\begin{align*}
    \partial^{-}P &\in \IPrm_{k}(\R^{d}), & \FI^{\circ}_{K}(\partial^{-}P - \partial^{-}S) &\leq \alpha,
\end{align*}
and conclusions (1)-(3) of Proposition \ref{pp:STQ} hold with $Q$ replaced by $P$ and $\varepsilon$ replaced by $\alpha$. Define $\Xbf$ as in \eqref{STQ1}, and apply Proposition \ref{pp:PolyQX} to $P$ with $A = \{0\} \times \R^{d}$ to obtain $\tilde{Q} \in \QPrm_{1+k}(\R^{1+d})$ with $\supp \tilde{Q} \subset [0,1] \times K$, $\supp \partial\tilde{Q} \subset \{0,1\} \times \R^{d}$, and
\begin{align*}
    \Xbf(\tilde{Q}-P) &\leq \alpha, & \partial \tilde{Q} \mvert \{0\} \times \R^{d} &= \partial P \mvert \{0\} \times \R^{d}.
\end{align*}
It follows that $\partial^{-}\tilde{Q} = \partial^{-}P$. Also notice that the estimates \eqref{STQ2}-\eqref{STQ7} hold with $Q$ replaced by $\tilde{Q}$. Since $\FI^{\circ}_{K}(\partial^{-}P-\partial^{-}S) \leq \alpha$, we can find $R \in \Irm_{k+1}(\R^{d})$ with $\supp R \subset K$ and 
\begin{align*}
    \partial R &= \partial^{-}P - \partial^{-}S, & \Mbf(R) &\leq \alpha
\end{align*}
Now set
\begin{align}\label{STOneSidedQ1}
    Q := \tilde{Q} + \delta_{0} \times R.
\end{align}
Then $\supp Q \subset [0,1] \times K$, $\partial Q \mvert (0,1) \times \R^{d} = 0$, $\partial^{-}Q = \partial^{-}S$, and $\partial^{+}Q = \partial^{+}\tilde{Q}$. It follows from \eqref{STQ2} and \eqref{STOneSidedQ1} that
\begin{align*}
    \FN_{[0,1] \times K}(Q-S) &\leq 3\alpha.
\end{align*}
We can use \eqref{STQ3}, Lemma \ref{lm:DissJumps}, and \eqref{STOneSidedQ1} to deduce that
\begin{align*}
    \Var_{\phi}(Q) &\leq \Var_{\phi}(S) + (2C_{\phi,K}+1)\alpha
\end{align*}
where $C_{\phi,K} > 0$ is defined as in \eqref{psiVsMass}. Since $\Vert \delta_{0} \times R \Vert_{\Lrm^{\infty}} = 0$, the estimates \eqref{STQ4} and \eqref{STQ5} follow from \eqref{STOneSidedQ1}. Finally, since $\partial^{+}Q = \partial^{+}\tilde{Q}$, the estimates
\begin{align*}
    \FN_{K}(\partial^{\pm}Q-\partial^{\pm}S) &\leq 2\alpha
\end{align*}
and
\begin{align*}
    \Mbf_{\psi}(\partial^{\pm}Q) &\leq \Mbf_{\psi}(\partial^{\pm}S) + (C_{\psi,K}+1)\alpha
\end{align*}
follow immediately from \eqref{STQ6} and \eqref{STQ7} respectively. The result follows after choosing $\alpha$ so that 
\begin{align*}
    \alpha &\leq \varepsilon \cdot \min\bigg\{\frac{1}{3},\frac{1}{C_{\psi,K}+1},\frac{1}{2C_{\phi,K}+1}\bigg\}. \qedhere
\end{align*}
\end{proof}

\section{Lipschitz Approximation}\label{s:Lipschitz}

Let $S \in \Nrm_{1+k}([0,1] \times \R^{d})$ be a space-time trajectory, so that $\partial S \mvert (0,1) \times \R^{d} = 0$. Recall from \cite[Section 4]{BonicattoRindlerTurnbull26} that $S$ is said to be a \textbf{Lipschitz trajectory} if $\tbf_{\#}\Vert S \Vert(\{0,1\}) = 0$, $\Vert S \Vert_{\Lrm^{\infty}} < \infty$, and the function $t \mapsto \Var(S;(0,t))$ is Lipschitz on $(0,1)$. In the previous sections, the approximants of space-time trajectories that we have constructed have not been Lipschitz; in this section, we shall rectify this by proving strict-type \textit{Lipschitz} approximation results for space-time trajectories. In particular, we shall establish the following general form of Theorem \ref{tm:AbbrvSTLipschitz}:

\begin{theorem}\label{tm:STLipschitz} Take flat integrands
\begin{align*}
    \psi &: \R^{d} \times \subext{k}\R^{d} \to [0,\infty), & \phi &: \R^{d} \times \subext{1+k}\R^{d} \to [0,\infty).
\end{align*}
Let $K \subset \R^{d}$ be compact and take $\varepsilon > 0$. Let $S \in \Nrm_{1+k}([0,1] \times \R^{d})$. Suppose that $\partial S \mvert (0,1) \times \R^{d} = 0$ and that $\supp S \subset [0,1] \times \Int K$. Then there is $V \in \NLip_{1+k}([0,1] \times \R^{d})$ such that $\partial V = \partial S$, $\supp V \subset [0,1] \times K$,
\begin{align*}
        \FN_{[0,1] \times K}(V-S) &\leq \varepsilon,
\end{align*}
and
\begin{align*}
    \Var_{\phi}(V) &\leq \Var_{\phi}(S) + \varepsilon, & \int \Mbf_{\psi}(V(t)) \dd t \leq \int \Mbf_{\psi}(S(t)) \dd t + \varepsilon.
\end{align*}
If $S \in \Irm_{1+k}([0,1] \times \R^{d})$, then we can additionally ensure that $V \in \ILip_{1+k}([0,1] \times \R^{d})$.
\end{theorem}

Notice that, unlike the previous strict-type approximation results that we have established for space-time trajectories, we here do not have uniform control over the slices of the approximants. 

\subsection{Polyhedral Trajectories} As in Section \ref{s:Q}, we shall first consider the case in which the space-time trajectory to be approximated is also a polyhedral chain. In particular, we shall approximate a polyhedral trajectory by a polyhedral Lipschitz trajectory. Since polyhedral chains are simple, it suffices to consider strict approximation results in this setting (see Proposition \ref{pp:SimpleSTS}).

We start with an auxiliary result which will allow us to discuss elements of $\Prm_{1+k}(\R^{1+d})$ within the context of the space-time framework developed in \cite{BonicattoRindlerTurnbull26,Rindler23}. 

\begin{lemma}\label{lm:SimpSlant} Take $v_{0},\dotsc,v_{1+k} \in \R^{1+d}$ and suppose that there are $i,j \in \{0,\dotsc,1+k\}$ with $\tbf(v_{i} - v_{j}) \neq 0$. Set $P := [\![v_{0},\dotsc,v_{1+k}]\!]$. Then there is $M \geq 0$ such that 
\begin{align*}
    \tbf_{\#}\Vert P \Vert([s,t]) \leq M \cdot \vert t - s \vert \hspace{1em} \text{for $s,t \in \R$ with $s \leq t$}.
\end{align*}
\end{lemma}

\begin{proof} First notice that $\Vert P \Vert = \mathscr{H}^{k+1} \mvert C$, where $C \subset \R^{1+d}$ is the convex hull of $\{v_{0},\dotsc,v_{1+k}\}$. Take a $(1+k)$-dimensional affine subspace $\Pi \subset \R^{1+d}$ with $C \subset \Pi$. Since $\tbf(v_{i}-v_{j}) \neq 0$, we must have $\tbf(\Pi) = \R$. It follows that there is an affine bijection $\Phi: \R \times \R^{k} \to \Pi$ with $\tbf(\Phi(t,x))=t$ for $(t,x) \in \R \times \R^{k}$. Fix $s,t \in \R$ with $s \leq t$. Then
\begin{align}\label{SimpSlant1}
    ([s,t] \times \R^{d}) \cap \Pi = \Phi([s,t] \times \R^{k}). 
\end{align}
Let $\tilde{\tbf}: \R \times \R^{k} \to \R$ and $\tilde{\pbf}: \R \times \R^{k} \to \R^{k}$ be the orthogonal projections onto the first and second components of the product $\R \times \R^{k}$ respectively. Notice that 
\begin{align}\label{SimpSlant2}
    C \subset \Phi(\tilde{\tbf}(\Phi^{-1}(C)) \times \tilde{\pbf}(\Phi^{-1}(C)))
\end{align}
Combining \eqref{SimpSlant1} and \eqref{SimpSlant2} yields
\begin{align*}
    ([s,t] \times \R^{d}) \cap C \subset \Phi([s,t] \times \tilde{\pbf}(\Phi^{-1}(C))).
\end{align*}
Since $\Phi$ is Lipschitz, we have
\begin{align*}
    \mathscr{H}^{1+k}(([s,t] \times \R^{d}) \cap C) \leq \Lip(\Phi)^{1+k}\cdot \mathscr{L}^{1+k}([s,t] \times \tilde{\pbf}(\Phi^{-1}(C)))
\end{align*}
and so
\begin{align*}
    \mathscr{H}^{1+k}(([s,t] \times \R^{d}) \cap C) \leq \Lip(\Phi)^{1+k}\cdot \vert t - s \vert \cdot \mathscr{L}^{k}(\tilde{\pbf}(\Phi^{-1}(C))).
\end{align*}
The result follows with $M = \Lip(\Phi)^{1+k} \cdot \mathscr{L}^{k}(\tilde{\pbf}(\Phi^{-1}(C)))$.
\end{proof}

\begin{corollary}\label{cr:PolyLInfty} Take $P \in \Prm_{1+k}(\R^{1+d})$. Then 
\begin{align*}
    \Vert P \Vert_{\Lrm^{\infty}} < \infty.
\end{align*}
\end{corollary}

\begin{proof} Take $v_{0},\dotsc,v_{1+k} \in \R^{1+d}$ and suppose that $P = [\![v_{0},\dotsc,v_{1+k}]\!]$. If $\tbf(v_{i}-v_{j}) = 0$ for all $i,j \in \{0,\dotsc,1+k\}$, then there is $t_{0} \in \R$ with $\tbf(\supp P) \subset \{t_{0}\}$ and so $\Vert P \Vert_{\Lrm^{\infty}} = 0$. Otherwise, we can use Lemma \ref{lm:SimpSlant} to find $M \geq 0$ with $\tbf_{\#}\Vert P \Vert([s,t]) \leq M \cdot \vert t - s \vert$ for $s,t \in \R$ with $s \leq t$. It follows that $\Theta^{1*}(\tbf_{\#}\Vert P \Vert,t) \leq M$ for all $t \in \R$. Since $\Mbf(P\vert_{t}) \leq \Theta^{*}(\tbf_{\#}\Vert P \Vert,t)$ for $\mathscr{L}^{1}$-a.e.\ $t \in \R$ (see \cite[Section 3.1]{BonicattoRindlerTurnbull26}), we have $\Vert P \Vert_{\Lrm^{\infty}} \leq M < \infty$. In the general case, the result follows from writing $P$ as a finite sum of oriented simplices. 
\end{proof}

Take a polyhedral trajectory $P \in \Prm_{1+k}([0,1] \times \R^{d})$, so that $\partial P \mvert (0,1) \times \R^{d} = 0$. It follows from Corollary \ref{cr:PolyLInfty} that $P$ is a Lipschitz trajectory if and only if $\tbf_{\#}\Vert P \Vert(\{0,1\}) = 0$ and the function $t \mapsto \Var(P;(0,t))$ is Lipschitz on $(0,1)$. We shall use this observation to prove the main result of this subsection:

\begin{proposition}\label{pp:PolyLipschitz} Let $K \subset \R^{d}$ be compact and take $\varepsilon > 0$. Let $P \in \Prm_{1+k}([0,1] \times \R^{d})$. Suppose that $\supp P \subset [0,1] \times K$ and that $\partial P \mvert (0,1) \times \R^{d} = 0$. Then there are 
\begin{align*}
    V \in \PLip_{1+k}([0,1] \times \R^{d}) &:= \Prm_{1+k}([0,1] \times \R^{d}) \cap \NLip_{1+k}([0,1] \times \R^{d})
\end{align*}
and $X \in \Prm_{k+2}([0,1] \times \R^{d})$
with $\supp V$, $\supp X \subset [0,1] \times K$ and 
\begin{align}\label{cc:PolyLipschitz}
    P &= V + \partial X, & \Mbf(V) &\leq \Mbf(P) + \varepsilon, & \Mbf(X) &\leq \varepsilon.
\end{align}
If $P \in \IPrm_{1+k}([0,1] \times \R^{d})$, then we can additionally ensure that $X \in \IPrm_{k+2}([0,1] \times \R^{d})$ and
\begin{align*}
    V \in \IPLip_{1+k}([0,1] \times \R^{d}) := \IPrm_{1+k}([0,1] \times \R^{d}) \cap \NLip_{1+k}([0,1] \times \R^{d}).
\end{align*}
\end{proposition}

We start by working with individual simplices (which may not be space-time trajectories):

\begin{lemma}\label{lm:SimpLipschitz} Let $K \subset \R^{d}$ be compact and take $\varepsilon > 0$. Take $v_{0},\dotsc,v_{1+k} \in \R^{1+d}$ and set $P = [\![v_{0},\dotsc,v_{1+k}]\!]$. Suppose that $\supp P \subset [0,1] \times K$. Then there are $V \in \IPrm_{1+k}(\R^{1+d})$ and $X \in \IPrm_{k+2}(\R^{1+d})$ with $\supp V$, $\supp X \subset [0,1] \times K$,
\begin{align}\label{cc:SimpLipschitz}
    P &= V + \partial X, & \Mbf(V) &\leq \Mbf(P) + \varepsilon, & \Mbf(X) &\leq \varepsilon
\end{align}
such that $t \mapsto \Var(V;(0,t))$ is Lipschitz on $(0,1)$ and $\tbf_{\#}\Vert V \Vert(\{t\}) = 0$ for $t \in \R$.
\end{lemma}

\begin{proof} If $P=0$, then we can take $V=0$ and $X=0$ and there is nothing to show. If there are $i,j \in \{0,\dotsc,1+k\}$ with $\tbf(v_{i}-v_{j}) \neq 0$. In this case, we can appeal to Lemma \ref{lm:SimpSlant} to find $M \geq 0$ with $\tbf_{\#}\Vert P \Vert([s,t]) \leq M \cdot \vert t - s \vert$ for $s,t \in \R$ with $s \leq t$. Then $\tbf_{\#}\Vert P \Vert(\{t\}) = 0$ for $t \in \R$ and
\begin{align*}
    \Var(P;(s,t)) \leq \Mbf(P \mvert (s,t) \times \R^{d}) \leq M \cdot \vert t - s \vert \hspace{1em} \text{for $s,t \in (0,1)$},
\end{align*}
so we can take $V=P$ and $X=0$ in this case. 

Otherwise, there is $\sigma \in [0,1]$ with 
\begin{align*}
    P = \bigg[\!\!\bigg[ \binom{\sigma}{\pbf(v_{0})}, \dotsc, \binom{\sigma}{\pbf(v_{1+k})} \bigg]\!\!\bigg].
\end{align*}
Fix $\tau \in [0,1]$ with $\sigma \neq \tau$ and define
\begin{align*}
    \tilde{X} := \bigg[\!\!\bigg[ \binom{\tau}{\pbf(v_{0})}, \binom{\sigma}{\pbf(v_{0})}, \dotsc, \binom{\sigma}{\pbf(v_{1+k})} \bigg]\!\!\bigg].
\end{align*}
Then $\tilde{X} \in \IPrm_{k+2}(\R^{1+d})$ and $\supp \tilde{X} \subset I \times K$, where $I \subset [0,1]$ is the smallest closed interval containing $\{\sigma,\tau\}$. Notice that
\begin{align}\label{SimpLipschitz1}
    \partial \tilde{X} = P - \sum_{j=0}^{k+1} (-1)^{j}\cdot\tilde{V}_{j}
\end{align}
where, for $j \in \{0,\dotsc,k+1\}$, $\tilde{V}_{j} \in \IPrm_{k+1}(\R^{1+d})$ is given by
\begin{align*}
    \tilde{V}_{j} := \bigg[\!\!\bigg[ \binom{\tau}{\pbf(v_{0})}, \binom{\sigma}{\pbf(v_{0})}, \dotsc,\binom{\sigma}{\pbf(v_{j-1})}, \binom{\sigma}{\pbf(v_{j+1})}, \dotsc, \binom{\sigma}{\pbf(v_{k+1})} \bigg]\!\!\bigg].
\end{align*}
Now set 
\begin{align*}
    \tilde{V} := \sum_{j=0}^{k+1} (-1)^{j}\cdot\tilde{V}_{j}.
\end{align*}
Then $\tilde{V} \in \IPrm_{k+1}(\R^{1+d})$ and $\supp \tilde{V} \subset I \times K$, and it follows from \eqref{SimpLipschitz1} that $P = \tilde{V} + \partial \tilde{X}$. 

Take $j \in \{0,\dotsc,k+1\}$ and notice that
\begin{align*}
    \tbf\bigg(\binom{\tau}{\pbf(v_{0})}-\binom{\sigma}{\pbf(v_{i})}\bigg) \neq 0 \hspace{1em} \text{for $i \in \{0,\dotsc,k+1\}$ with $i \neq j$}.
\end{align*}
We can therefore appeal to Lemma \ref{lm:SimpSlant} to find $M_{j} \geq 0$ with $\tbf_{\#}\Vert \tilde{V}_{j} \Vert([s,t]) \leq M_{j} \cdot \vert t - s \vert$ for $s,t \in \R$ with $s \leq t$. Then, as above, we have $\tbf_{\#}\Vert \tilde{V}_{j} \Vert(\{t\}) = 0$ for $t \in \R$ and 
\begin{align*}
    \Var(\tilde{V}_{j};(s,t)) \leq M_{j} \cdot \vert t - s \vert \hspace{1em} \text{for $s,t \in (0,1)$}.
\end{align*}
It follows that $\tbf_{\#}\Vert \tilde{V} \Vert(\{t\}) = 0$ for $t \in \R$ and that the map $t \mapsto \Var(\tilde{V};(0,t))$ is Lipschitz on $\R$. Now notice that $\pbf_{*}\tilde{X} = 0$, that $\pbf_{*}\tilde{V}_{j} = 0$ for $j \in \{1,\dotsc,k+1\}$, and that $\pbf_{*}\tilde{V}_{0} = \pbf_{*}P$. We can therefore appeal to \cite[Lemma 4.3]{BonicattoRindlerTurnbull26} to deduce that
\begin{align*}
    \Var(\tilde{X}) &= 0, & \Var(\tilde{V}) &\leq \Mbf(\pbf_{*}P).
\end{align*}
Since $P = \delta_\sigma \times \pbf_{*}P$, we have $\Mbf(\pbf_{*}P) = \Mbf(P)$ and so $\Var(\tilde{V}) \leq \Mbf(P)$. Now take $\delta \in (0,1)$ and define $a: \R \to \R$ by
\begin{align*}
    a(t) := \sigma + \delta(t-\sigma) \hspace{1em} \text{for $t \in \R$}.
\end{align*}
Set
\begin{align*}
    V &:= a_{*}\tilde{V}, & X &:= a_{*}\tilde{X}.
\end{align*}
Since $a$ is affine, we have $V \in \IPrm_{1+k}(\R^{1+d})$ and $X \in \IPrm_{k+2}(\R^{1+d})$. Since $\tbf_{\#}\Vert \tilde{V} \Vert(\{t\}) = 0$ for $t \in \R$, we have $\tbf_{\#}\Vert V \Vert(\{t\}) = 0$ for $t \in \R$. Since $\Lip(a) = \delta$ and the map $t \mapsto \Var(\tilde{V};(0,t))$ is Lipschitz on $\R$, we can use Proposition \ref{pp:DissRescale} to deduce that the map $t \mapsto \Var(V;(0,t))$ is Lipschitz on $\R$. Since $a(\sigma)=\sigma$, we have $a_{*}P = P$ and so 
\begin{align*}
    P = V + \partial X.
\end{align*}
Also notice that $\supp V$, $\supp X \subset a(I) \times K$; since $\delta \in (0,1)$, we have $a(I) \subset I$ and so $\supp X$, $\supp V \subset [0,1] \times K$. It remains to estimate $\Mbf(X)$ and $\Mbf(V)$. To this end, we use Proposition \ref{pp:DissRescale} to deduce that
\begin{align*}
    \Var(V) = \Var(\tilde{V}) &\leq \Mbf(P), & \Var(X) &= \Var(\tilde{X}) = 0
\end{align*}
and \cite[Proposition 3.6]{BonicattoRindlerTurnbull26} to deduce that
\begin{align*}
    \Vert V \Vert_{\Lrm^{\infty}} &= \Vert \tilde{V} \Vert_{\Lrm^{\infty}}, & \Vert X \Vert_{\Lrm^{\infty}} &= \Vert \tilde{X} \Vert_{\Lrm^{\infty}}.
\end{align*}
Since $\supp V$, $\supp X \subset a(I) \times \R^{d}$, we can use \cite[Lemma 3.10]{BonicattoRindlerTurnbull26} to estimate
\begin{align*}
    \Mbf(V) &\leq \Mbf(P) + \mathscr{L}^{1}(a(I)) \cdot \Vert V \Vert_{\Lrm^{\infty}}, & \Mbf(X) &\leq \mathscr{L}^{1}(a(I)) \cdot \Vert X \Vert_{\Lrm^{\infty}}.
\end{align*}
We therefore have
\begin{align*}
    \Mbf(V) &\leq \Mbf(P) + \delta\cdot\mathscr{L}^{1}(I) \cdot \Vert V \Vert_{\Lrm^{\infty}}, & \Mbf(X) &\leq \delta\cdot\mathscr{L}^{1}(I) \cdot \Vert X \Vert_{\Lrm^{\infty}}.
\end{align*}
The result follows after choosing $\delta$ so that
\begin{align*}
    \delta \cdot \mathscr{L}^{1}(I) \cdot \max\{\Vert V \Vert_{\Lrm^{\infty}},\Vert X \Vert_{\Lrm^{\infty}}\} &\leq \varepsilon \qedhere.
\end{align*}
\end{proof}

\begin{proof}[Proof of Proposition \ref{pp:PolyLipschitz}] Since $P \in \Prm_{1+k}(\R^{1+d})$, there is $N \in \N$, $\alpha \in \R^{N}$ with $\alpha_{j} \neq 0$ for $j \in \{1,\dotsc,N\}$, and a collection $\{P_{j}\}_{j=1}^{N}$ of $\mathscr{H}^{1+k}$-disjoint oriented $(1+k)$-simplices in $\R^{1+d}$ such that 
\begin{align*}
    P = \sum_{j=1}^{N} \alpha_{j}P_{j}.
\end{align*}
Fix $j \in \{1,\dotsc,N\}$ and notice that $\supp P_{j} \subset [0,1] \times K$. We can therefore appeal to Lemma \ref{lm:SimpLipschitz} to deduce that there are $V_{j} \in \IPrm_{1+k}(\R^{1+d})$ and $X_{j} \in \IPrm_{k+2}(\R^{1+d})$ with $\supp V_{j}$, $\supp X_{j} \subset [0,1] \times K$,
\begin{align}\label{PolyLipschitz}
    P_{j} &= V_{j} + \partial X_{j}, & \Mbf(V_{j}) &\leq \Mbf(P_{j}) + \frac{\varepsilon}{\vert \alpha_{j} \vert\cdot N}, & \Mbf(X_{j}) &\leq \frac{\varepsilon}{\vert \alpha_{j} \vert\cdot N}
\end{align}
such that $\tbf_{\#}\Vert V_{j} \Vert(\{0,1\}) = 0$ and the map $t \mapsto \Var(V_{j};(0,t))$ is Lipschitz on $(0,1)$. Set 
\begin{align*}
    V &:= \sum_{j=1}^{N}\alpha_{j}V_{j}, & X &:= \sum_{j=1}^{N}\alpha_{j}X_{j}.
\end{align*}
Then $V \in \Prm_{1+k}(\R^{1+d})$, $X \in \Prm_{k+2}(\R^{1+d})$, and $\supp V$, $\supp X \subset [0,1] \times K$. Furthermore, \eqref{cc:PolyLipschitz} follows from \eqref{PolyLipschitz} and the fact that the collection $\{P_{j}\}_{j=1}^{N}$ is $\mathscr{H}^{k+1}$-disjoint. One can also check that $\tbf_{\#}\Vert V \Vert(\{0,1\}) = 0$ and that the map $t \mapsto \Var(V;(0,t))$ is Lipschitz on $(0,1)$. Since $V \in \Prm_{1+k}(\R^{1+d})$, we can use Corollary \ref{cr:PolyLInfty} to deduce that $\Vert V \Vert_{\Lrm^{\infty}} < \infty$, and since $\partial V = \partial P$, we have $\partial V \mvert (0,1) \times \R^{d} = \partial P \mvert (0,1) \times \R^{d} = 0$. We therefore have
\begin{align*}
    V \in \PLip_{k+1}([0,1] \times \R^{d}).
\end{align*}
Finally, if $P \in \IPrm_{1+k}(\R^{1+d})$, then notice that $\alpha \in \Z^{N}\backslash\{0\}$ and so $V \in \IPrm_{1+k}(\R^{1+d})$ and $X \in \IPrm_{k+2}(\R^{1+d})$.
\end{proof}

\subsection{Proof of Theorem \ref{tm:STLipschitz}} To prove Theorem \ref{tm:STLipschitz}, we shall need the following `gluing' lemma:

\begin{lemma}\label{lm:Gluing} Take flat integrands
\begin{align*}
    \psi &: \R^{d} \times \subext{k}\R^{d} \to [0,\infty), & \phi &: \R^{d} \times \subext{1+k}\R^{d} \to [0,\infty). 
\end{align*}
Let $K \subset \R^{d}$ be compact and take $\varepsilon > 0$. Let $S$, $S^{\pm} \in \NLip_{1+k}([0,1] \times K)$. Suppose that $\partial S \mvert (0,1) \times \R^{d} = \partial S^{\pm} \mvert (0,1) \times \R^{d} = 0$ and
\begin{align*}
    \partial^{\pm}S &= \partial^{\mp}S^{\pm}, & \Var(S^{\pm}) &\leq \varepsilon.
\end{align*}
Then there is $V \in \NLip_{1+k}([0,1] \times K)$ with $\partial^{\pm}V = \partial^{\pm}S^{\pm}$,
\begin{align*}
    \Fbf_{[0,1] \times K}(V-S) &\leq 3\varepsilon,
\end{align*}
and
\begin{align*}
    \Var_{\phi}(V) &\leq \Var_{\phi}(S) + 2C_{\phi,K}\cdot\varepsilon, & \int \Mbf_{\psi}(V(t)) \dd t \leq \int \Mbf_{\psi}(S(t)) \dd t + \varepsilon,
\end{align*}
where $C_{\phi,K} > 0$ is defined as in \eqref{psiVsMass}. If $S,S^{\pm} \in \ILip_{1+k}([0,1] \times K)$, then we can additionally ensure that $V \in \ILip_{1+k}([0,1] \times K)$. 
\end{lemma}

\begin{proof} Fix $\alpha \in (0,1/2)$ (to be chosen later), and define $a: [0,1] \to [0,1]$ and $a^{\pm}: [0,1] \to [0,1]$ via, for $t \in [0,1]$,
\begin{align*}
    a(t) &:= (1-2\alpha)t + \alpha, & a^{-}(t) &:= \alpha t, & a^{+}(t) := \alpha t + (1-\alpha).
\end{align*}
Set
\begin{align}\label{Gluing1}
    V := (a^{-})_{*}S^{-} + a_{*}S + (a^{+})_{*}S^{+}.
\end{align}
Since $S,S^{\pm} \in \NLip_{1+k}([0,1] \times K)$, one can check that $V \in \NLip_{1+k}([0,1] \times K)$. Since $\partial^{\pm}S = \partial^{\mp}S^{\pm}$, we have 
\begin{align*}
    \partial V &= (a^{-})_{*}(\delta_{1} \times \partial^{-}S - \delta_{0} \times \partial^{-}S^{-}) \\
    &\qquad + a_{*}(\delta_{1} \times \partial^{+}S - \delta_{0} \times \partial^{-}S) \\
    &\qquad + (a^{+})_{*}(\delta_{1} \times \partial^{+}S^{+} - \delta_{0} \times \partial^{+}S) \\
    &= \delta_{1} \times \partial^{+}S^{+} - \delta_{0} \times \partial^{-}S^{-}.
\end{align*}
We now appeal to Proposition \ref{pp:DissRescale} and \cite[Proposition 3.6]{BonicattoRindlerTurnbull26} to deduce that 
\begin{align}\label{Gluing1.5}
    \Var_{\phi}((a^{\pm})_{*}S^{\pm}) &= \Var_{\phi}(S^{\pm}), & \Vert (a^{\pm})_{*}S^{\pm} \Vert_{\Lrm^{\infty}} &= \Vert S^{\pm} \Vert_{\Lrm^{\infty}}.
\end{align}
Notice that $\supp (a^{\pm})_{*}S^{\pm} \subset a^{\pm}([0,1]) \times K$; we then use \cite[Lemma 3.10]{BonicattoRindlerTurnbull26} and the fact that $\diam a^{\pm}([0,1]) = \alpha$ to estimate
\begin{align}\label{Gluing2}
    \Mbf((a^{\pm})_{*}S^{\pm}) &\leq \Var(S^{\pm}) + \alpha \cdot \Vert S^{\pm}\Vert_{\Lrm^{\infty}} \leq \varepsilon + \alpha \cdot \Vert S^{\pm}\Vert_{\Lrm^{\infty}}.
\end{align}
Since $\Lip(a \times \mathrm{id}_{\R^{d}}) = 1$ and $\Vert a \times \mathrm{id}_{\R^{d}} - \mathrm{id}_{[0,1] \times \R^{d}} \Vert_{\Lrm^{\infty}} \leq \alpha$, we can estimate
\begin{align}\label{Gluing3}
    \FN_{[0,1] \times K}(a_{*}S-S) &\leq \alpha \cdot \Nbf(S).
\end{align}
Combining \eqref{Gluing2} and \eqref{Gluing3} yields
\begin{align*}
    \FN_{[0,1] \times K}(V-S) \leq \alpha(\Nbf(S) + \Vert S^{-} \Vert_{\Lrm^{\infty}} + \Vert S^{+} \Vert_{\Lrm^{\infty}}) + 2\varepsilon.
\end{align*}
We then use \eqref{Gluing1.5}, Proposition \ref{pp:DissRescale}, and the fact that $\Var(S^{\pm}) \leq \varepsilon$ to obtain
\begin{align*}
    \Var_{\phi}(V) \leq \Var_{\phi}(S) + 2C_{\phi,K}\cdot \varepsilon
\end{align*}
where $C_{\phi,K}$ is defined as in \eqref{psiVsMass}. We then appeal to \cite[Proposition 3.6]{BonicattoRindlerTurnbull26} once more to estimate
\begin{align*}
    \int \Mbf_{\psi}(V(t)) \dd t &\leq \int_{\alpha}^{1-\alpha} \Mbf_{\psi}(S(a^{-1}(t))) \dd t + C_{\psi,K}\cdot \alpha(\Vert S^{-} \Vert_{\Lrm^{\infty}}+\Vert S^{+} \Vert_{\Lrm^{\infty}}) \\
    &= (1-2\alpha)\int \Mbf_{\psi}(S(t)) \dd t + C_{\psi,K}\cdot \alpha(\Vert S^{-} \Vert_{\Lrm^{\infty}}+\Vert S^{+} \Vert_{\Lrm^{\infty}}).
\end{align*}
Now choose $\alpha$ so that
\begin{align*}
    \alpha \cdot \max\{\Nbf(S)+\Vert S^{-} \Vert_{\Lrm^{\infty}}+\Vert S^{+} \Vert_{\Lrm^{\infty}},C_{\psi,K}\cdot(\Vert S^{-} \Vert_{\Lrm^{\infty}}+\Vert S^{+} \Vert_{\Lrm^{\infty}})\} &\leq \varepsilon. 
\end{align*}
Finally, if $S,S^{\pm} \in \ILip_{1+k}([0,1] \times K)$, then it follows from \eqref{Gluing1} that $V \in \ILip_{1+k}([0,1] \times K)$. 
\end{proof}

\begin{proof}[Proof of Theorem \ref{tm:STLipschitz}] Fix $\alpha > 0$ (to be chosen later), and take an intermediate compact set $A \subset \Int K$ with $\supp S \subset [0,1] \times \Int A$. We first apply Theorem \ref{tm:STPoly} to $S$ to obtain $P \in \Prm_{1+k}([0,1] \times \R^{d})$ such that $\supp P \subset [0,1] \times A$, $\partial P \mvert (0,1) \times \R^{d} = 0$,
\begin{align*}
    \FN_{[0,1] \times A}(P-S) &\leq \alpha, & \FN^{\circ}_{A}(\partial^{\pm}P-\partial^{\pm}S) &\leq \alpha,
\end{align*}
and
\begin{align*}
    \Var_{\phi}(P) &\leq \Var_{\phi}(S) + \alpha, & \int \Mbf_{\psi}(P(t)) \dd t &\leq \int \Mbf_{\psi}(S(t)) \dd t + \alpha.
\end{align*}
We then use Proposition \ref{pp:PolyLipschitz} together with Proposition \ref{pp:SimpleSTS} to find $V_{P} \in \PLip_{1+k}([0,1] \times \R^{d})$ with $\supp V_{P} \subset [0,1] \times A$, $\partial V_{P} = \partial P$, 
\begin{align*}
    \FN_{[0,1] \times A}(V_{P}-P) &\leq \alpha
\end{align*}
and
\begin{align*}
    \Var_{\phi}(V_{P}) &\leq \Var_{\phi}(P) + \alpha, & \int \Mbf_{\psi}(V_{P}(t)) \dd t &\leq \int \Mbf_{\psi}(P(t)) \dd t + \alpha.
\end{align*}
Since $\FN^{\circ}_{A}(\partial^{\pm}P-\partial^{\pm}S) \leq \alpha$, we can use \cite[Lemma 5.1]{BonicattoRindlerTurnbull26} to construct $S^{\pm} \in \NLip_{1+k}([0,1] \times \R^{d})$ with $\supp S^{\pm} \subset [0,1] \times K$, $\Var(S^{\pm}) \leq 2\alpha$, and 
\begin{align*}
    \partial S^{+} &= \delta_{1} \times \partial^{+}S - \delta_{0} \times \partial^{+}P, \\
    \partial S^{-} &= \delta_{1} \times \partial^{-}P - \delta_{0} \times \partial^{-}S.
\end{align*}
Since $\Var(S^{\pm}) \leq 2\alpha$, we can appeal to Lemma \ref{lm:Gluing} (with $2\alpha$ in place of $\varepsilon$) to obtain $V \in \NLip_{1+k}([0,1] \times K)$ with $\partial^{\pm}V = \partial^{\pm}S^{\pm} = \partial^{\pm}S$, 
\begin{align*}
    \FN_{[0,1] \times K}(V-V_{P}) &\leq 6\alpha,
\end{align*}
and
\begin{align*}
    \Var_{\phi}(V) &\leq \Var_{\phi}(V_{P}) + 4C_{\phi,K}\cdot\alpha, & \int \Mbf_{\psi}(V(t)) \dd t \leq \int \Mbf_{\psi}(V_{P}(t)) \dd t + 2\alpha.
\end{align*}
Now choose $\alpha$ so that
\begin{align*}
    \alpha &\leq \varepsilon \cdot \min\bigg\{\frac{1}{8},\frac{1}{2(1+2C_{\phi,K})}\bigg\}.
\end{align*}
Combining the estimates above yields the result. If $S \in \Irm_{1+k}([0,1] \times \R^{d})$, then we follow the argument above, noting that conclusions (4) and (5) of Theorem \ref{tm:STPoly}, the additional conclusion of Proposition \ref{pp:PolyLipschitz}, the integral version of \cite[Lemma 5.1]{BonicattoRindlerTurnbull26} (see also \cite[Theorem 5.4]{Rindler23}), and the additional conclusion of Lemma \ref{lm:Gluing} hold in this case.
\end{proof}

\section{Results with Controlled Support}\label{s:Support}

Many of the approximation results that we have discussed above share the property that the support of the approximant may grow in some way. In what follows, we shall demonstrate how to rectify this if we assume that the supporting set has some additional regularity. In particular, we consider compact sets with the following property:

\begin{definition}\label{df:PullIn} Let $K \subset \R^{d}$ be compact. We say that $K$ has the \textbf{pull-in property} if, for $\varepsilon > 0$, there is a Lipschitz $f_{\varepsilon,K}: K \to \Int K$ such that
\begin{align*}
    \Lip(f_{\varepsilon,K}-\mathrm{id}_{K}) &\leq \varepsilon, & \Vert f_{\varepsilon,K} -\mathrm{id}_{K} \Vert_{\Lrm^{\infty}} &\leq \varepsilon
\end{align*}
and the affine homotopy $H_{\varepsilon,K}: [0,1] \times K \to \R^{d}$ from $\mathrm{id}_{K}$ to $f_{\varepsilon,K}$ satisfies
\begin{align*}
    H_{\varepsilon,K}([0,1] \times K) \subset K.
\end{align*}
\end{definition}

There are a wide variety of compact sets $K \subset \R^{d}$ with the pull-in property: For example, it follows from \cite[Proposition 2.10]{BonicattoRindlerTurnbull26} that the closure of any bounded $\Crm^{0}$ domain, and hence the closure of any bounded Lipschitz domain, has this property.

Before proceeding, we record the following elementary estimate:

\begin{lemma}\label{lm:MultiEst} Take linear maps $A,L: \R^{d} \to \R^{d}$. Then
\begin{align}\label{MultiEst}
    \Vert (A+L)\xi - A\xi \Vert \leq ((\Vert A \Vert+\Vert L \Vert)^{k}-\Vert A \Vert^{k})\cdot \Vert \xi \Vert \hspace{1em} \text{for $\xi \in \subext{k}\R^{d}$}.
\end{align}
\end{lemma}

\begin{proof} We proceed by induction on $k$. If $k=0$, then the estimate \eqref{MultiEst} follows immediately from the observation that $(A + L)\xi - \xi = \xi - \xi = 0$ for $\xi \in \subext{0}\R^{d}$. Otherwise, suppose that
\begin{align*}
    \Vert (A+L)\eta - A\eta \Vert \leq ((\Vert A \Vert + \Vert L \Vert)^{k-1}-\Vert A \Vert^{k-1})\cdot \Vert \eta \Vert \hspace{1em} \text{for $\eta \in \subext{k-1}\R^{d}$}.
\end{align*}
Take $\xi \in \subext{k}\R^{d}$. First suppose that $\xi$ is simple; without loss of generality, we can ensure that $\xi = \eta \wedge v$ for some $\eta \in \subext{k-1}\R^{d}$ and $v \in \R^{d}$ with $\mathrm{span}\hspace{0.2em} \eta \perp v$. Then
\begin{align*}
    (A+L)\xi - A\xi &= ((A+L)\eta) \wedge (Av + Lv) - A\eta \wedge Av \\
    &= ((A+L)\eta - A\eta) \wedge (Av+Lv) + A\eta \wedge Lv
\end{align*}
and we can estimate
\begin{align}\label{MultiEst1}
\begin{aligned}
    \Vert (A+L)\xi - A\xi \Vert &\leq ((\Vert A \Vert+\Vert L \Vert)^{k-1}-\Vert A \Vert^{k-1}) \Vert \eta \Vert (\Vert A \Vert+\Vert L \Vert) \Vert v \Vert \\ &+ \Vert A \Vert^{k-1}\Vert \eta \Vert \Vert L \Vert \Vert v \Vert.
\end{aligned}
\end{align}
Since $\mathrm{span}\hspace{0.2em}\eta \perp v$, we have $\Vert \xi \Vert = \Vert\eta \Vert \Vert v \Vert$ (see \cite[1.8.4]{Federer69book}) and so \eqref{MultiEst1} reduces to \eqref{MultiEst}. In the general case, the estimate follows from the fact that $\Vert \cdot \Vert^{\flat} = \Vert \cdot \Vert$ (see \cite[1.8.1]{Federer69book}).
\end{proof}

\subsection{Normal and Integral Currents} We shall need the following lemma:

\begin{lemma}\label{lm:PullInto} Let $\psi: \R^{d} \times \subext{k}\R^{d} \to [0,\infty)$ be a positive convex integrand. Take a compact set $K \subset \R^{d}$ with the pull-in property and let $\varepsilon > 0$. Then there is $\eta > 0$ with the following property: For $T \in \Nrm_{k}(K)$ and $\delta > 0$ with $\delta \leq \eta$, we have
\begin{align}\label{cc:PullInto}
    \Mbf_{\psi}((f_{\delta,K})_{*}T) &\leq (1+\varepsilon)\Mbf_{\psi}(T)
\end{align}
where, for $\delta > 0$, $f_{\delta,K}$ is defined as in Definition \ref{df:PullIn}.
\end{lemma}

\begin{proof} Fix $\alpha > 0$ (to be chosen later). Since $K \subset \R^{d}$ is compact, there is $\eta > 0$ such that, for $x_{1},x_{2} \in K$,
\begin{align}\label{PullInto1}
    \vert x_{2} - x_{1} \vert \leq \eta &\Rightarrow \vert \psi(x_{2},\xi)-\psi(x_{1},\xi) \vert \leq \alpha \cdot \Vert \xi \Vert \hspace{1em} \text{for $\xi \in \subext{k}\R^{d}$}.
\end{align}
After possibly adjusting $\eta$, we can also ensure that $(1+\eta)^{k} -1 \leq \alpha$. Take $\delta > 0$ with $\delta \leq \eta$. We shall now use Lemma \ref{lm:psiMassPushFwd} to estimate $\Mbf_{\psi}((f_{\delta,K})_{*}T)$; in what follows, we abbreviate $f = f_{\delta,K}$:
\begin{align*}
    \Mbf_{\psi}(f_{*}T) &\leq \int \psi(f(x),D^{T}f(x)\vec{T}(x)) \dd \Vert T \Vert(x) \\
    &\leq \int \psi(f(x),\vec{T}(x)) + \psi(f(x),D^{T}f(x)\vec{T}(x)-\vec{T}(x)) \dd \Vert T \Vert(x)\\
    &\leq \int \psi(f(x),\vec{T}(x)) + C_{\psi,K}\cdot \Vert D^{T}f(x)\vec{T}(x)-\vec{T}(x) \Vert \dd \Vert T \Vert(x).
\end{align*}
Here, $C_{\psi,K} > 0$ is defined as in \eqref{psiVsMass}. Since $\Lip(f-\mathrm{id}_{K}) \leq \delta$, we can use Lemma \ref{lm:MultiEst} to estimate
\begin{align}\label{PullInto2}
    \Mbf_{\psi}(f_{*}T) &\leq \int \psi(f(x),\vec{T}(x)) \dd \Vert T \Vert + ((1+\delta)^{k}-1) C_{\psi,K} \cdot \Mbf(T).
\end{align}
Since $\Vert f - \mathrm{id}_{K} \Vert_{\Lrm^{\infty}} \leq \delta$, we can use \eqref{PullInto1} to estimate
\begin{equation}\label{PullInto3}
\begin{aligned}
    \int &\psi(f(x),\vec{T}(x)) \dd \Vert T \Vert(x) \\
    &\leq \Mbf_{\psi}(T) + \int \vert \psi(f(x),\vec{T}(x)) - \psi(x,\vec{T}(x)) \vert \dd \Vert T \Vert(x) \leq \Mbf_{\psi}(T) + \alpha \cdot \Mbf(T).
\end{aligned}
\end{equation}
Combining \eqref{PullInto2} and \eqref{PullInto3} and using the fact that $(1+\delta)^{k}-1 \leq \alpha$ yields
\begin{align*}
    \Mbf_{\psi}(f_{*}T) &\leq \Mbf_{\psi}(T) + \alpha(C_{\psi,K} + 1) \cdot \Mbf(T) \\
    &\leq (1+\alpha \cdot C_{\psi,K}(C_{\psi,K}+1)) \cdot \Mbf_{\psi}(T).
\end{align*}
The result follows after choosing $\alpha$ so that
\begin{align*}
    \alpha &\leq \frac{\varepsilon}{C_{\psi,K}(C_{\psi,K}+1)}. \qedhere
\end{align*}
\end{proof}

By combining Lemma \ref{lm:PullInto} with the strict-type approximation results for normal and integral currents obtained in Section \ref{s:ClassicalPolyStrict} and Section \ref{s:Q}, one can obtain analogous results adapted to the setting in which the supporting compact set has the pull-in property. For example, we have the following result:

\begin{proposition}\label{pp:*BdrylessQStrict} Let $\psi: \R^{d} \times \subext{k}\R^{d} \to [0,\infty)$ be a flat integrand. Take a compact set $K \subset \R^{d}$ with the pull-in property and let $\varepsilon > 0$. Take $T \in \Nrm_{k}(K)$ and suppose that $\partial T = 0$. Then there is $Q \in \Qrm_{k}(K)$ with
\begin{align*}
    \partial Q &= 0, & \Fbf_{K}(Q-T) &\leq \varepsilon, & \Mbf_{\psi}(Q) &\leq \Mbf_{\psi}(T) + \varepsilon.
\end{align*}
\end{proposition}

\begin{proof} If $T = 0$, then we can take $Q = 0$ and there is nothing to show. Otherwise, fix $\alpha \in (0,1)$ (to be chosen later) and define $f_{\alpha,K}: K \to \Int K$ as in Definition \ref{df:PullIn}. In what follows, we abbreviate $f = f_{\alpha,K}$. Since $\partial f_{*}T = 0$ and $\supp f_{*}T \subset \Int K$, we can appeal to Proposition \ref{pp:QBdryQStrict} to find $Q \in \Qrm_{k}(\R^{d})$ with $\supp Q \subset K$ and
\begin{align*}
    \partial Q &= 0, & \Fbf_{K}(Q-f_{*}T) &\leq \frac{\varepsilon}{2}, & \Mbf_{\psi}(Q) \leq \Mbf_{\psi}(f_{*}T) + \frac{\varepsilon}{2}.
\end{align*}
Now let $H$ be the affine homotopy from $\mathrm{id}_{K}$ to $f$. Then $H([0,1] \times K) \subset K$ and so
\begin{align*}
    \Fbf_{K}(f_{*}T-T) &\leq \Mbf(H_{*}([\![0,1]\!] \times T)).
\end{align*}
Since $\Lip(f-\mathrm{id}_{K}) \leq \alpha < 1$ and $\Vert f-\mathrm{id}_{K} \Vert_{\Lrm^{\infty}} \leq \alpha$, we have
\begin{align*}
    \Mbf(H_{*}([\![0,1]\!] \times T)) \leq \alpha 2^{k}\cdot \Mbf(T).
\end{align*}
We now use Lemma \ref{lm:PullInto} to find $\eta > 0$ such that, for all $\delta > 0$ with $\delta \leq \eta$, we have
\begin{align*}
    \Mbf_{\psi}((f_{\delta,K})_{*}T) \leq \bigg(1+\frac{\varepsilon}{2 \cdot \Mbf_{\psi}(T)}\bigg) \cdot \Mbf_{\psi}(T). 
\end{align*}
where, for $\delta > 0$, $f_{\delta,K}$ is defined as in Definition \ref{df:PullIn}. The result follows after choosing $\alpha$ so that
\begin{align*}
    \alpha &\leq \min\bigg\{\frac{\varepsilon}{2^{k+1}\cdot\Mbf(T)}, \eta \bigg\}. \qedhere
\end{align*}
\end{proof}

\subsection{Space-Time Trajectories} We start by working as in the proof of Lemma \ref{lm:PullInto} to obtain the following estimate:

\begin{lemma}\label{lm:STPullInto} Let $\phi: \R^{d} \times \subext{1+k}\R^{d} \to [0,\infty)$ be a positive convex integrand. Take a compact $K \subset \R^{d}$ with the pull-in property and let $\varepsilon > 0$. Then there is $\eta > 0$ with the following property: For $S \in \Nrm_{1+k}([0,1] \times K)$ and $\delta > 0$ with $\delta \leq \eta$, we have
\begin{align}\label{cc:STPullInto}
    \Var_{\phi}((\overline{f_{\delta,K}})_{*}S) &\leq (1+\varepsilon)\Var_{\phi}(S).
\end{align}
Here, for $\delta > 0$, we set $\overline{f_{\delta,K}} = \mathrm{id}_{[0,1]} \times f_{\delta,K}$, where $f_{\delta,K}$ is defined as in Definition \ref{df:PullIn}.
\end{lemma}

\begin{proof} Fix $\alpha > 0$ (to be chosen later). Since $K \subset \R^{d}$ is compact, there is $\eta > 0$ such that, for $x_{1},x_{2} \in K$, 
\begin{align}\label{STPullInto1}
    \vert x_{2} - x_{1} \vert \leq \eta &\Rightarrow \vert \phi(x_{2},\xi)-\phi(x_{1},\xi) \vert \leq \alpha \cdot \Vert \xi \Vert \hspace{1em} \text{for $\xi \in \subext{1+k}\R^{d}$}.
\end{align}
After possibly adjusting $\eta$, we can also ensure that $(1+\eta)^{k+1} -1 \leq \alpha$. Take $\delta > 0$ with $\delta \leq \eta$. We shall now use Lemma \ref{lm:PushFwdDiss} to estimate $\Var_{\phi}
((\overline{f_{\delta,K}})_{*}S)$; in what follows, we abbreviate $f = f_{\delta,K}$ and $\overline{f} = \overline{f_{\delta,K}}$:
\begin{align*}
    \Var_{\phi}&(\overline{f}_{*}S) \leq \int \phi(f(x),\pbf(D^{S}\overline{f}(t,x)\vec{S}(t,x))) \dd \Vert S \Vert(t,x) \\
    &\leq \int \phi(f(x),\pbf\vec{S}(t,x)) + \phi(f(x),\pbf(D^{S}\overline{f}(t,x)\vec{S}(t,x)-\vec{S}(t,x)) \dd \Vert S \Vert(t,x) \\
    &\leq \int \phi(f(x),\pbf\vec{S}(t,x)) + C_{\phi,K}\cdot \Vert \pbf(D^{S}\overline{f}(t,x)\vec{S}(t,x)-\vec{S}(t,x)) \Vert \dd \Vert S \Vert(t,x) 
\end{align*}
where $C_{\phi,K} > 0$ is defined as in \eqref{psiVsMass}. Take $(t,x) \in \R^{1+d}$ and suppose that $D^{S}\overline{f}(t,x)$ exists. Set $V := \pbf(\mathrm{span}\hspace{0.2em}\vec{S}(t,x))$ and notice that
\begin{align*}
    \pbf(D^{S}\overline{f}(t,x)v) = \lim_{h \to 0} \frac{f(x+h\pbf v)-f(x)}{h} \hspace{1em} \text{for $v \in \mathrm{span}\hspace{0.2em}\vec{S}(t,x)$.}
\end{align*}
We can therefore define a linear map $g(t,x): V \to \R^{d}$ via
\begin{align*}
    g(t,x)(\pbf v) := \pbf(D^{S}\overline{f}(t,x)v) \hspace{1em} \text{for $v \in \mathrm{span}\hspace{0.2em}\vec{S}(t,x)$}.
\end{align*}
It follows that
\begin{align}\label{STPullInto1.5}
    \pbf(D^{S}\overline{f}(t,x)\vec{S}(t,x))=g(t,x)(\pbf\vec{S}(t,x)).
\end{align}
One can check that
\begin{align*}
    \Lip(g(t,x) - \mathrm{id}_{V}) \leq \Lip(f-\mathrm{id}_{\R^{d}}) \leq \delta 
\end{align*}
and we can use Lemma \ref{lm:MultiEst} in conjunction with \eqref{STPullInto1.5} to deduce that
\begin{align*}
    \Vert \pbf(D^{S}\overline{f}(t,x)\vec{S}(t,x)-\vec{S}(t,x)) \Vert \leq ((1+\delta)^{k+1}-1) \Vert \pbf\vec{S}(t,x) \Vert.
\end{align*}
We therefore have
\begin{align}\label{STPullInto2}
    \Var_{\phi}(\overline{f}_{*}S) \leq \int \phi(f(x),\pbf\vec{S}(t,x)) \dd \Vert S \Vert(t,x) + ((1+\delta)^{k+1}-1)C_{\phi,K}\cdot \Var(S)
\end{align}
Since $\Vert f - \mathrm{id}_{K} \Vert_{\Lrm^{\infty}} \leq \delta$, we can use \eqref{STPullInto1} to estimate
\begin{equation}\label{STPullInto3}
\begin{aligned}
    \int \phi&(f(x),\pbf\vec{S}(t,x)) \dd \Vert S \Vert(t,x) \\
    &\leq \Var_{\phi}(S) + \int \vert \phi(f(x),\pbf\vec{S}(t,x)) -\phi(x,\pbf\vec{S}(t,x)) \vert \dd \Vert S \Vert(t,x) \\
    &\leq \Var_{\phi}(S) + \alpha \cdot \Var(S) \vphantom{\int \vert \phi(f(x),\pbf\vec{S}(t,x))}.
\end{aligned}
\end{equation}
Combining \eqref{STPullInto2} and \eqref{STPullInto3} and using the fact that $(1+\delta)^{k+1} - 1 \leq \alpha$ yields
\begin{align*}
    \Var_{\phi}(\overline{f}_{*}S) &\leq \Var_{\phi}(S) + \alpha(C_{\phi,K}+1)\cdot \Var(S) \\
    &\leq (1+\alpha \cdot C_{\phi,K}(C_{\phi,K}+1))\cdot\Var_{\phi}(S).
\end{align*}
The result follows after choosing $\alpha > 0$ so that
\begin{align*}
    \alpha &\leq \frac{\varepsilon}{C_{\phi,K}(C_{\phi,K}+1)}. \qedhere
\end{align*}
\end{proof}

The following result will allow us to adapt the strict-type approximation results for space-time trajectories established in Section \ref{s:Q} and Section \ref{s:Lipschitz} to prove analogous results in the setting where the supporting set has the pull-in property:

\begin{lemma}\label{lm:STPullIntoEsts} Take positive convex integrands 
\begin{align*}
    \psi &: \R^{d} \times \subext{k}\R^{d} \to [0,\infty), & \phi &: \R^{d} \times \subext{1+k}\R^{d} \to [0,\infty). 
\end{align*}
Take a compact set $K \subset \R^{d}$ with the pull-in property and let $\varepsilon > 0$. Let $S \in \Nrm_{1+k}([0,1] \times K)$ and suppose that $\partial S \mvert (0,1) \times \R^{d} = 0$. Then there is $\eta > 0$ such that, for $\delta > 0$ with $\delta \leq \eta$, we have
\begin{align*}
    \supp (\overline{f_{\delta,K}})_{*}S &\subset [0,1] \times \Int K, & \partial(\overline{f_{\delta,K}})_{*}S \mvert (0,1) \times \R^{d} &= 0
\end{align*}
and the following conclusions hold:
\begin{enumerate}
    \item We have
    \begin{align*}
        \FN_{[0,1] \times K}((\overline{f_{\delta,K}})_{*}S-S) &\leq \varepsilon, & \Var_{\phi}((\overline{f_{\delta,K}})_{*}S) &\leq \Var_{\phi}(S) + \varepsilon,
    \end{align*}
    and
    \begin{align*}
        \int \Mbf_{\psi}(((\overline{f_{\delta,K}})_{*}S)(t)) \dd t \leq \int \Mbf_{\psi}(S(t)) \dd t + \varepsilon.
    \end{align*}
    \item We have
    \begin{align*}
        \FN^{\circ}_{K}(\partial^{\pm}((\overline{f_{\delta,K}})_{*}S)-\partial^{\pm}S) &\leq \varepsilon, & \Mbf_{\psi}(\partial^{\pm}((\overline{f_{\delta,K}})_{*}S)) &\leq \Mbf_{\psi}(\partial^{\pm}S) + \varepsilon.
    \end{align*}
    \item We have
    \begin{align*}
    \esssup_{t \in [0,1]} \Mbf_{\psi}(((\overline{f_{\delta,K}})_{*}S)(t)) \leq \esssup_{t \in [0,1]} \Mbf_{\psi}(S(t)) + \varepsilon.
    \end{align*}
    \item If $\partial^{\pm}S \in \Irm_{k}(\R^{d})$, then $\partial^{\pm}((\overline{f_{\delta,K}})_{*}S) \in \Irm_{k}(\R^{d})$ and 
    \begin{align*}
        \FI_{K}^{\circ}(\partial^{\pm}((\overline{f_{\delta,K}})_{*}S)-\partial^{\pm}S) &\leq \varepsilon.
    \end{align*}
    \item If $S \in \Irm_{1+k}([0,1] \times \R^{d})$, then $((\overline{f_{\delta,K}})_{*}S) \in \Irm_{1+k}([0,1] \times \R^{d})$ and
    \begin{align*}
        \FI_{[0,1] \times K}((\overline{f_{\delta,K}})_{*}S-S) &\leq \varepsilon.
    \end{align*}
\end{enumerate} 
Here, for $\delta > 0$, we set $\overline{f_{\delta,K}} := \mathrm{id}_{[0,1]} \times f_{\delta,K}$, where $f_{\delta,K}$ is defined as in Definition \ref{df:PullIn}.
\end{lemma}

\begin{proof} Fix $\alpha \in (0,1)$ (to be chosen later). Let $\eta > 0$ be the minimum of the two constants from Lemma \ref{lm:PullInto} and Lemma \ref{lm:STPullInto}. By possibly adjusting $\eta$, we can ensure that $\eta \leq \alpha$. Take $\delta > 0$ with $\delta \leq \eta$. In what follows, we abbreviate $f = f_{\delta,K}$ and $\overline{f} = \overline{f_{\delta,K}}$. 

Since $\overline{f}([0,1] \times K) \subset [0,1] \times \Int K$, we have $\supp \overline{f}_{*}S \subset [0,1] \times \Int K$. One can check that
\begin{align*}
    \partial(\overline{f}_{*}S) = \delta_{1} \times f_{*}\partial^{+}S - \delta_{0} \times f_{*}\partial^{-}S.
\end{align*}
It follows that $\partial(\overline{f}_{*}S) \mvert (0,1) \times \R^{d} = 0$ and that $\partial^{\pm}(\overline{f}_{*}S) = f_{*}\partial^{\pm}S$.

Let $H$ and $\overline{H}$ be the affine homotopies from $\mathrm{id}_{K}$ to $f$ and from $\mathrm{id}_{[0,1] \times K}$ to $\overline{f}$ respectively. Since $H([0,1] \times K) \subset K$, we have $\overline{H}([0,1] \times ([0,1] \times K)) \subset [0,1] \times K$; combining this observation with the homotopy formula \eqref{Homotopy} yields
\begin{align*}
    \FN_{[0,1] \times K}(\overline{f}_{*}S-S) &\leq \Mbf(\overline{H}_{*}([\![0,1]\!] \times S)) + \Mbf(\overline{H}_{*}([\![0,1]\!] \times \partial S)).
\end{align*}
Since $\Lip(\overline{f}-\mathrm{id}_{[0,1] \times K}) \leq \alpha < 1$ and $\Vert \overline{f}-\mathrm{id}_{[0,1] \times K}\Vert_{\Lrm^{\infty}} \leq \alpha$, we can estimate
\begin{align}\label{STPullIntoEsts1}
    \Mbf(\overline{H}_{*}([\![0,1]\!] \times S)) + \Mbf(\overline{H}_{*}([\![0,1]\!] \times \partial S)) \leq 2^{k+1}\alpha \cdot \Nbf(S)
\end{align}
and so 
\begin{align}\label{STPullIntoEsts2}
    \FN_{[0,1] \times K}(\overline{f}_{*}S-S) \leq 2^{k+1}\alpha \cdot \Nbf(S).
\end{align}
It follows from Lemma \ref{lm:STPullInto} that
\begin{align}\label{STPullIntoEsts3}
    \Var_{\phi}(\overline{f}_{*}S) &\leq (1+\alpha) \Var_{\phi}(S)
\end{align}
and, since $(\overline{f}_{*}S)(t) = f_{*}S(t)$ for $\mathscr{L}^{1}$-a.e.\ $t \in \R$ (see \cite[4.3.2 (7)]{Federer69book} or \cite[Proposition 3.9]{BonicattoRindlerTurnbull26}), we can use Lemma \ref{lm:PullInto} to deduce that
\begin{align}\label{STPullIntoEsts4}
    \int \Mbf_{\psi}((\overline{f}_{*}S)(t)) \dd t \leq \int \Mbf_{\psi}(f_{*}S(t)) \dd t \leq (1+\alpha)\int \Mbf_{\psi}(S(t)) \dd t.
\end{align}
Conclusion (1) now follows from \eqref{STPullIntoEsts2}, \eqref{STPullIntoEsts3}, and \eqref{STPullIntoEsts4} provided that we ensure that
\begin{align*}
    \alpha \cdot \max\bigg\{\Var_{\phi}(S),2^{k+1}\cdot \Nbf(S),\int \Mbf_{\psi}(S(t)) \dd t \bigg\} \leq \varepsilon.
\end{align*}

Since $\partial^{\pm}(\overline{f}_{*}S) = f_{*}\partial^{\pm}S$, we can use the homotopy formula \eqref{Homotopy} to deduce that
\begin{align*}
    \partial^{\pm}(\overline{f}_{*}S) - \partial^{\pm}S = \partial H_{*}([\![0,1]\!] \times \partial^{\pm}S),
\end{align*}
and so we can work as in \eqref{STPullIntoEsts1} to estimate
\begin{align}\label{STPullIntoEsts5}
    \Fbf^{\circ}_{K}(\partial^{\pm}(\overline{f}_{*}S)-\partial^{\pm}S) \leq \Mbf(H_{*}([\![0,1]\!] \times \partial^{\pm}S)) \leq 2^{k}\alpha \cdot \Mbf(\partial^{\pm}S).
\end{align}
On the other hand, we can use Lemma \ref{lm:PullInto} to deduce that
\begin{align}\label{STPullIntoEsts6}
    \Mbf_{\psi}(\partial^{\pm}(\overline{f}_{*}S)) &\leq (1+\alpha)\Mbf_{\psi}(\partial^{\pm}S).
\end{align}
Conclusion (2) now follows from \eqref{STPullIntoEsts5} and \eqref{STPullIntoEsts6} provided that we ensure that
\begin{align*}
    \alpha \cdot \max\{2^{k}\cdot\Mbf(\partial^{+}S),2^{k}\cdot\Mbf(\partial^{-}S),\Mbf_{\psi}(\partial^{+}S),\Mbf_{\psi}(\partial^{-}S)\} \leq \varepsilon.
\end{align*}
To obtain conclusion (3), first notice that there is nothing to show in the case that 
\begin{align*}
    \esssup_{t \in [0,1]} \Mbf_{\psi}(S(t)) = \infty.
\end{align*}
Otherwise, we use Lemma \ref{lm:PullInto} (as in \eqref{STPullIntoEsts4}) to estimate
\begin{align*}
    \esssup_{t \in [0,1]} \Mbf_{\psi}((\overline{f}_{*}S)(t)) \leq (1+\alpha) \esssup_{t \in [0,1]} \Mbf_{\psi}(S(t)).
\end{align*}
and the conclusion follows provided that we ensure that
\begin{align*}
    \alpha \cdot \esssup_{t \in [0,1]} \Mbf_{\psi}(S(t)) \leq \varepsilon.
\end{align*}
We obtain conclusion (4) after noting that $H_{*}([\![0,1]\!] \times \partial^{\pm}S) \in \Irm_{k+1}(K)$ in this case. Similarly, we obtain conclusion (5) after noting that $\overline{H}_{*}([\![0,1]\!] \times S) \in \Irm_{k+2}([0,1] \times K)$ and $\overline{H}_{*}([\![0,1]\!] \times \partial S) \in \Irm_{k+1}([0,1] \times K)$ in this case.
\end{proof}

Finally, we give three examples of adaptations of results from Section \ref{s:Q} and Section \ref{s:Lipschitz} to the setting where the supporting set has the pull-in property.

\begin{proposition}\label{pp:*STQ} Take flat integrands
\begin{align*}
    \psi &: \R^{d} \times \subext{k}\R^{d} \to [0,\infty), & \phi &: \R^{d} \times \subext{1+k}\R^{d} \to [0,\infty). 
\end{align*}
Take a compact $K \subset \R^{d}$ with the pull-in property and let $\varepsilon > 0$. Let $S \in \Nrm_{1+k}([0,1] \times K)$ and suppose that $\partial S \mvert (0,1) \times \R^{d} = 0$. Then there is $Q \in \QPrm_{1+k}([0,1] \times K)$ such that $\partial Q \mvert (0,1) \times \R^{d} = 0$ and the following conclusions hold:
\begin{enumerate}
    \item We have
    \begin{align*}
        \FN_{[0,1] \times K}(Q-S) &\leq \varepsilon
    \end{align*}
    and
    \begin{align*}
        \Var_{\phi}(Q) &\leq \Var_{\phi}(S) + \varepsilon, & \int \Mbf_{\psi}(Q(t)) \dd t \leq \int \Mbf_{\psi}(S(t)) \dd t + \varepsilon.
    \end{align*}
    \item We have
    \begin{align*}
        \FN_{K}(\partial^{\pm}Q - \partial^{\pm}S) &\leq \varepsilon, & \Mbf_{\psi}(\partial^{\pm}Q) &\leq \Mbf_{\psi}(\partial^{\pm}S) + \varepsilon.
    \end{align*}
    \item We have
    \begin{align*}
        \esssup_{t \in [0,1]} \Mbf_{\psi}(Q(t)) &\leq \esssup_{t \in [0,1]}\Mbf_{\psi}(S(t)) + \varepsilon.
    \end{align*}
\end{enumerate}
\end{proposition}

\begin{proof} Let $\eta > 0$ be the constant from Lemma \ref{lm:STPullIntoEsts} (with $\varepsilon/2$ in place of $\varepsilon$), and let $\overline{f_{\eta,K}} := \mathrm{id}_{[0,1]} \times f_{\eta,K}$, where $f_{\eta,K}$ is defined as in Definition \ref{df:PullIn}. It then follows from Lemma \ref{lm:STPullIntoEsts} that $\supp (\overline{f_{\eta,K}})_{*}S \subset [0,1] \times \Int K$, $\partial (\overline{f_{\eta,K}})_{*}S \mvert (0,1) \times \R^{d} = 0$, and that conclusions (1)-(3) hold with $Q$ replaced by $(\overline{f_{\eta,K}})_{*}S$ and $\varepsilon$ replaced by $\varepsilon/2$. The result follows by applying Proposition \ref{pp:STQ} to $(\overline{f_{\eta,K}})_{*}S$ (with $\varepsilon/2$ in place of $\varepsilon$).
\end{proof}

\begin{proposition}\label{pp:*STOneSidedQ} Take flat integrands
\begin{align*}
    \psi &: \R^{d} \times \subext{k}\R^{d} \to [0,\infty), & \phi &: \R^{d} \times \subext{1+k}\R^{d} \to [0,\infty). 
\end{align*}
Take a compact $K \subset \R^{d}$ with the pull-in property and let  $\varepsilon > 0$. Let $S \in \Nrm_{1+k}([0,1] \times K)$. Suppose that $\partial S \mvert (0,1) \times \R^{d} = 0$ and that $\partial^{-}S \in \Qrm_{k}(\R^{d})$. Then there is $Q \in \Qrm_{1+k}([0,1] \times K)$ such that $\partial Q \mvert (0,1) \times \R^{d} = 0$, conclusions (1)-(3) of Proposition \ref{pp:*STQ} hold, and $\partial^{-}Q = \partial^{-}S$.
\end{proposition}

\begin{proof} Fix $\alpha > 0$ (to be chosen later). It is not restrictive to assume that $\partial^{-}S \in \Irm_{k}(\R^{d})$. Let $\eta$ be the constant from Lemma \ref{lm:STPullIntoEsts} (with $\alpha$ in place of $\varepsilon$), and set $\overline{f_{\eta,K}} = \mathrm{id}_{[0,1]} \times f_{\eta,K}$, where $f_{\eta,K}$ is defined as in Definition \ref{df:PullIn}. It follows from Lemma \ref{lm:STPullIntoEsts} that $\supp (\overline{f_{\eta,K}})_{*}S \subset [0,1] \times \Int K$, $\partial (\overline{f_{\eta,K}})_{*}S \mvert (0,1) \times \R^{d} = 0$, and that conclusions (1)-(3) of Proposition \ref{pp:*STQ} hold with $Q$ replaced by $(\overline{f_{\eta,K}})_{*}S$ and $\varepsilon$ replaced by $\alpha$. Furthermore, since $\partial^{-}S \in \Irm_{k}(\R^{d})$, we have conclusion (4) of Lemma \ref{lm:STPullIntoEsts}; in particular, we have
\begin{align}\label{*STOneSidedQ1}
    \FI_{K}^{\circ}(\partial^{-}(\overline{f_{\eta,K}})_{*}S-\partial^{-}S) &\leq \alpha.
\end{align}
We now appeal to Theorem \ref{tm:STOneSidedQ} to find a $\tilde{Q} \in \Qrm_{1+k}([0,1] \times K)$ such that $\partial \tilde{Q} \mvert (0,1) \times \R^{d} = 0$, $\partial^{-}\tilde{Q}=\partial^{-}((\overline{f_{\eta,K}})_{*}S)$, and conclusions (1)-(3) of Proposition \ref{pp:*STQ} hold with $S$ replaced by $(\overline{f_{\eta,K}})_{*}S$, $Q$ replaced by $\tilde{Q}$, and $\varepsilon$ replaced by $\alpha$. It follows from \eqref{*STOneSidedQ1} that
\begin{align*}
    \FI_{K}^{\circ}(\partial^{-}\tilde{Q}-\partial^{-}S) &\leq \alpha
\end{align*}
and so we can find $R \in \Irm_{k+1}(K)$ with $\partial R = \partial^{-}\tilde{Q}-\partial^{-}S$ and $\Mbf(R) \leq \alpha$. Now set
\begin{align*}
    Q = \tilde{Q} + \delta_{0} \times R.
\end{align*}
Then $\supp Q \subset [0,1] \times K$, $\partial Q \mvert (0,1) \times \R^{d} = 0$, $\partial^{-}Q = \partial^{-}S$, and $\partial^{+}Q = \partial^{+}\tilde{Q}$. Notice that 
\begin{align*}
    \FN_{[0,1] \times K}(Q-S) &\leq 3\alpha
\end{align*}
and we can use Lemma \ref{lm:DissJumps} to estimate
\begin{align*}
    \Var_{\phi}(Q) &\leq \Var_{\phi}(S) + (C_{\phi,K}+2)\alpha
\end{align*}
where $C_{\phi,K}$ is defined as in \eqref{psiVsMass}. Since $\Vert \delta_{0} \times R \Vert_{\Lrm^{\infty}} = 0$, we have
\begin{align*}
    \int \Mbf_{\psi}(Q(t)) \dd t \leq \int \Mbf_{\psi}(S(t)) \dd t + 2\alpha
\end{align*}
and
\begin{align*}
    \esssup_{t \in [0,1]} \Mbf_{\psi}(Q(t)) &\leq \esssup_{t \in [0,1]}\Mbf_{\psi}(S(t)) + 2\alpha.
\end{align*}
Finally, since $\partial^{+}Q = \partial^{+}\tilde{Q}$, we have
\begin{align*}
    \FN_{K}(\partial^{\pm}Q-\partial^{\pm}S) &\leq 2\alpha
\end{align*}
and
\begin{align*}
    \Mbf_{\psi}(\partial^{\pm}Q) &\leq \Mbf_{\psi}(\partial^{\pm}S) + 2\alpha.
\end{align*}
The result follows after choosing $\alpha$ so that
\begin{align*}
    \alpha &\leq \varepsilon \cdot \min\bigg\{\frac{1}{3}, \frac{1}{C_{\phi,K}+2}\bigg\}. \qedhere
\end{align*}
\end{proof}

\begin{proposition} Take flat integrands
\begin{align*}
    \psi &: \R^{d} \times \subext{k}\R^{d} \to [0,\infty), & \phi &: \R^{d} \times \subext{1+k}\R^{d} \to [0,\infty). 
\end{align*}
Take a compact set $K \subset \R^{d}$ with the pull-in property and let $\varepsilon > 0$. Let $S \in \Nrm_{1+k}([0,1] \times K)$ and suppose that $\partial S \mvert (0,1) \times \R^{d} = 0$. Then there is $V \in \NLip_{1+k}([0,1] \times K)$ such that $\partial V = \partial S$, 
\begin{align*}
        \FN_{[0,1] \times K}(V-S) &\leq \varepsilon,
\end{align*}
and
\begin{align*}
    \Var_{\phi}(V) &\leq \Var_{\phi}(S) + \varepsilon, & \int \Mbf_{\psi}(V(t)) \dd t \leq \int \Mbf_{\psi}(S(t)) \dd t + \varepsilon.
\end{align*}
If $S \in \Irm_{1+k}([0,1] \times \R^{d})$, then we can additionally ensure that $V \in \ILip_{1+k}([0,1] \times \R^{d})$.
\end{proposition}

\begin{proof} Fix $\alpha > 0$ (to be chosen later). Let $\eta$ be the constant from Lemma \ref{lm:STPullIntoEsts} (with $\alpha$ in place of $\varepsilon$), and set $\overline{f_{\eta,K}} := \mathrm{id}_{[0,1]} \times f_{\eta,K}$, where $f_{\eta,K}$ is defined as in Definition \ref{df:PullIn}. It follows from Lemma \ref{lm:STPullIntoEsts} that $\supp (\overline{f_{\eta,K}})_{*}S \subset [0,1] \times \Int K$, $\partial (\overline{f_{\eta,K}})_{*}S \mvert (0,1) \times \R^{d} = 0$, and that conclusions (1)-(3) of Lemma \ref{lm:STPullIntoEsts} hold with $\varepsilon$ replaced with $\alpha$. We then appeal to Theorem \ref{tm:STLipschitz} to find $\tilde{V} \in \NLip_{1+k}([0,1] \times \R^{d})$ with $\partial \tilde{V} = \partial (\overline{f_{\eta,K}})_{*}S$, $\supp \tilde{V} \subset [0,1] \times K$,
\begin{align*}
        \FN_{[0,1] \times K}(\tilde{V}-(\overline{f_{\eta,K}})_{*}S) &\leq \alpha, & \Var_{\phi}(\tilde{V}) &\leq \Var_{\phi}((\overline{f_{\eta,K}})_{*}S) + \alpha,
\end{align*}
and
\begin{align*}
    \int \Mbf_{\psi}(\tilde{V}(t)) \dd t \leq \int \Mbf_{\psi}(((\overline{f_{\eta,K}})_{*}S)(t)) \dd t + \alpha.
\end{align*}
Notice that we have the following estimate from conclusion (2) of Lemma \ref{lm:STPullIntoEsts}:
\begin{align*}
    \FN^{\circ}_{K}(\partial^{\pm}((\overline{f_{\eta,K}})_{*}S)-\partial^{\pm}S) &\leq \alpha,
\end{align*}
We can therefore appeal to \cite[Theorem 5.3]{BonicattoRindlerTurnbull26} (which holds not only for closures of $\Crm^{0}$ domains but for all compact sets with the pull-in property) to construct $S^{\pm} \in \NLip_{1+k}([0,1] \times K)$ with $\Var(S^{\pm}) \leq 2\alpha$ and
\begin{align*}
    \partial S^{+} &= \delta_{1} \times \partial^{+}S - \delta_{0} \times \partial^{+}(\overline{f_{\eta,K}})_{*}S, \\
    \partial S^{-} &= \delta_{1} \times \partial^{-}(\overline{f_{\eta,K}})_{*}S - \delta_{0} \times \partial^{-}S
\end{align*}
We now appeal to Lemma \ref{lm:Gluing} (with $2\alpha$ in place of $\varepsilon$) to obtain $V \in \NLip_{1+k}([0,1] \times K)$ with $\partial^{\pm}V=\partial^{\pm}S^{\pm} = \partial^{\pm}S$,
\begin{align*}
    \FN_{[0,1] \times K}(V-\tilde{V}) &\leq 6\alpha,
\end{align*}
and 
\begin{align*}
    \Var_{\phi}(V) &\leq \Var_{\phi}(\tilde{V}) + 4C_{\phi,K}\cdot\alpha, & \int \Mbf_{\psi}(V(t)) \dd t \leq \int \Mbf_{\psi}(\tilde{V}(t)) \dd t + 2\alpha.
\end{align*}
Now choose $\alpha$ so that
\begin{align*}
    \alpha &\leq \varepsilon \cdot \min\bigg\{\frac{1}{8},\frac{1}{2(1+2C_{\phi,K})}\bigg\}.
\end{align*}
Combining the estimates above yields the result. If $S \in \Irm_{1+k}([0,1] \times \R^{d})$, then we follow the argument above, noting that conclusions (4) and (5) of Lemma \ref{lm:STPullIntoEsts}, the integral version of \cite[Theorem 5.3]{BonicattoRindlerTurnbull26} (see also \cite[Theorem 5.4]{Rindler23}), and the additional conclusion of Lemma \ref{lm:Gluing} hold in this case.
\end{proof}

\appendix

\section{Reshetnyak's Theorem in Normed Spaces}\label{s:Reshetnyak}

In what follows, we fix a finite-dimensional normed space $(X,\vert \cdot \vert)$ and take a compact $K \subset \R^{d}$. We shall use $\Mcal(K;X)$ to denote the space of finite $X$-valued Borel measures on $K$. Set
\begin{align*}
    \Bbb_{X} &:= \{v \in X: \vert v \vert \leq 1\}, & \Sbb_{X} := \{v \in X: \vert v \vert = 1\}.
\end{align*}
Recall that a convex subset $F \subset \Bbb_{X}$ is a \textbf{face} of $\Bbb_{X}$ if, for $x,y \in \Bbb_{X}$ and $\lambda \in (0,1)$, we have
\begin{align*}
    (1-\lambda)x+\lambda y \in F &\Rightarrow x,y \in F.
\end{align*}
We say that a face $F \subset \Bbb_{X}$ is \textbf{proper} if $F \neq \Bbb_{X}$. One can show that, if $F$ is a proper face of $\Bbb_{X}$, then $F \subset \Sbb_{X}$.

The following result is a slight generalisation of Reshetnyak's continuity theorem: 

\begin{proposition}\label{pp:ReshetnyakCty} Take a sequence $(\mu_{j})_{j} \subset \Mcal(K;X)$, let $\mu \in \Mcal(K;X)$, and suppose that 
\begin{align}\label{MeasureStrictConv}
    \mu_{j} &\wkto \mu, & \vert \mu_{j} \vert(K) &\to \vert \mu \vert(K).
\end{align}
Take a continuous $f: K \times \Sbb_{X} \to \R$ and suppose that $f(x,\cdot)$ is affine on the proper faces of $\Bbb_{X}$ for $\vert \mu \vert$-a.e.\ $x \in K$. Then
\begin{align}\label{Reshetnyak}
    \int f\bigg(x,\frac{\mathrm{d}\mu_{j}}{\mathrm{d}\vert \mu_{j} \vert}(x)\bigg) \dd\vert \mu_{j} \vert(x) \to \int  f\bigg(x,\frac{\mathrm{d}\mu}{\mathrm{d}\vert \mu \vert}(x)\bigg) \dd\vert \mu \vert(x).
\end{align}
\end{proposition} 
\noindent Indeed, if $(X,\vert \cdot \vert)$ is strictly convex, then the second condition on the integrand $f$ in Proposition \ref{pp:ReshetnyakCty} is vacuous and we recover the classical Reshetnyak continuity theorem (see, for example, \cite[Theorem 2.39]{AmbrosioFuscoPallara00book}). Proposition \ref{pp:ReshetnyakCty} follows from \cite[Theorem 0.2]{Delladio91}, but the arguments there can be simplified in our case; we therefore give a complete proof.

\begin{proof}[Proof of Proposition \ref{pp:ReshetnyakCty}] Take a (not relabelled) subsequence of $(\mu_{j})_{j}$. For $j \in \N$, let 
\begin{align*}
    \alpha_{j} := \bigg(\mathrm{id}_{K}, \hspace{0.2em} \frac{\mathrm{d}\mu_{j}}{\mathrm{d}\vert \mu_{j} \vert}\bigg)_{\#}\vert \mu_{j}\vert,
\end{align*}
that is, $\alpha_{j}$ is the finite positive Borel measure on $K \times \Sbb_{X}$ given by
\begin{align*}
    \int \phi \dd \alpha_{j} = \int \phi\bigg(x,\frac{\mathrm{d}\mu_{j}}{\mathrm{d}\vert \mu_{j} \vert}(x)\bigg) \dd\vert\mu_{j}\vert(x)
\end{align*}
for all continuous $\phi: K \times \Sbb_{X} \to \R$. Notice that $\alpha_{j}(K \times \Sbb_{X}) = \vert \mu_{j} \vert(K)$ for $j \in \N$ and so $(\alpha_{j}(K \times \Sbb_{X}))_{j}$ is uniformly bounded. Take a (not relabelled) subsequence of $(\alpha_{j})_{j}$ with $\alpha_{j} \wkto \alpha$ for some finite positive Borel measure $\alpha$ on $K \times \Sbb_{X}$. In particular, we have 
\begin{align}\label{Reshetnyak0}
    \int \phi \dd\alpha_{j} \to \int \phi \dd\alpha
\end{align}
for all continuous $\phi: K \times \Sbb_{X} \to \R$. 

Let $\pi: K \times \Sbb_{X} \to K$ where $\pi(x,v)=x$ for $(x,v) \in K \times \Sbb_{X}$. Notice that $\pi_{\#}\alpha_{j} = \vert \mu_{j} \vert$ for $j \in \N$, and we can use \cite[Corollary 13.2]{Rindler26book} to deduce that $\vert \mu_{j} \vert \wkto \vert \mu \vert$; it follows that $\pi_{\#}\alpha = \vert \mu \vert$ and so $\alpha_{j}(K \times \Sbb_{X}) \to \alpha(K \times \Sbb_{X})$.

We now investigate the structure of $\alpha$. To this end, we disintegrate $\alpha$ with respect to $\pi$ (see, for example, \cite[Theorem 5.3.1]{AmbrosioGigliSavare05book}) to obtain a Borel family $\{\lambda_{x}\}$ 
of probability measures on $\Sbb_{X}$ such that 
\begin{align}\label{Reshetnyak1}
    \int \phi \dd \alpha = \int_{K} \int_{\Sbb_{X}} \phi(x,v) \dd\lambda_{x}(v) \dd\pi_{\#}\alpha(x)
\end{align}
for all bounded Borel $\phi: K \times \Sbb_{X} \to \R$. Take $\psi \in \Crm(K)$. For $j \in \N$, we have
\begin{align}\label{Reshetnyak2}
    \int \psi(x) v \dd \alpha_{j}(x,v) = \int \psi(x) \dd \mu_{j}(x).
\end{align}
We then take $j \to \infty$ in \eqref{Reshetnyak2} and use the weak$^{*}$ convergences $\alpha_{j} \wkto \alpha$ and $\mu_{j} \wkto \mu$ to obtain
\begin{align*}
    \int \psi(x)v \dd \alpha(x,v) = \int \psi(x) \dd\mu(x).
\end{align*}
We can use \eqref{Reshetnyak1} to deduce that
\begin{align*}
    \int \psi(x) \dd\mu(x)=\int \psi(x)v \dd \alpha(x,v) = \int_{K} \psi(x) \int_{\Sbb_{X}} v \dd\lambda_{x}(v) \dd\vert \mu \vert(x).
\end{align*}
Since $\psi$ was arbitrary, we have
\begin{align}\label{Reshetnyak3}
    \int_{\Sbb_{X}} v \dd\lambda_{x}(v) = \frac{\mathrm{d}\mu}{\mathrm{d}\vert \mu \vert}(x) \hspace{1em} \text{for $\vert \mu \vert$-a.e.\ $x \in K$}.
\end{align}

Fix an $x \in K$ for which $f(x,\cdot)$ is affine on the proper faces of $\Bbb_{X}$, \eqref{Reshetnyak3} holds, and
\begin{align*}
    \bigg\vert \frac{\mathrm{d}\mu}{\mathrm{d}\vert \mu \vert}(x) \bigg\vert = 1.
\end{align*}
Notice that this is the case for $\vert\mu\vert$-a.e. $x \in K$. Set
\begin{align*}
    w := \int_{\Sbb_{X}} v \dd \lambda_{x}(v).
\end{align*}
Notice that $\vert w \vert = 1$, so there is $\Phi \in X^{*}$ with dual norm $\vert \Phi \vert = 1$ and $\Phi(w) = 1$. Then
\begin{align*}
    1 = \Phi(w) = \int_{\Sbb_{X}} \Phi(v) \dd\lambda_{x}(v).
\end{align*}
Since $\Phi(v) \leq 1$ for all $v \in \Sbb_{X}$, we must have $\Phi(v) = 1$ for $\lambda_{x}$-a.e.\ $v \in \Sbb_{X}$. In particular, we must have 
\begin{align*}
    \supp \lambda_{x} \subset \Pi := \{v \in \Sbb_{X}: \Phi(v) = 1\}.
\end{align*}
Notice that $\Pi = \{v \in \Bbb_{X}: \Phi(v) = 1\}$ and so $\Pi \subset \Bbb_{X}$ is convex. In fact, $\Pi$ is a face of $\Bbb_{X}$. To see this, let $\lambda \in (0,1)$ and take $v_{1},v_{2} \in \Bbb_{X}$. If $(1-\lambda)v_{1}+\lambda v_{2} \in \Pi$, then $1 = (1-\lambda)\Phi(v_{1})+\lambda \Phi(v_{2})$. Since $\Phi(v_{1}),\Phi(v_{2}) \leq 1$, we must have $\Phi(v_{1})=\Phi(v_{2})=1$ and so $v_{1},v_{2} \in \Pi$. Since $0 \not\in \Pi$, $\Pi$ is a proper face of $\Bbb_{X}$ and so $f(x,\cdot)$ is affine on $\Pi$. We can therefore use Jensen's inequality with $f(x,\cdot)$ and $-f(x,\cdot)$ together with \eqref{Reshetnyak3} to yield
\begin{align}\label{Reshetnyak4}
    \int_{\Sbb_{X}} f(x,v) \dd \lambda_{x}(v) = f\bigg(x,\frac{\mathrm{d}\mu}{\mathrm{d}\vert \mu \vert}(x)\bigg).
\end{align}
Integrating \eqref{Reshetnyak4} with respect to $\vert \mu \vert$, we arrive at
\begin{align*}
    \int f \dd\alpha = \int f\bigg(x,\frac{\mathrm{d}\mu}{\mathrm{d}\vert \mu \vert}(x)\bigg) \dd\vert \mu \vert(x).
\end{align*}
We can now obtain \eqref{Reshetnyak} using \eqref{Reshetnyak0} and the fact that the initial choice of subsequence was arbitrary.
\end{proof}

It turns out that Proposition \ref{pp:ReshetnyakCty} is sharp in the following sense:

\begin{proposition}\label{pp:ReshetnyakSharp} Suppose that $K = [0,1]^{d}$. Take a continuous $f: K \times \Sbb_{X} \to \R$ and suppose that, for all sequences $(\mu_{j})_{j} \subset \Mcal(K;X)$ and $\mu \in \Mcal(K;X)$ satisfying \eqref{MeasureStrictConv}, we have \eqref{Reshetnyak}. Then $f(x,\cdot)$ is affine on the proper faces of $\Bbb_{X}$ for all $x \in K$.
\end{proposition}

\begin{proof} For $t \in [0,1]$, define
\begin{align*}
    A_{j}(t) := \bigcup_{i=0}^{j-1} \bigg[\frac{i}{j},\frac{i+t}{j}\bigg] \subset [0,1] \hspace{1em} \text{for $j \in \N$}.
\end{align*}
Let $g: [0,1] \to \R$ be continuous. One can work as in \cite[Example 4.9]{Rindler26book} to show that
\begin{align}\label{ReshetnyakSharp1}
    \lim_{j \to \infty} \int_{A_{j}(t)} g(x) \dd x &= t \cdot \int_{0}^{1} g(x) \dd x, 
\end{align}
and
\begin{align}\label{ReshetnyakSharp2}
    \lim_{j \to \infty} \int_{[0,1] \backslash A_{j}(t)} g(x) \dd x &= (1-t) \cdot \int_{0}^{1} g(x) \dd x.
\end{align}

Fix a proper face $F$ of $\Bbb_{X}$. Take $v_{1},v_{2} \in F$ and let $t \in [0,1]$. Let $\varphi \in \Crm^{\infty}_{c}((0,1)^{d})$ be non-negative and, for $j \in \N$, define $\mu_{j} \in \Mcal([0,1]^{d};X)$ where
\begin{align*}
    \mu_{j} = (\mathbf{1}_{A_{j}(1-t) \times [0,1]^{d-1}} v_{1} +\mathbf{1}_{([0,1]\backslash A_{j}(1-t)) \times [0,1]^{d-1}}v_{2})\varphi \cdot \mathscr{L}^{d} \mvert [0,1]^{d}.
\end{align*}
One can use \eqref{ReshetnyakSharp1} and \eqref{ReshetnyakSharp2} to deduce that $\mu_{j} \wkto \mu$ where $\mu \in \Mcal([0,1]^{d};X)$ is given by
\begin{align*}
    \mu = ((1-t)v_{1}+tv_{2})\varphi \cdot \mathscr{L}^{d} \mvert [0,1]^{d}.
\end{align*}
Since $v_{1},v_{2} \in F$ and $F \subset \Sbb_{X}$, we have $\vert \mu_{j} \vert = \varphi \cdot \mathscr{L}^{d} \mvert [0,1]^{d}$. Since $F$ is convex, we have $(1-t)v_{1}+tv_{2} \in F$ and so $\vert \mu \vert = \varphi \cdot \mathscr{L}^{d} \mvert [0,1]^{d}$. It follows that $\mu_{j} \to \mu$ strictly, and so we can use the assumed statement \eqref{Reshetnyak} to deduce that
\begin{align*}
    \int_{[0,1]^{d}} f\bigg(x,\frac{\mathrm{d}\mu_{j}}{\mathrm{d}\vert \mu_{j} \vert}(x)\bigg) \dd\vert \mu_{j} \vert(x) \to \int_{[0,1]^{d}} f(x,(1-t)v_{1}+tv_{2}) \varphi(x)\dd x.
\end{align*}
On the other hand, it follows from \eqref{ReshetnyakSharp1} and \eqref{ReshetnyakSharp2} that
\begin{align*}
    \int_{[0,1]^{d}} f\bigg(x,\frac{\mathrm{d}\mu_{j}}{\mathrm{d}\vert \mu_{j} \vert}(x)\bigg) \dd\vert \mu_{j} \vert(x) \to \int_{[0,1]^{d}} ((1-t)f(x,v_{1}) + tf(x,v_{2})) \varphi(x)\dd x.
\end{align*}
Since the choice of $\varphi$ was arbitrary and $f$ is continuous, we have
\begin{align*}
    f(x,(1-t)v_{1}+tv_{2}) &= (1-t)f(x,v_{1}) + t f(x,v_{2}) \hspace{1em} \text{for $x \in [0,1]^{d}$}.
\end{align*}
The result follows.
\end{proof}

\begin{remark}\label{rm:Reshetnyak} Notice that, in the case $(X,\vert \cdot \vert) = (\subext{k}\R^{d},\Vert \cdot \Vert)$, the sequence $(\mu_{j})_{j}$ constructed in Proposition \ref{pp:ReshetnyakSharp} and its limit $\mu$ can be identified with normal $k$-currents in $\R^{d}$. 
\end{remark}

\bibliographystyle{plain}

\end{document}